\documentclass[11pt]{amsart}
\usepackage{fullpage}
\usepackage{amsmath, amssymb, amsthm}
\usepackage{amsaddr}

\usepackage{graphicx}

\usepackage{color}
\usepackage[colorlinks=true,linkcolor=red,citecolor=blue,urlcolor=cyan]{hyperref}

\usepackage[ruled,vlined]{algorithm2e}

\newtheorem{thm}{Theorem}[section]
\theoremstyle{remark}
\newtheorem{rem}[thm]{Remark}

\newcommand{\ie}{{\it i.e.}}
\newcommand{\eg}{{\it e.g.}}

\usepackage[backend=bibtex,giveninits=true,sorting=nyt,style=alphabetic,natbib=true,maxcitenames=2,maxbibnames=10,url=false,doi=true,backref=false]{biblatex}
\renewbibmacro{in:}{\ifentrytype{article}{}{\printtext{\bibstring{in}\intitlepunct}}}

\begin{document}
\title{Dynamics and stationary configurations \\  of heterogeneous foams}

\author{Dong Wang, Andrej Cherkaev, and Braxton Osting}
\address{Department of Mathematics, University of Utah, Salt Lake City, UT}
\email{ \{dwang,cherk,osting\}@math.utah.edu}
\thanks{B. Osting acknowledges partial support from NSF DMS 16-19755 and 17-52202.}

\subjclass[2010]{ 49Q05, 
58E12, 
53A10, 
35K15. 
}

\keywords{Minimal surface; foam; bubble; threshold dynamics method}

\date{\today}

\begin{abstract} 
We consider the variational foam model, where the goal is to minimize the total surface area of a collection of bubbles subject to the constraint that the volume of each bubble is prescribed. We apply sharp interface methods to develop an efficient computational method for this problem. In addition to simulating time dynamics, we also report on stationary states of this flow for $\leq 21$ bubbles in two dimensions and $\leq 17$ bubbles in three dimensions. For small numbers of bubbles, we recover known analytical results, which we briefly discuss. In two dimensions, we also recover the previous numerical results of Cox et. al. (2003), computed using other methods. Particular attention is given to locally optimal foam configurations and heterogeneous foams, where the volumes of the bubbles are not equal. Configurational transitions are reported for the quasi-stationary flow where the volume of one of the bubbles is varied and, for each volume, the stationary state is computed. The results from these numerical experiments are described and accompanied by many figures and videos.  
\end{abstract}

\maketitle

\section{Introduction}
We consider the model for a $d$-dimensional foam  ($d=2,3$) comprised of $n$ bubbles,  $\{\Omega_i \}_{i=1}^n$, each with a prescribed volume, $ \mathcal H^d(\Omega_i) = V_i$, that arrange themselves as to minimize the total surface area, 
\begin{equation}
\label{e:min}
\min_{ \mathcal H^d(\Omega_i) = V_i} \ \mathcal H^{d-1}( \cup_{i=1}^n \partial \Omega_i). 
\end{equation}
Here we have denoted the $d$-dimensional Hausdorff measure by $\mathcal H^d$. Note that in \eqref{e:min}, the interfaces between bubbles and the interface between the  bubbles, $ \cup_{i=1}^n \partial \Omega_i$ and the rest of Euclidean space, $\mathbb R^d \setminus \cup_{i=1}^n \partial \Omega_i$, receive equal weight. We refer to stationary solutions of \eqref{e:min} as \emph{stationary $n$-foams}. If the areas are all equal, we say the foam is \emph{equal-area} and otherwise we say the foam is \emph{heterogeneous}.  
The isoperimetric variational problem \eqref{e:min} is classical; its history and the state of known results can be found in the recent book \cite{Morgan_1988}. The two-dimensional problem is discussed in \cite{Foisy_1993,Cox_2003,Wichiramala_2004}, the three-dimensional problem is discussed in  \cite{Taylor_1976,Hutchings_2002}, and the higher-dimensional $n=2$ problem is discussed in \cite{Lawlor_2012}. We'll further review the most relevant of these results in Section~\ref{s:background}.

\medskip

In this paper, we apply sharp interface methods from computational geometry to investigate \eqref{e:min}; these methods are described in Section~\ref{s:Meth}. 
In particular, we study an approximate gradient flow of \eqref{e:min} in dimensions $d=2,3$, which gives the time-evolution of a foam for a given initial configuration. This corresponds to a volume-constrained mean curvature flow of the interfaces between bubbles. 
An example of such time-dynamics for an equal-area, two-dimensional, $n=12$-foam is given in Figure~\ref{fig:dynamics}. 
An example for an equal-area, three-dimensional $n=8$-foam is given in Figure~\ref{fig:dynamics3d}.

We also study stationary foams of the gradient flow. In two dimensions, we recover many of the results from \cite{Cox_2003}, where candidate solutions for the equal-area problem \eqref{e:min} for many values of $n$ were found using very different computational methods then the present work. See Figure~\ref{f:eq-area} for stationary configurations of two-dimensional, equal-area $n$-foams for $n=2,\ldots,21$.
 Particular emphasis is given to the existence of multiple stationary foams that correspond to geometrically distinct configurations but have similar total surface areas. For example, a second two-dimensional, $n=16$-foam with slightly larger total perimeter than the configuration in Figure~\ref{f:eq-area} is given in Figure~\ref{f:162d}. Our computational methods also extend to three dimensions; stationary foams for the equal-area problem for $n=1,\ldots,17$ are displayed in Figure~\ref{f:eq-volume}. As far as we know,  these results are new for $n\geq 5$.

To further study multiple stationary foams, we consider heterogeneous foams. In particular, we study the \emph{quasi-stationary flow} where the area of one of the bubbles is slowly varied and for each area, the stationary solution is computed. We observe \emph{configurational transitions} where there are sudden changes in the stationary foams in this quasi-stationary flow. Examples of this can be seen in Figure~\ref{fig:increasing}. A comparison of two different quasi-stationary flows between an $n=6$ and $n=7$ equal-area foam is given in Figure~\ref{fig:hysteresis}. 

We conclude in Section~\ref{s:disc} with a discussion. 

\section{Background} \label{s:background}
In this section, we review some relevant previous results in two and three dimensions. 

\subsection{Two dimensional results}\label{s:b2d}
In 1993, Foisy, Alfaro, Brock, Hodges, and Zimba proved that the equal-area 2-foam in two dimensions is given by two intersecting discs separated by a line so that all angles are $120^\circ$ \cite{Foisy_1993}. 
In 2004, Wichiramala showed that the equal-area 3-foam in two dimensions is given by three intersecting discs so that all angles are $120^\circ$ \cite{Wichiramala_2004}. 
For a two-dimensional $n$-foam with $n\geq 4$, the optimal domain isn't known analytically, but, for small values of $n$, candidate solutions have been computed numerically \cite{Cox_2003}. 

For all $n$, necessary conditions for any minimizer are given by Plateau's laws:
\begin{itemize}
\item[(i)] each interface between bubbles has constant curvature and 
\item[(ii)] interfaces between bubbles meet in threes at vertices with equal angles. 
\end{itemize} 
We give a brief and formal derivation of Plateau's laws here; our goal is to give an accessible discussion that we can refer to when analyzing the numerical results.  

Given $n$ bubbles, $\Omega_1, \ldots, \Omega_n$, with given areas $V_1, \ldots, V_n$,  \ie, $\int_{\Omega_i} dx = V_i$, 
our variational problem is to find the configuration that has minimal total length of the interface $\Gamma$ that separates the bubbles. We assume that the bubbles are all contained in a region $\Omega$, and denote the complement of the bubbles in $\Omega$ by $\Omega_0 = \Omega \setminus \cup_{i=1}^n \Omega_i$. 
The interface $\Gamma$ is the union of the shared boundaries $\Gamma_{i,j}$ between all neighboring domains $\Omega_i$ and  $\Omega_j$ as well as the outer interfaces $\Gamma_{i0}$ of the external domains $\Omega_i$ with $\Omega_0$. The total length of the boundary is
$
J(\Omega_i) =\sum_{\substack{i, j=0 \\ i \neq j}}^n  \ \ \int_{\Gamma_{i,j}} ds. 
$
The variational problem can then be written 
\begin{subequations} \label{e:VarProb}
\begin{align}
\min_{\Omega_i}  & \ J(\Omega_i) \\
\textrm{subject to} & \  \int_{\Omega_i} dx = V_i, \quad i =1,\ldots, n.
\end{align}
\end{subequations}

\subsection*{Interfaces have constant curvature}
Introducing the Lagrange multipliers $\lambda_i$,  we formulate the Lagrangian for \eqref{e:VarProb}, 
\begin{align}
\label{ext-funct}
L(\Omega_i) &= J(\Omega_i) + \sum_{i=1}^n \lambda_i \left( \int_{\Omega_i} dx -V_i \right)
\end{align}

To see how the Lagrangian in \eqref{ext-funct} changes as we vary the domains $\Omega_i$, we first recall the formulas for the shape derivative of the area and perimeter with respect to changes in the domain.
Consider a domain $\Omega $ with piecewise smooth boundary $\Gamma$ and let  $s$ be the distance along the boundary.  
We consider the infinitesimal deformation of the domain in the direction of a velocity field $V$, which moves a point $x$ on the boundary $\Gamma$ to the point $x + \varepsilon \left( V(x) \cdot \hat n (x) \right) \hat n(x)$, where $\varepsilon> 0$ is small and $\hat n$ is the normal vector to $\Gamma$.  In other words, the point $x$ on the boundary of $\Gamma$ is moving in the normal direction at speed $\varepsilon c(x)$ where $c(x) = V(x) \cdot \hat n (x) $ .  
The resulting change in the area of $\Omega$,  $\delta |\Omega|$, and the change in the arc length of $\Gamma$, $\delta |\Gamma| $, are  given by 
$$
\delta |\Omega| = \varepsilon \int_\Gamma c(x) \ ds + o(\varepsilon), 
\quad \textrm{and} \quad 
\delta  |\Gamma| = \varepsilon \int_\Gamma  \kappa(s) c(x)  \ ds + o(\varepsilon), 
$$
where $\kappa(s)$ denotes the curvature of $\Gamma$. 

Using these shape derivatives, and looking for critical points of the Lagrangian $L$ in \eqref{ext-funct} due to a variation of the boundary $\Gamma_{i,j}$ between $\Omega_i $  and $\Omega_j $, we arrive at the condition  
$$
\int_{\Gamma_{i,j}} (\lambda_i- \lambda_j + \kappa_{i,j}) c(s) \ ds =0, 
\qquad \qquad \forall i,j =1,\ldots, n
$$
where $\kappa_{i,j}$ is the curvature of the boundary between domains $i$ and $j$ and $c(s)$ is the speed of variation at the point $s$ on an interface $\Gamma_{i,j}$. 
Since this condition should hold for all $c(s)$,  we arrive at the optimality condition
\begin{equation}
\kappa_{i,j}= \lambda_i-\lambda_j= \mbox{constant}
\label{curv-const}
\end{equation}
The optimality condition \eqref{curv-const} implies that
(i) the outer interfaces $\Gamma_{i0}$ of the external domains $\Omega_i$ with $\Omega_0$ are arcs of circles, 
(ii) the shared boundaries $\Gamma_{i,j}$ between all neighboring domains $\Omega_i$ and  $\Omega_j$ are arcs of circles, and, in particular,  
(iii) the interfaces between congruent bubbles are straight lines. 
The value of the Lagrange multiplier, $\lambda_i$, depends on the size of the domain $\Omega_i$ as well as on its position in the foam. In particular, the interface between a larger and smaller bubble should ``bend towards'' the larger shape. We have that   
$\lambda_i \to  
\begin{cases} 
0 & |\Omega_i| \to  \infty \\
\infty & |\Omega_i| \to  0 
\end{cases}$.

\subsection*{Triple junctions have equal angles} 
 Finding optimal angles between the arcs of three domains that meet at a single point requires a separate variational argument, analogous to the Weierstrass test \cite{Young1969}. 
 Assume that three boundary arcs $\Gamma_1$,  $\Gamma_2$, and $\Gamma_3$ meet at a point $z$ and consider a ball $B_\varepsilon$ of radius $\varepsilon$ centered at $z$.  
We now fix the ball $B_\varepsilon$ and the three points $x_i = \Gamma_i \cap \partial B_\varepsilon$ for $i=1,2,3$. 
We will minimize the Lagrangian, $L$, in \eqref{ext-funct} by varying the position of $z\in B_\varepsilon$. 
The change of the areas of the domains $\Omega_i$ is $O(\varepsilon^2)$ while the variation of the boundary lengths is $O(\varepsilon)$; therefore the contribution of the  increment of areas within $B_\varepsilon$ can be neglected. 
Next, the variation of the interface lengths are approximated (up to $o(\varepsilon)$) by the variation of distances $|x_i-z|$. 
We  arrive at the local problem:
$$
\min_z \ j(z),  \qquad \textrm{where} \quad  j(z) = \sum_{i=1}^3 |x_i-z|.  
$$

First, we observe that sum of any two angles between $\Gamma_1$, $\Gamma_2$ and $\Gamma_3$ is smaller than $180^\circ$. If an angle is larger than  $180^\circ$, than all  three circumferential points $x_1$ and $x_2$ and $x_3$ lie on one side of the ball in a half-disc. Such a configuration cannot be optimal because all three lengths can be decreased by simply shifting the point $z$ towards the middle point, $x_2$. 

If the angles are such that any two of them are smaller that $180^\circ$, the optimal intersection point $z$ is in the ball $B$ and may be found from the condition
$$
\nabla j(z)= 0  
\qquad \implies \qquad 
\sum_{i=1}^3 \frac{x_i-z}{|x_i-z|} = 0. 
$$
That is, the sum of the three unit vectors is zero, which implies that the angle between any two of them is $120^\circ$. One can also show that in a stationary foam, four or more bubbles cannot meet at a single point. 

\begin{rem}
The honeycomb structure satisfies the necessary conditions for optimality and is the optimal configuration of equal-area bubbles as $n\to \infty$ \cite{Hales_2001}.  
\end{rem}

\subsection{Three dimensional results} \label{s:b3d}
In three dimensions, less is known about optimal foam configurations. 
The double bubble conjecture was proven in 2002 by M. Hutchings, F. Morgan, M. Ritore, and A. Ros \cite{Hutchings_2002}. 
The necessary conditions for any minimizer are referred to as Plateau's laws:
\begin{itemize}
\item[(i)] interfaces between bubbles have constant mean curvature, 
\item[(ii)] bubbles can meet in threes  at $120^\circ$ angles along smooth curves, called Plateau borders, and 
\item[(iii)] bubbles can meet in fours and the four corresponding Plateau borders meet pairwise at angles of $\cos^{-1}(-1/3) \approx 109^\circ$. 
\end{itemize}
In what follows, we give a brief and formal derivation of these conditions here; a rigorous proof was given by Taylor \cite{Taylor_1976}. 

As in the two-dimensional case, we consider $n$ bubbles $\Omega_1$, ...  $\Omega_n$, with given volumes $V_1, \ldots, V_n$, \ie, $\int_{\Omega_i} dx = V_i$.
Our goal is to find the configuration that has minimal total surface area of the interfaces,  $\Gamma=\cup \Gamma_{i,j}$, that separate the bubbles. 
Again, the interfaces consists of the shared components $\Gamma_{i,j}$ of two neighboring domains $\Omega_i$ and  $\Omega_j$ for $i,j =1,\ldots, n$ and the interfaces $\Gamma_{i,0}$ of an external bubble $\Omega_i$ with the complement, $\Omega_0$. Introducing a multiplier $\lambda_i$ for each volume constraint, the Lagrangian for this variational problem is given by 
\begin{equation}
\label{ext-funct3d}
L(\Omega_i) =  \sum_{\substack{i, j=0 \\ i \neq j}}^n  \ \ \int_{\Gamma_{i,j}} ds + \sum_{i=1}^n \lambda_i \left( \int_{\Omega_i} dx -V_i \right),
\end{equation}
where $ds$ is an element of the interface $\Gamma_{i,j}$.

\subsection*{Interfaces have constant mean curvature} 
Taking the shape derivative of the Lagrangian in \eqref{ext-funct3d} and looking for critical points, an similar argument to the one given for two dimensions yields the stationary conditions
 $$
 \kappa_{i,j} = \lambda_i-\lambda_j \quad \mbox{on } \Gamma_{i,j}.
 $$
Here $\kappa_{i,j}$ is the mean curvature of the interface of $\Gamma_{i,j}$ (compare with \eqref{curv-const}). This condition states that the mean curvature of each interface, $\Gamma_{i,j}$, is constant. 

\begin{rem}
Minimal surfaces are a special case of the problem under study. Here, the constraints on volumes are not imposed; therefore the minimal surface problem corresponds to $\lambda_i=0, ~\forall i$ and has the well-known optimality condition: $\kappa = 0$. 
\end{rem}
 
\subsection*{Three bubbles meeting along a curve} 
We consider three smooth boundaries $\partial \Omega_{i}$, $\partial \Omega_{j}$, and $\partial \Omega_{k}$  intersecting along a curve $\gamma$, referred to as a \emph{Plateau border}. 
The conditions of optimality at the Plateau border $\gamma$  can  be found from  local variations inside an infinitesimal cylinder around the curve. 
The variation in an infinitesimal cylinder results in a change in the surface area that dominates the change in volume.  
Therefore, the necessary condition is identical to the corresponding well-studied condition for the minimal surface problem.
At any point of the Plateau border $\gamma$, the sum of the three normal vectors $\hat n_{i}$ to the intersecting surfaces $\partial \Omega_{i}$ is zero and these vectors are orthogonal to the tangent $\hat t$ of  $\gamma$:
$$
\sum_{i=1}^3 \hat n_i=0, \qquad \hat n_{i} \cdot \hat t=0, \quad i=1,2,3.
$$
 This implies  that   $ \hat n_i \cdot \hat n_j= - \frac{1}{2},~ ( i \neq j) $ and the angle between the normals is $120^\circ$.  

\subsection*{Four bubbles meeting at a point}
Similarly, we can consider a vertex where four bubbles intersect. Again, taking variations inside an infinitesimal ball around the vertex, we find that the sum of the four tangential vectors  $\hat t_i$ to the Plateau borders is zero:
$\sum_{i=1}^4 \hat t_i=0.$
This condition implies that the tangential vectors are the directions from the center of a regular tetrahedron to its vertices. Thus, 
  $\hat t_i \cdot \hat t_j= - \frac{1}{3},~ (i \neq j) $
and  the angle $\phi$ between any two tangent vectors  is
$\phi= \arccos(- \frac{1}{3}) \approx 109^\circ$. 

\begin{rem}
Kelvin's packing of truncated octahedra satisfy the necessary conditions for optimality  \cite{thompson1887}. 
The Weaire--Phelan structure also satisfies the necessary conditions for optimality and is the partition of three dimensional space with smallest known total surface area; it has 0.3\% smaller total surface area than Kelvin's structure \cite{WeairePhelan}. 
\end{rem}

\section{Computational Methods} \label{s:Meth}
In this section, we discuss  computational methods for the foam model problem \eqref{e:min}. 
Here, the goal is to find interfaces between adjacent bubbles such that the total interfacial area is minimal with the constraint that the volume of each bubble is fixed. 
To design a numerical algorithm for \eqref{e:min}, the first consideration is the method to represent the interfaces between bubbles. 
For contrast, we review several choices before describing the method used in the present work. 

\subsection{Previous Results} One method, known as the \emph{front tracking method}, uses a discrete set of points to represent the interfaces \cite{womble1989front}. Then, the energy is minimized by moving the points in the normal direction of the interface subject to some constraints. Although this idea is simple and straightforward, a number of difficult and complicated issues arise when dealing with multiple bubbles and possible topological changes, especially in three-dimensional simulations.  

In \cite{Brakke_1992}, the author developed and implemented\footnote{\url{http://facstaff.susqu.edu/brakke/evolver/evolver.html}} a method, referred to as the \emph{Surface Evolver}, for solving a class of problems, including \eqref{e:min}. A surface in this method is represented by the union of simplices and physical quantities (\eg, surface tension, crystalline integrands, and curvature) are computed using  finite elements. The surface evolver iteratively moves the vertices using the gradient descent method, thus changing the surface. 
Although this idea is simple and straightforward,  a number of difficult and complicated issues arise when dealing with multiple bubbles and possible topological changes. 

Another approach is the \emph{level set method}, where the interfaces is represented by the zero-level-set of a function $\varphi$ \cite{osher1988fronts}. This function is evolves in time according to a partial differential equation of Hamilton-Jacobi type, 
\[
\frac{\partial \varphi}{\partial t} = V_n|\nabla \varphi|. 
\]
Here, $|\cdot|$ is the Euclidean norm, $\nabla$ denotes the spatial gradient, and $V_n$ is the normal velocity. 
This type of method can easily handle topology changes because the interface is implicitly determined by the zero-level-set of the function $\varphi$. 
However, it is difficult to deal with the interface motion near multiple junctions and this type of method also needs to be reinitialized at each step or after every few steps. 

Another option is to use the \emph{phase field approach} where the interface is represented by a level-set of an order parameter function, $\phi$;  see, \eg, \cite{yue2004diffuse}.  Here, $\phi$ takes 
two distinct values (\eg, $\pm1$) for the two-phase case or 
several distinct vectors in the multiple-phase case. 
The function $\phi$ then evolves according to the Cahn-Hillard or Allen-Cahn equation, where a  
potential enforce that the function $\phi$ smoothly changes between the distinct values (or vectors) in a thin $\varepsilon$-neighborhood of the interface.
This approach is simple and insensitive to topological changes. 
However, if the evolution of multiple junctions with arbitrary surface tensions needs to be resolved, it is difficult to find a suitable multi-well potential. 
Also, since  it is desirable for $\varepsilon$ to be small, a very find mesh is needed to resolve the interfacial layer of width $\varepsilon$. 
Consequently, this algorithm is computationally expensive.

In \cite{Cox_2003}, the authors iterated a shuffling-and-relaxation procedure to gradually find a candidate foam. 
At each iteration, they selected the shortest side and applied to it a neighbor-swapping topological process followed by relaxing the configuration
in a quadratic mode. 
In the two-dimensional case, many nice candidates for various $n$ are presented in \cite{Cox_2003}. 
However, this method requires a careful choice of both the initial configuration and the shuffling procedure is heuristic. 
The candidate configuration highly relied on the initial ``circular'' configuration. 
Also, it appears that this procedure needs a large number of iterations to reach a stationary candidate. 
It would be challenging to apply these ideas to the three-dimensional or heterogeneous foams.

\subsection{Computational method}
In this paper, we use computational methods that are based on the \emph{threshold dynamics methods} developed in \cite{merriman1992diffusion,MBO1993,merriman1994motion,esedoglu2015threshold}. 
Here, $n$ indicator functions are used to denote the respective regions of each bubble in an $n-$foam. 
Additionally, we fix a rectangular box, $\Omega \subset \mathbb R^d$ ($d= 2, 3$), which contains the supports of these $n$ indicator functions and add an $(n+1)$-th indicator function to denote the complement of the $n-$foam.
Let $u = (u_1,u_2,\cdots,u_{n+1})$ denote these indicator functions. We define 
\begin{align*}
\mathcal{B}=\left\{u\in BV(\Omega)\colon u_i(x)=\{0, 1\}, \ \sum_{i=1}^{n+1}u_i=1,\ a.e.\ x \in\Omega, \ \textrm { and } \int_{\Omega} u_i(x) \ dx=V_i, \ i\in [n+1] \right\}, 
\end{align*}
where $V_i$ is the prescribed volume of the $i-$th bubble for $i\in [n+1]$. 
The constraints that $u_i(x) \in \{0,1\}$ and $\sum_i u_i = 1$ together force the indicator functions to have disjoint support---which is equivalent to their representative domains being disjoint.  
We approximate the surface area of the interface between the $i$-th and $j$-th bubbles by $\mathcal H^{d-1}( \partial \Omega_i \cap \partial \Omega_j) \approx L(u_i, u_j)$, with
\begin{align}\label{lengthapprox}
L(u_i, u_j) := \frac{\sqrt{\pi}}{\sqrt{\tau}}\int_{\Omega} u_i(x) (G_{\tau}*u_j)(x) \  dx, 
\quad \textrm{ where } \ \ 
G_{\tau}(x)=\frac{1}{(4\pi \tau)^{\frac d 2 }}\exp\left(-\frac{|x|^2}{4 \tau} \right).
\end{align}
The $\Gamma$ convergence of \eqref{lengthapprox} to the interfacial area was proven in \cite{alberti1998non,miranda2007short,esedoglu2015threshold}. 
Using \eqref{lengthapprox}, the optimization problem \eqref{e:min} can be approximated as  
\begin{equation} \label{pro:approx}
\min_{u \in \mathcal B} \ \mathcal E^\tau(u), 
\qquad \textrm{ where } \quad 
\mathcal E^\tau(u) = \sum_{\substack{i,j = 0 \\ i \neq j}}^{n+1} L(u_i, u_j).
\end{equation}

Since the energy functional $\mathcal E^\tau(u)$ is concave, we can relax the constraint set in \eqref{pro:approx} to obtain the equivalent problem \cite{esedoglu2015threshold,OstingWang2017,OstingWang2018}, 
\begin{equation} 
\label{pro:approx}
\min_{u \in \mathcal K} \ \mathcal E^\tau(u),
\end{equation}
where
\[
\mathcal K = \left\{u\in BV(\Omega)\colon u_i(x)\in [0, 1],  \  \sum_i^{n+1}u_i=1,\ a.e.\ x \in\Omega,  \ \textrm{ and }  \ \int_{\Omega} u_i(x) \ dx=A_i, i\in [n+1] \right\}
\]
is the convex hull of $\mathcal{B}$.  The sequential linear programming approach to minimizing $\mathcal E^\tau(u)$ is to consider a sequence of functions $\{u^s:= (u_1^s, u_2^s, \cdots, u_{n+1}^s)\}_{s=0}^\infty$which satisfies
\begin{equation}\label{pro:lin}
u^{s+1} = \arg \min_{u\in \mathcal K} \ \mathcal L_{u^s}(u)
\end{equation}
where $\mathcal L_{u^s}^\tau (u)$ is the linearization of $\mathcal E^\tau$. In this case, 
\[\mathcal L_{u^s}^\tau (u) = \sum_{i=1}^{n+1} \int_{\Omega}  \Psi_i^s(x) u_i(x) \ dx, 
\qquad \textrm{ where } \quad 
 \Psi_i^s = \sum_{\substack{j=1 \\ j\neq i}}^{n+1} G_{\tau}*u_j^s= G_{\tau}*\left(1-u_i^s\right).
 \]
Since $u^s$ is given, \eqref{pro:lin} is a linear minimization problem. 
If we were to neglect the volume constraints, \eqref{pro:lin} could be solved point-wisely by setting
\begin{equation}
u_i^{s+1}(x) = \begin{cases} 1 & \textrm{if} \  \Psi_i^s(x)  =  \min_{k\in [n+1]}  \Psi_k^s(x);  \\
0 & \textrm{otherwise}.
\end{cases}
\end{equation}
However, this solutions generally doesn't satisfy the volume constraints. 

Motivated by the schemes for the volume-preserving, two-phase flow \cite{ruuth2003simple,xu2016efficient,Elsey_2017}, 
to find a solution $u^{s+1}\in \mathcal{B}$ (\ie, each $u_i^{s+1}$ satisfies the corresponding volume constraint), 
Jacobs et. al. proposed an efficient auction dynamics scheme to impose the volume constraints for the multiphase problem \cite{jacobsauction}. 
In particular, they developed a membership auction scheme to find $n+1$ constants $\lambda_i$, $i\in [n+1]$ such that the solution $u^{s+1} \in \mathcal B$ can be solved by 
\begin{equation}
u_i^{s+1}(x) = \begin{cases} 1, & \textrm{if} \  \Psi_i^s(x)+\lambda_i  = \min_{k\in [n+1]}  \ (\Psi_k^s(x)+\lambda_k)  \\
0, & \textrm{otherwise}.
\end{cases}
\end{equation}
The algorithm is summarized in Algorithm \ref{a:MBO} and we refer to  \cite{jacobsauction} for details of the derivation. The algorithm was also proven to be unconditionally stable for any $\tau>0$  \cite{jacobsauction}.

\begin{algorithm}[h]
\DontPrintSemicolon
 \KwIn{Let $\Omega_n$ be the discretization of the domain $\Omega$, 
 $n$ be the number of grid points, 
 $u^0 = (u^0_1,\ldots,u^0_{n+1})$ be the indicator functions for an initial $n+1$-partition, 
 $\tau>0$ be the time step, 
 $V_i$ for $i \in [n+1]$ be the prescribed volumes, 
 $\varepsilon_0$ be the initial value of $\varepsilon$, 
 $\alpha$ be the $\varepsilon$-scaling factor, and
  $\varepsilon_{\min}$ be the auction error tolerance.  }
 \KwOut{ $ u^S \in \mathcal B$ that minimizes  \eqref{pro:approx}.}
 
Set $s=1$ \;
Set $\bar{\varepsilon} = \varepsilon_{\min}/n$\;

 \While{not converged}{
{\bf 1. (Diffusion step)} Compute the  coefficient functions, 
$$
\Phi_i = 1-\Psi_i^s = G_{\tau}*u_i^s, \qquad \qquad i \in [n+1]
$$ 

{\bf 2. (Find $\lambda$ using auction dynamics)} \;
Set $\lambda_i = 0 $ for $i \in [n+1]$ \;
Set $\varepsilon = \varepsilon_0$ \;  

\While{$\varepsilon > \bar{\varepsilon}$} {
Mark all $x \in \Omega_n$ as unassigned\; 
Set $u_i = 0$ for $i \in [n+1]$ \; 
\While{some x is marked as unassigned}{\For{each unassigned $x \in \Omega_n$}{
Calculate $i^* \in  \arg \max_{i\in[n+1]} \Phi_i(x) - \lambda_i$ \;
Calculate $j^* \in \arg\max_{j \neq i^*} \Phi_j(x)-\lambda_j$ and set 
$$
b(x) = \lambda_{i^*}+\varepsilon+(\Phi_{i^*}(x)-\lambda_{i^*})-(\Phi_{j^*}(x)-\lambda_{j^*})
$$

\If{$\sum_{x} u_{i^*}(x) = V_{i^*}$}{
Find $y = \arg\min_{z\in u_{i^*}^{-1}(1)}b(z)$ \;
Set $u_{i^*}(y) = 0 $ and set $u_{i^*}(x) = 1$ \;
Mark $y$ as unassigned and mark $x$ as assigned \;
Set $\lambda_{i^*}= \min_{z\in u_{i^*}^{-1}(1)}b(z)$}
\Else{
Set $u_{i^*}(x) = 1$ \;
\If{$\sum_{x} u_{i^*}(x) = V_{i^*}$}{Set $\lambda_{i^*}= \min_{z\in u_{i^*}^{-1}(1)}b(z)$}}
}}
Set $\varepsilon = \varepsilon/ \alpha$ \; 
\If{$\varepsilon < \bar{\varepsilon}$}{Set $u^{s+1} = u$ }}

Set $s = s+1$} 
\caption{Auction dynamics algorithm for solving \eqref{pro:approx}  \cite[Algorithm 1, 2]{jacobsauction}.} 
\label{a:MBO} 
\end{algorithm}

\section{Two-dimensional numerical examples} \label{s:2d}

\subsection{Time-evolution of foams}\label{s:2dDynamics} 
\begin{figure}[t]
\centering
\includegraphics[width=1 \textwidth, clip, trim=5cm 0cm 5cm 2cm]{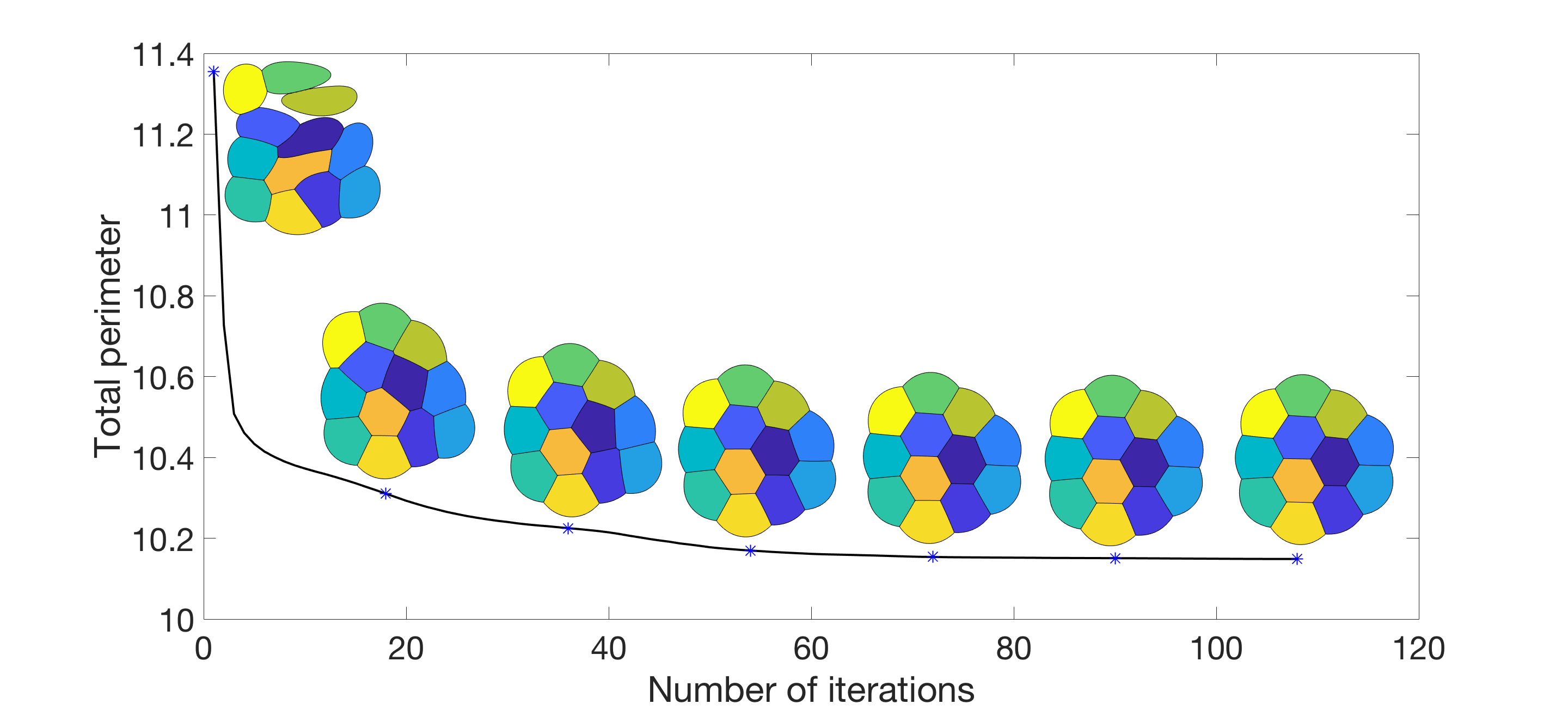}
\caption{A plot of the energy as a $n=12$-foam evolves from a random initialization together with the foam configuration at various iterations. See Section~\ref{s:2dDynamics}.}\label{fig:dynamics}
\end{figure}

For an equal-areal, $n=12$-foam, we show the time evolution corresponding to the gradient flow of the  total energy with a random initialization. In the subsequence, we generate the random initialization with volume constraints as the following: 
\begin{enumerate}
\item Generate a random $n$-Voronoi tessellation in a smaller box contained in the whole computational domain and set the complement as $n+1-$th Voronoi domain. 
\item Set $u_i = 1/\tilde{V}_i$, $i \in [N+1]$ where $\tilde{V}_i$ is the volume of the $i-$th Voronoi domain. 
\item Run Algorithm \ref{a:MBO} once to get a $n+1-$partition in the computational domain and set the corresponding indicator functions as the random initial condition.
\end{enumerate}
The energy at each iteration is plotted in Figure \ref{fig:dynamics} with the foam configuration at various iterations. Note that the energy decays very fast; in 108 iterations, the configuration is stationary in the sense that no grid points are changing bubble membership. After $\approx 50$ iterations, the foam configuration changes very little.

\subsection{Stationary solutions}  \label{s:2dStationary}
\begin{figure}[t!]
\begin{center}
\includegraphics[width=.22\textwidth,clip,trim=5cm 2cm 5cm 2cm]{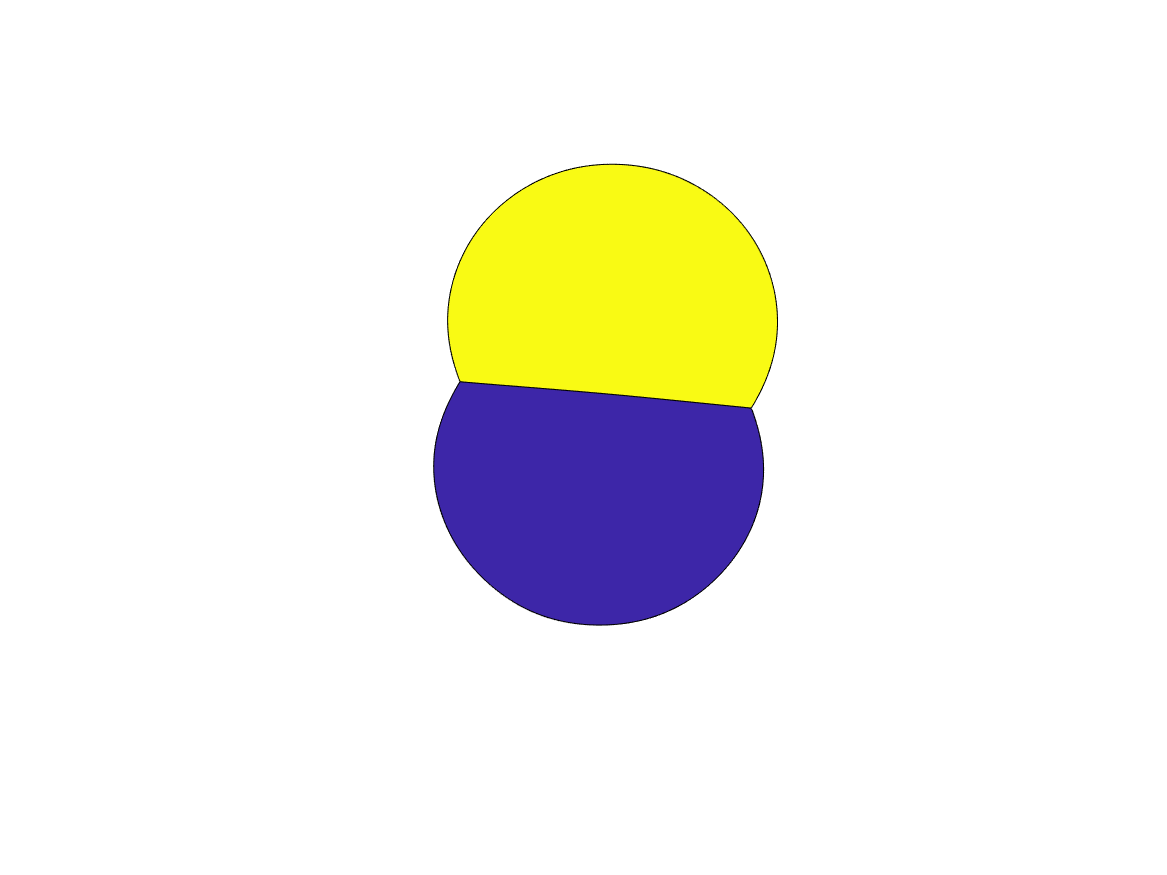}
\includegraphics[width=.22\textwidth,clip,trim= 5cm 2cm 5cm 2cm]{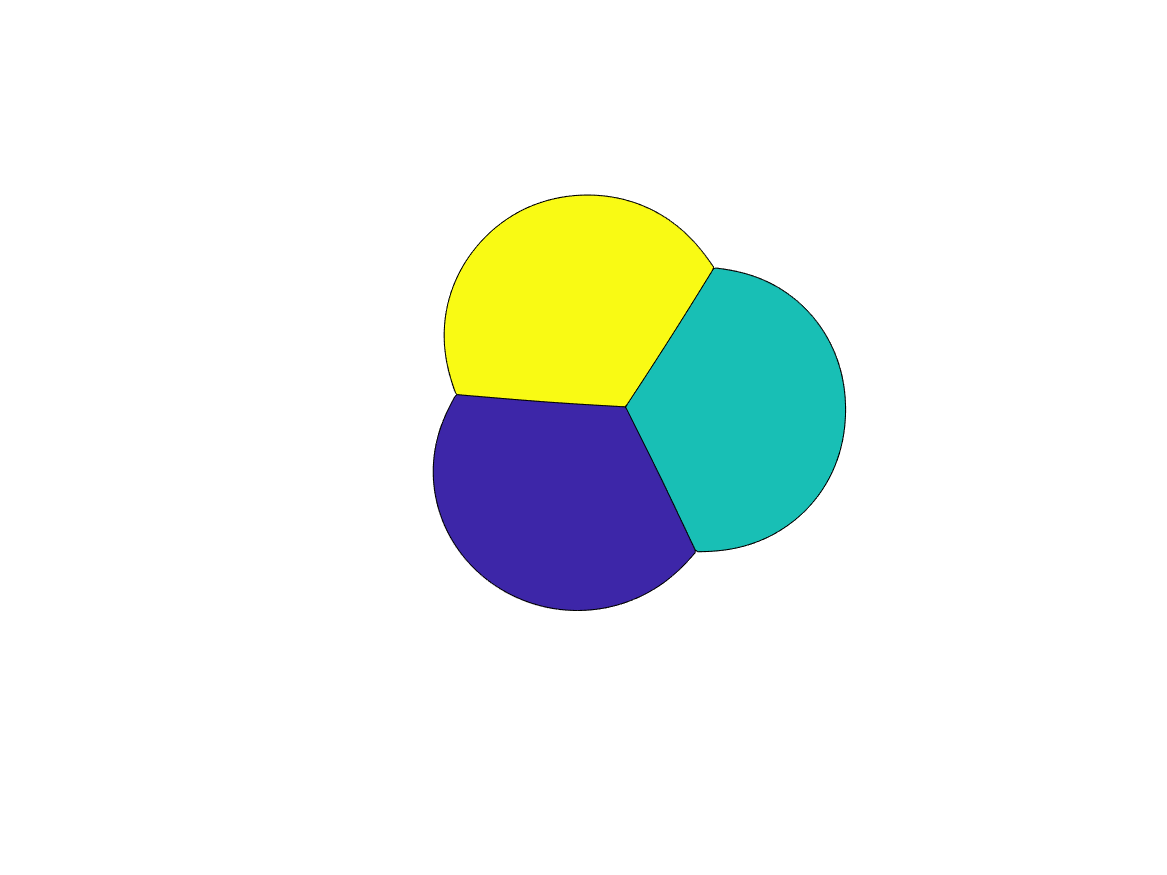}
\includegraphics[width=.22\textwidth,clip,trim= 5cm 2cm 5cm 2cm]{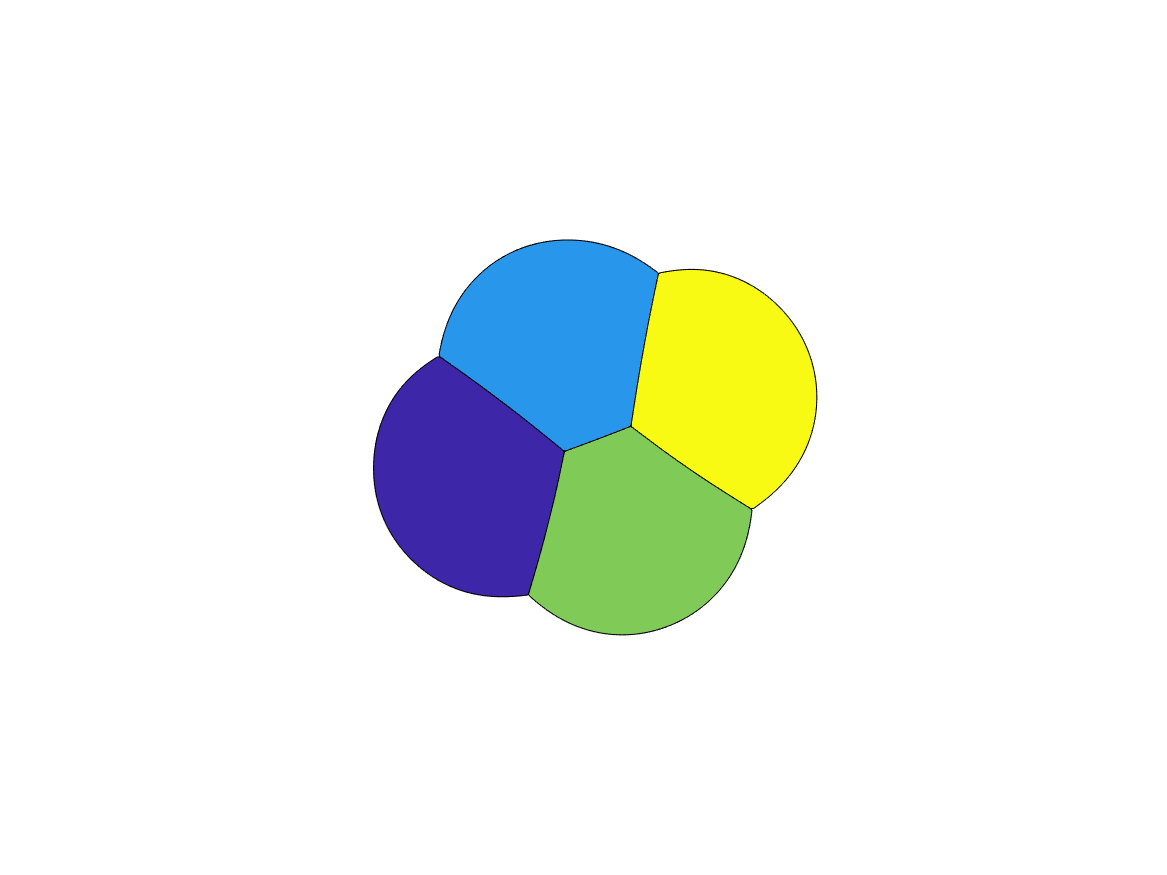} 
\includegraphics[width=.22\textwidth,clip,trim= 5cm 2cm 5cm 2cm]{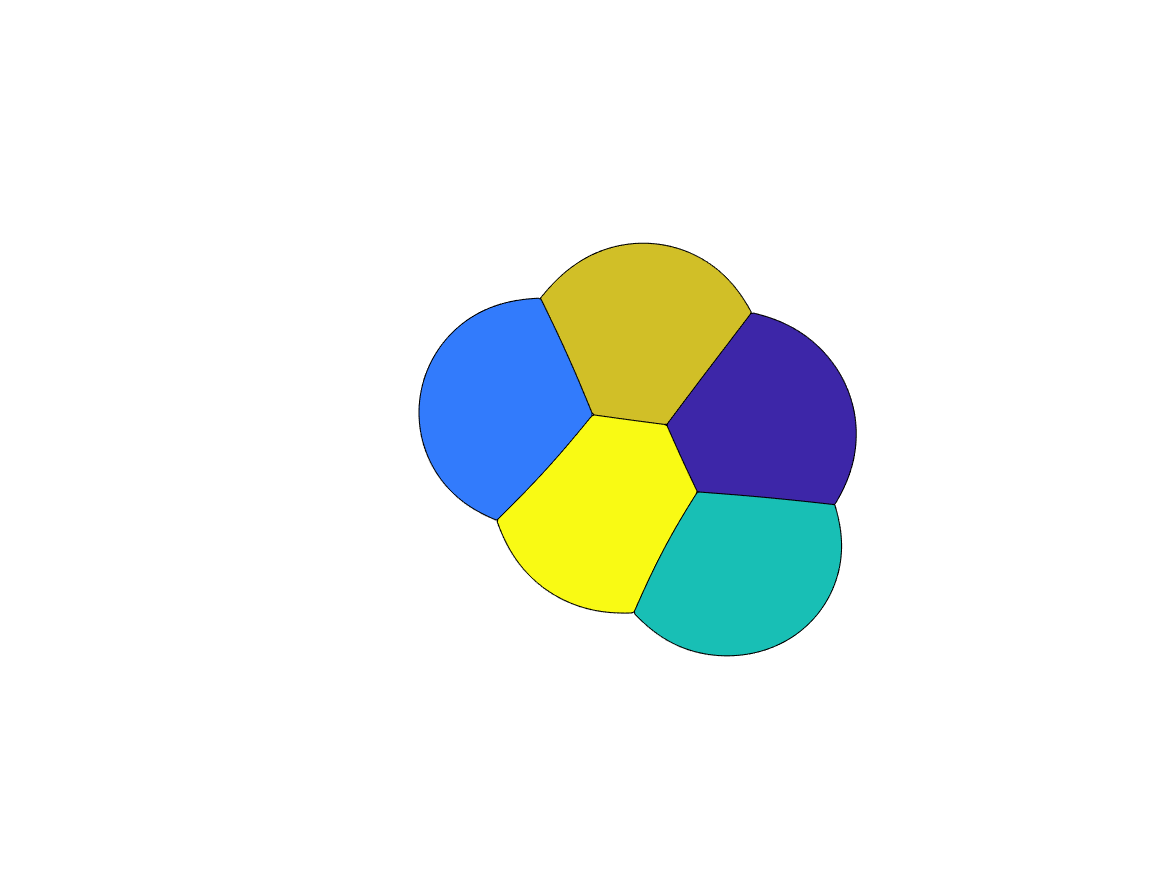} 
\includegraphics[width=.22\textwidth,clip,trim= 5cm 2cm 5cm 2cm]{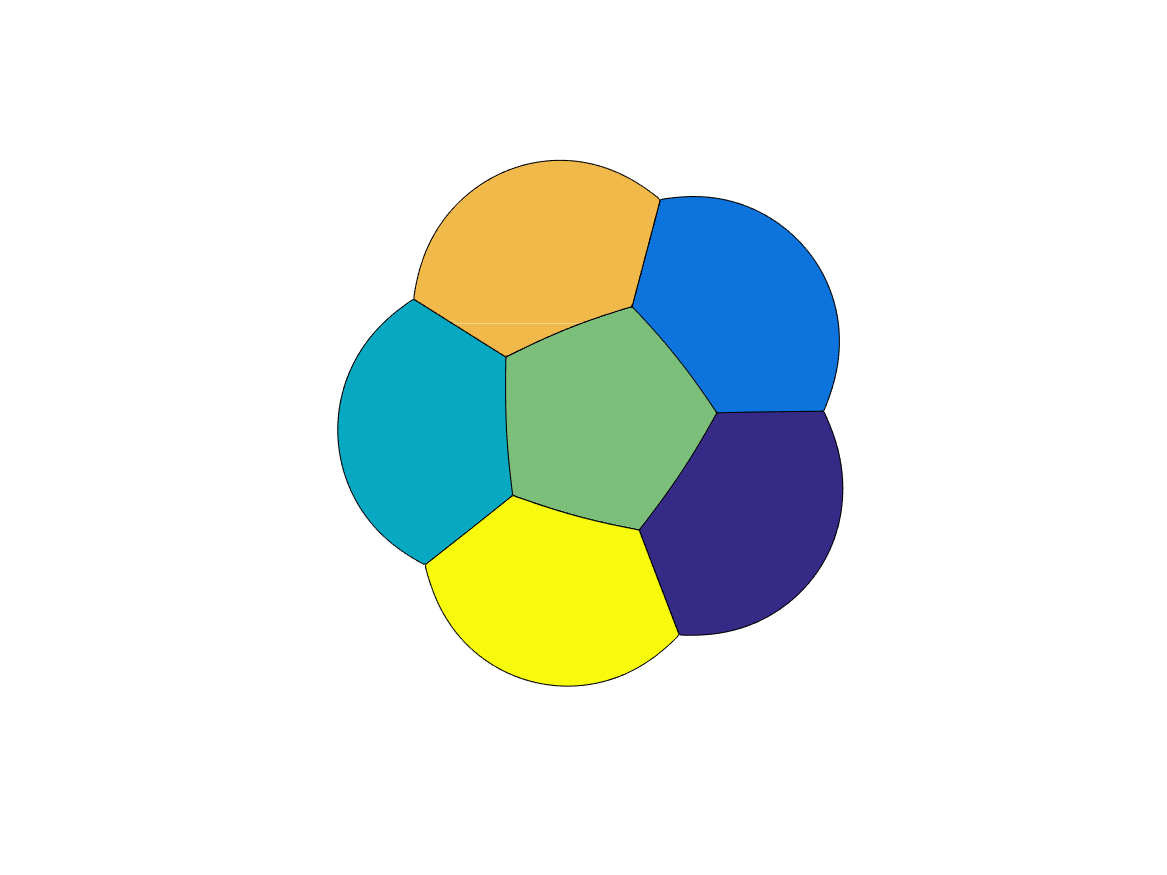}
\includegraphics[width=.22\textwidth,clip,trim=5cm 2cm 5cm 2cm]{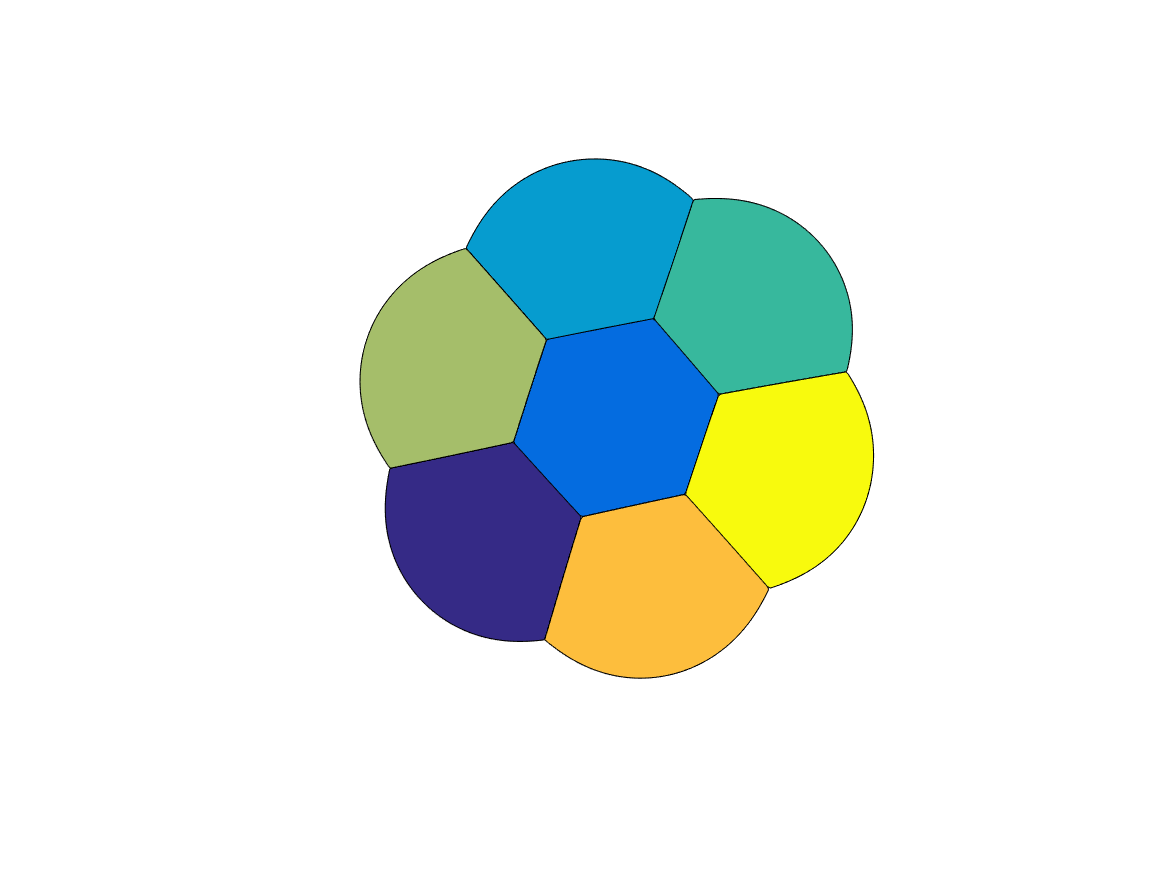}
\includegraphics[width=.22\textwidth,clip,trim= 5cm 2cm 5cm 2cm]{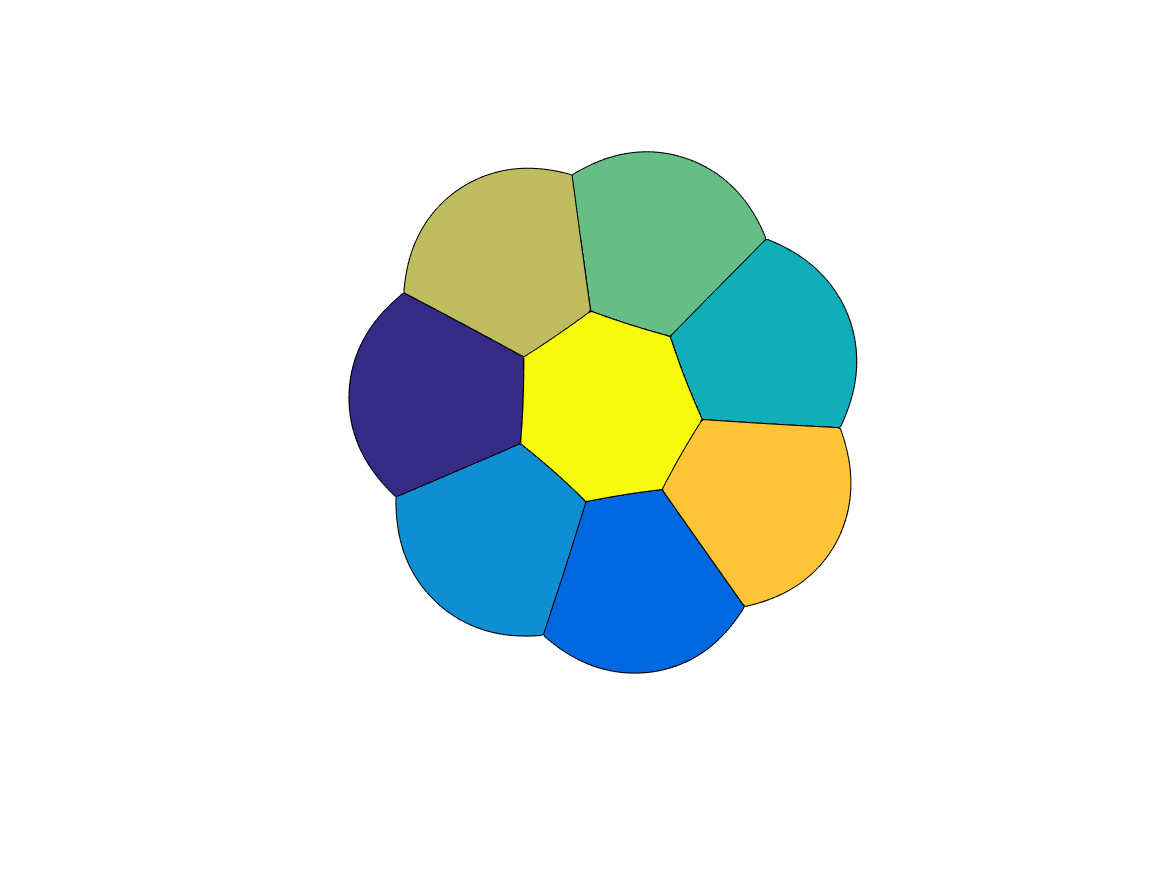}
\includegraphics[width=.22\textwidth,clip,trim= 5cm 2cm 5cm 2cm]{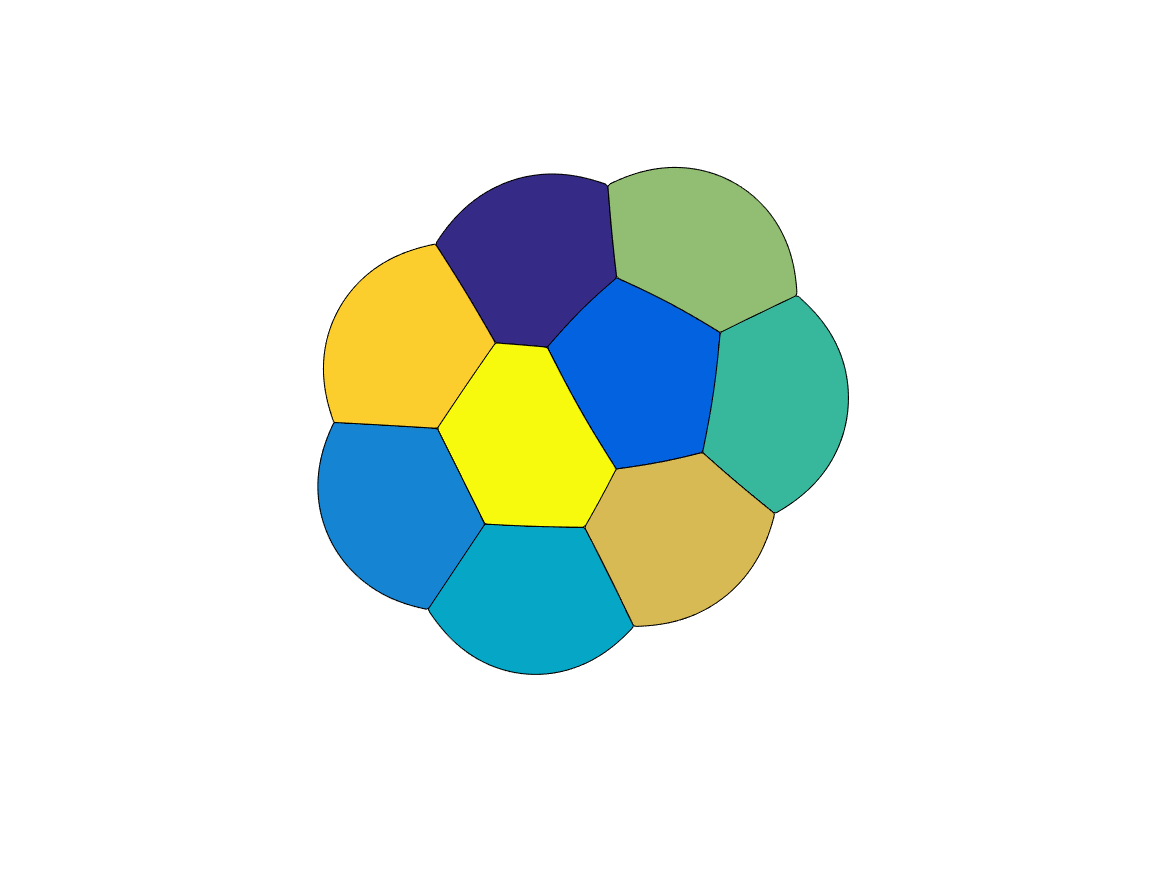}
\includegraphics[width=.22\textwidth,clip,trim= 5cm 2cm 5cm 2cm]{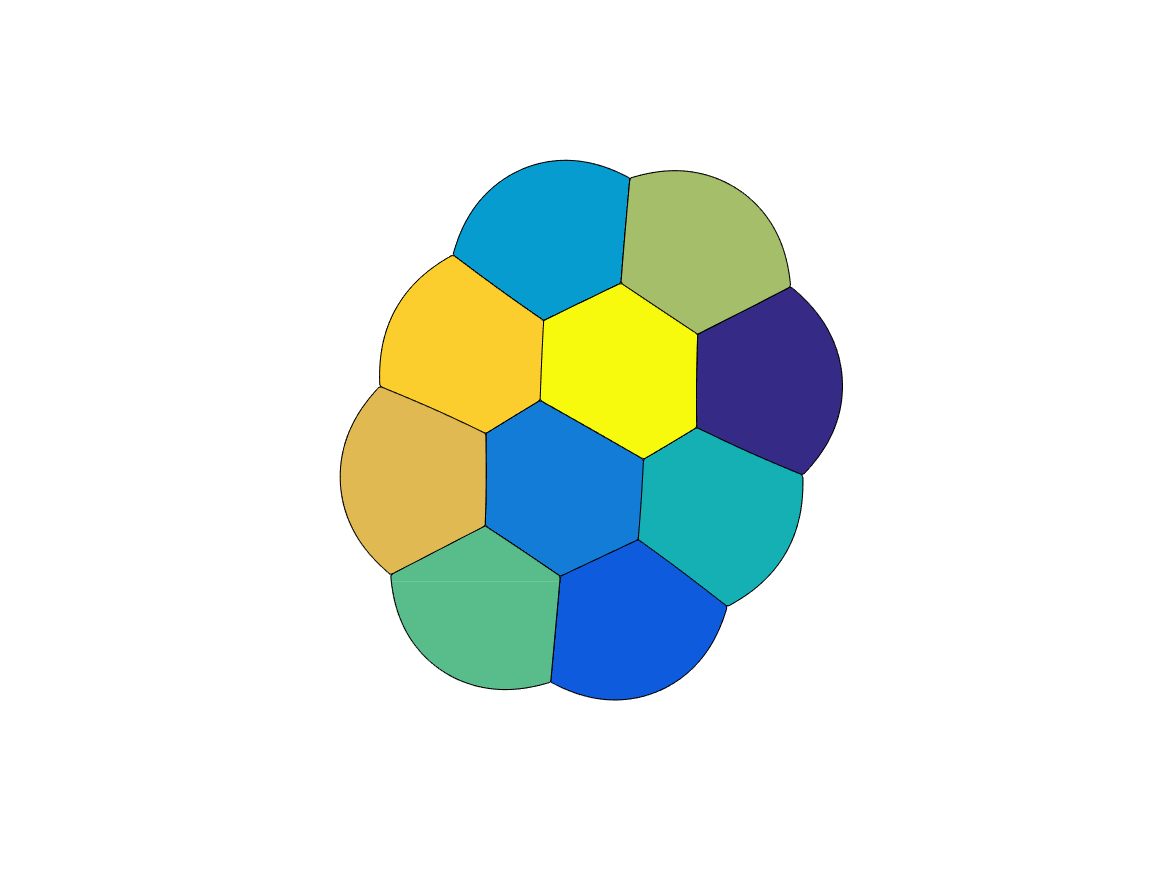}
\includegraphics[width=.22\textwidth,clip,trim= 5cm 2cm 5cm 2cm]{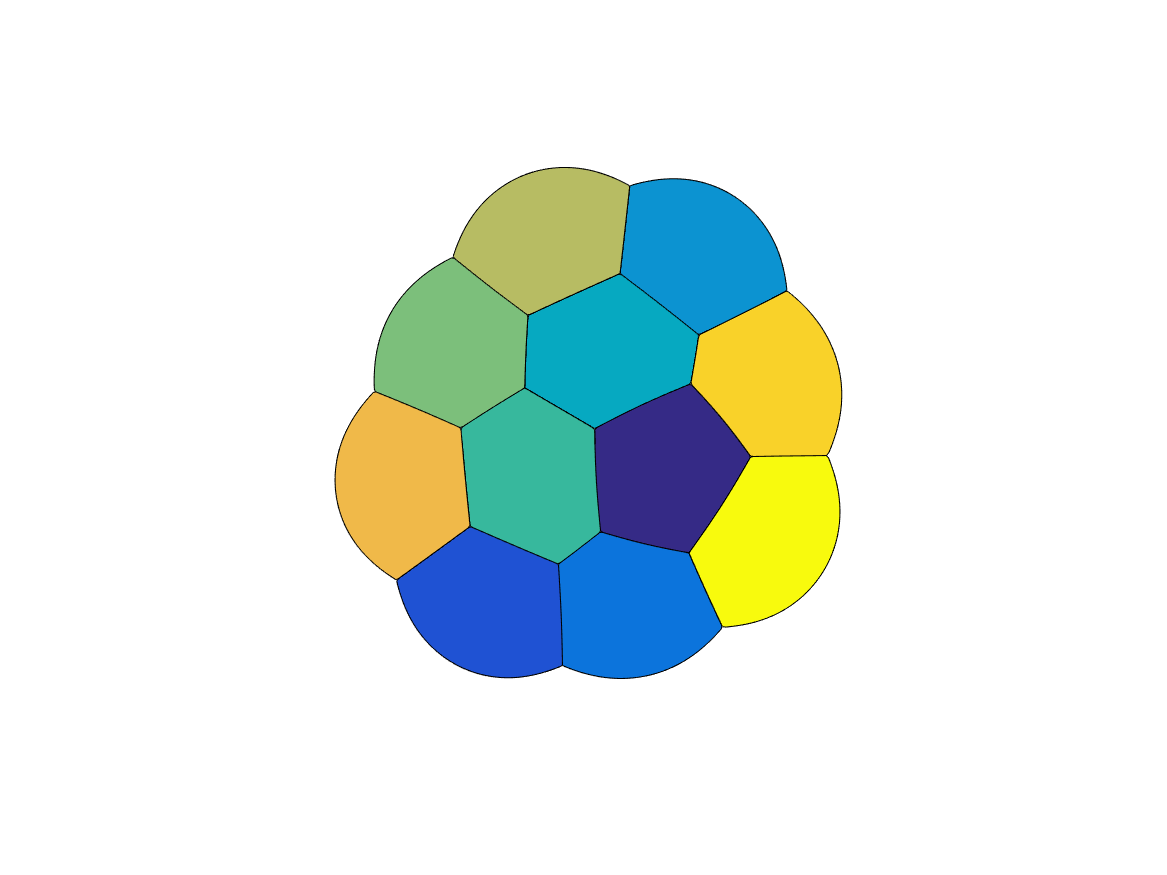}
\includegraphics[width=.22\textwidth,clip,trim= 5cm 2cm 5cm 2cm]{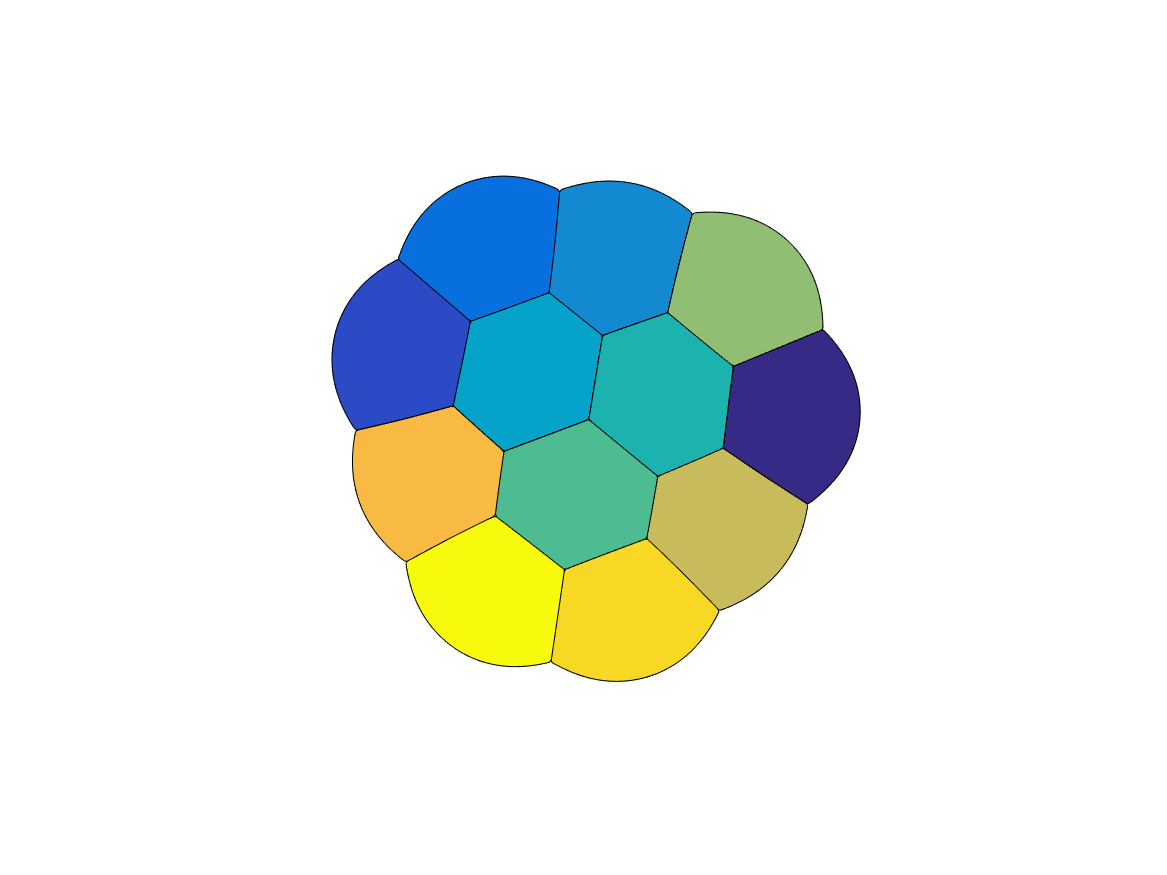}
\includegraphics[width=.22\textwidth,clip,trim= 5cm 2cm 5cm 2cm]{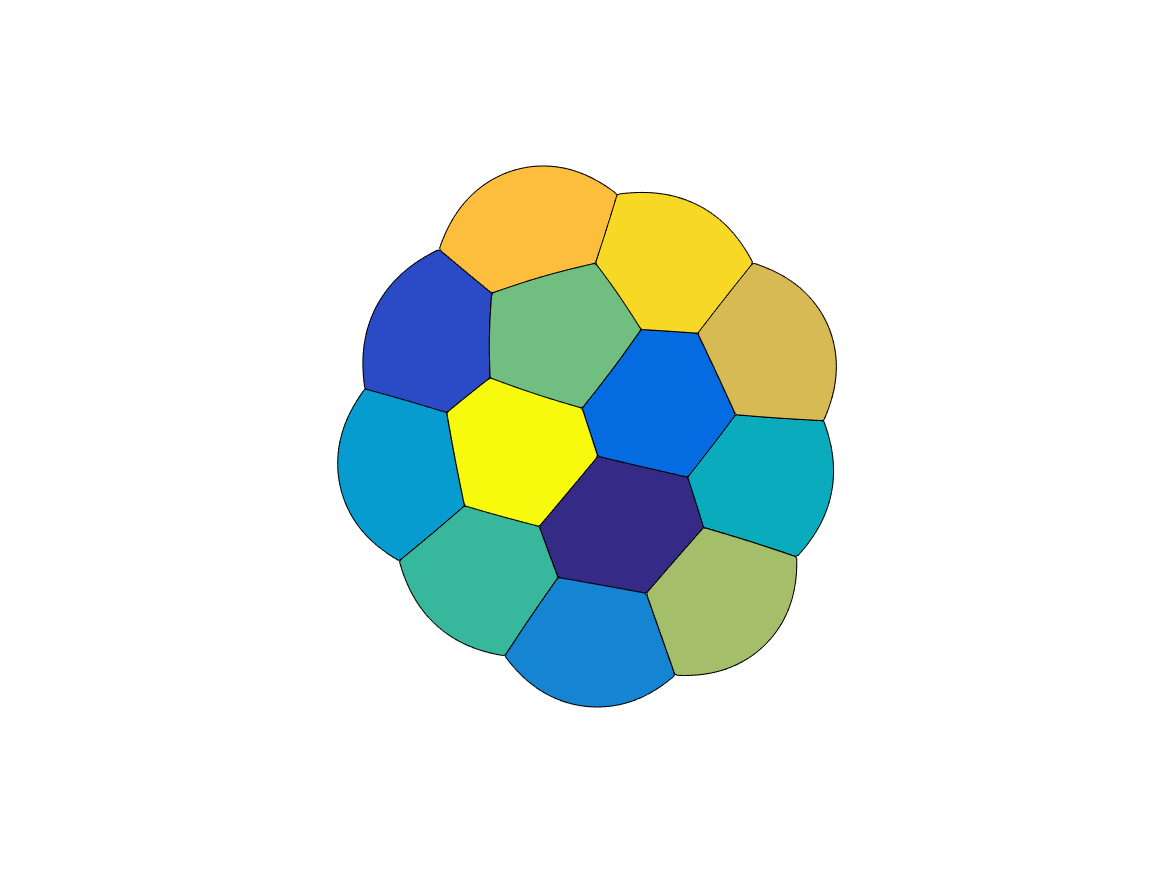}
\includegraphics[width=.22\textwidth,clip,trim= 5cm 2cm 5cm 2cm]{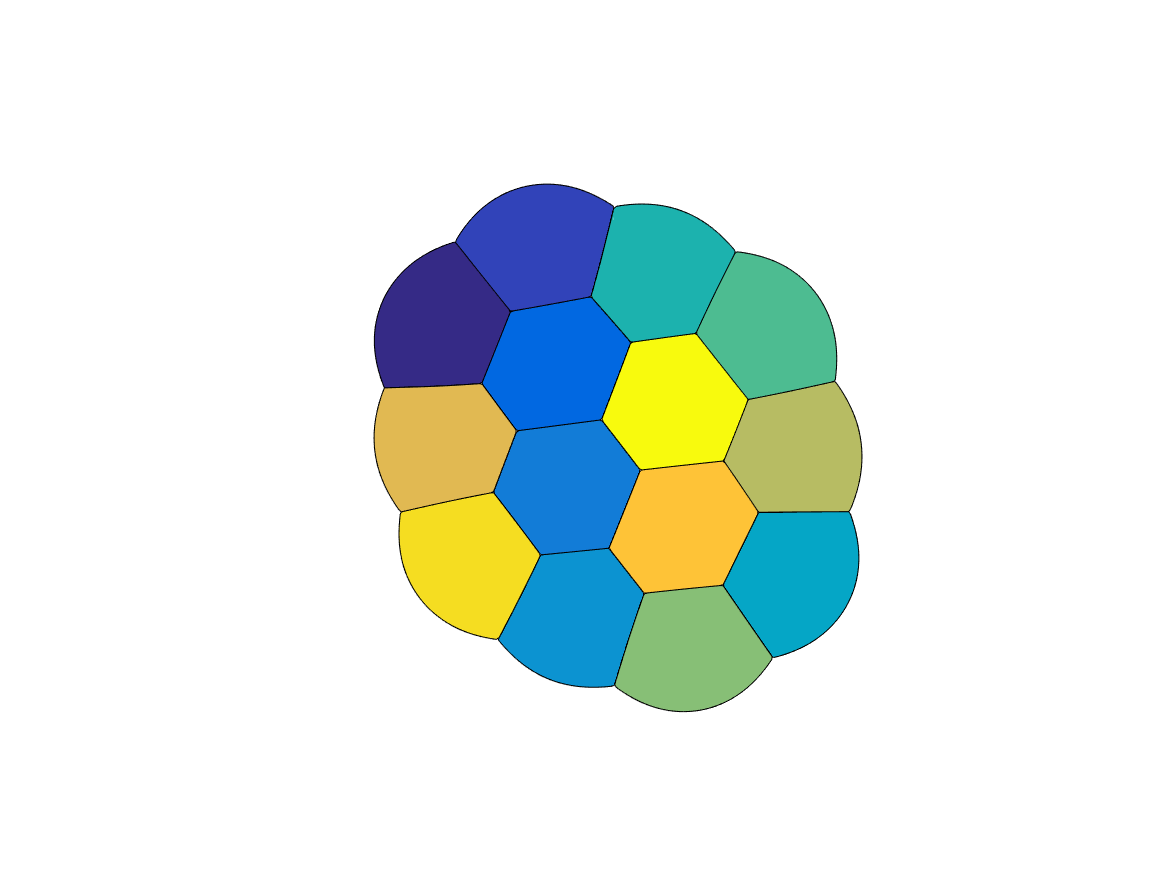}
\includegraphics[width=.22\textwidth,clip,trim= 5cm 2cm 5cm 2cm]{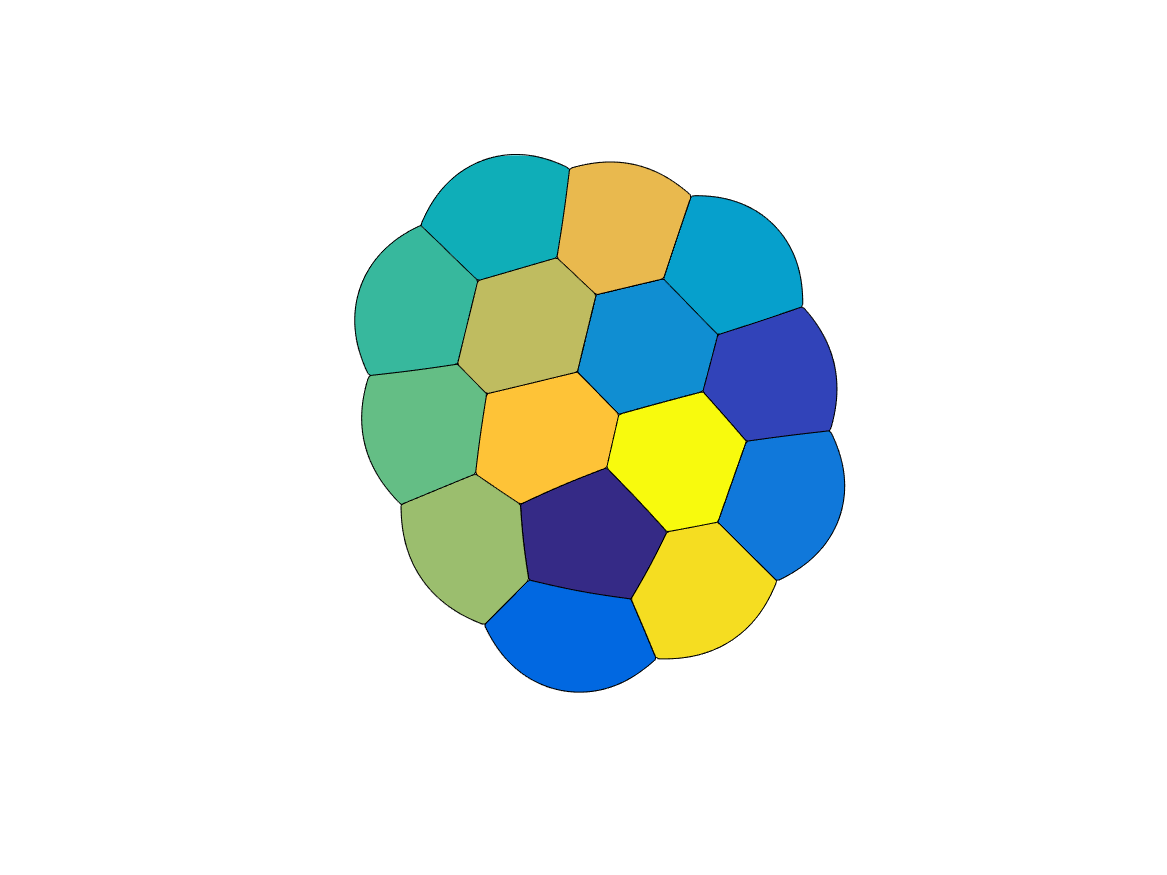}
\includegraphics[width=.22\textwidth,clip,trim= 5cm 2cm 5cm 2cm]{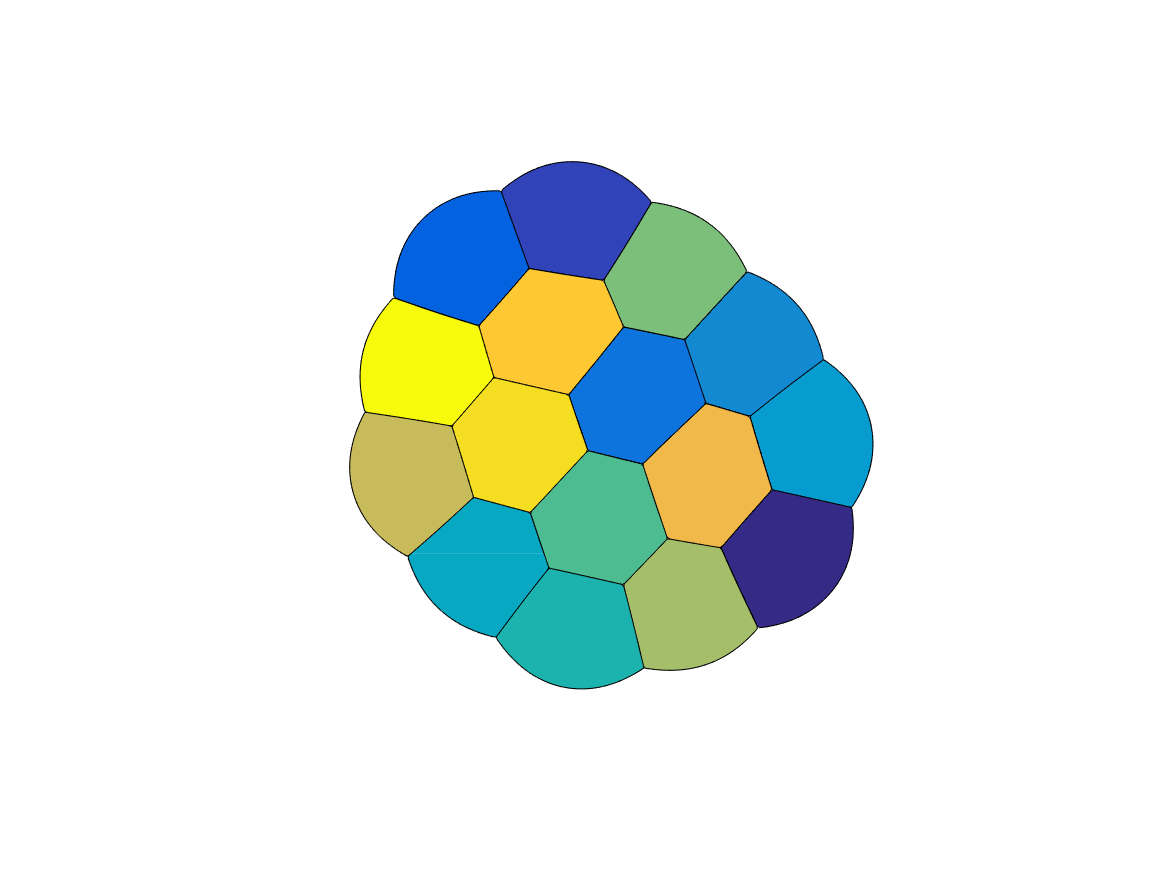}
\includegraphics[width=.22\textwidth,clip,trim= 5cm 2cm 5cm 2cm]{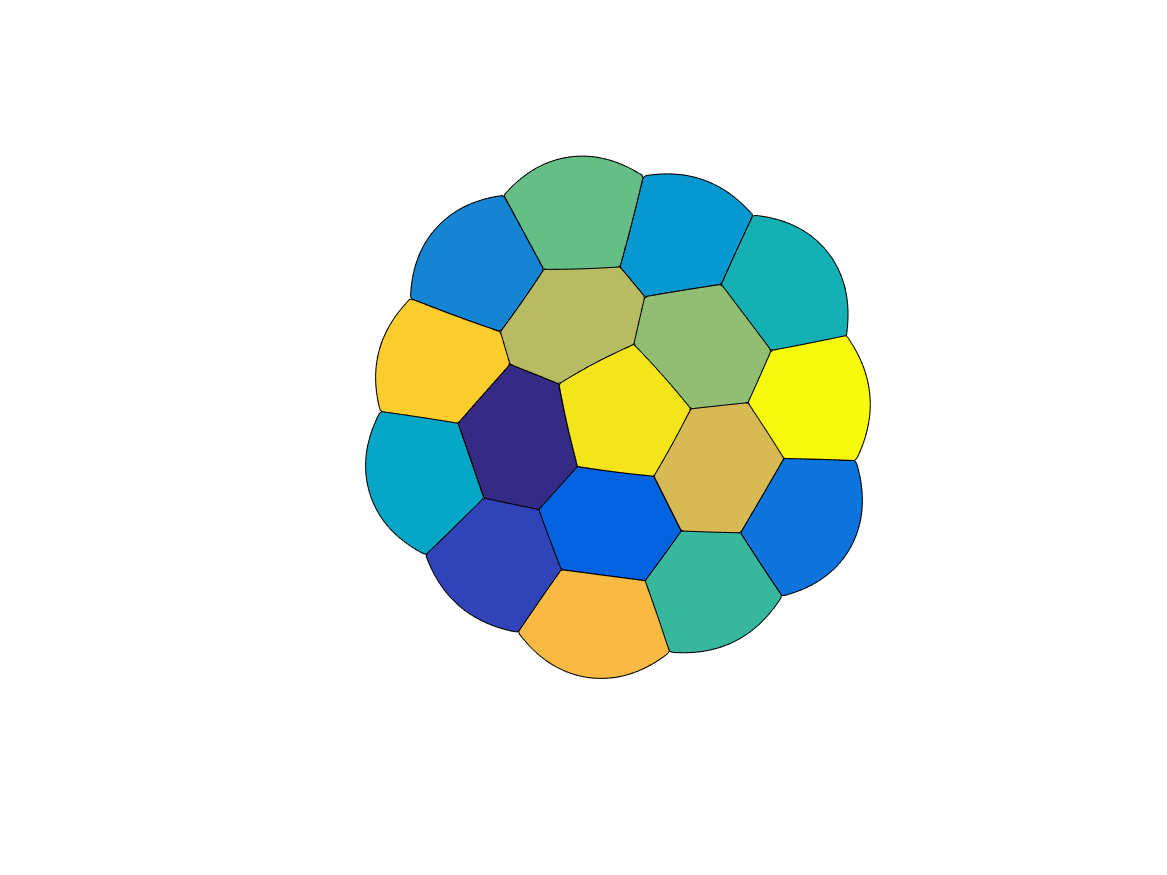}
\includegraphics[width=.22\textwidth,clip,trim= 5cm 2cm 5cm 2cm]{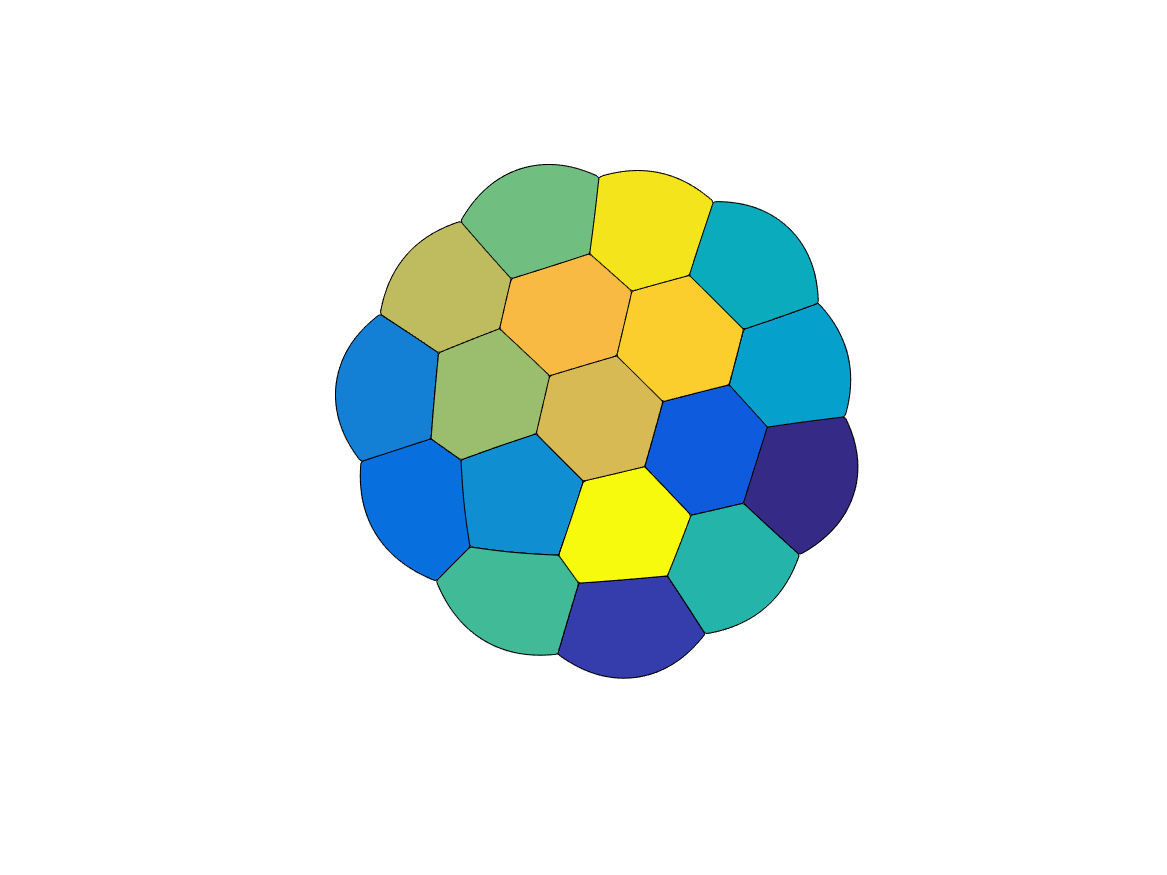}
\includegraphics[width=.22\textwidth,clip,trim= 5cm 2cm 5cm 2cm]{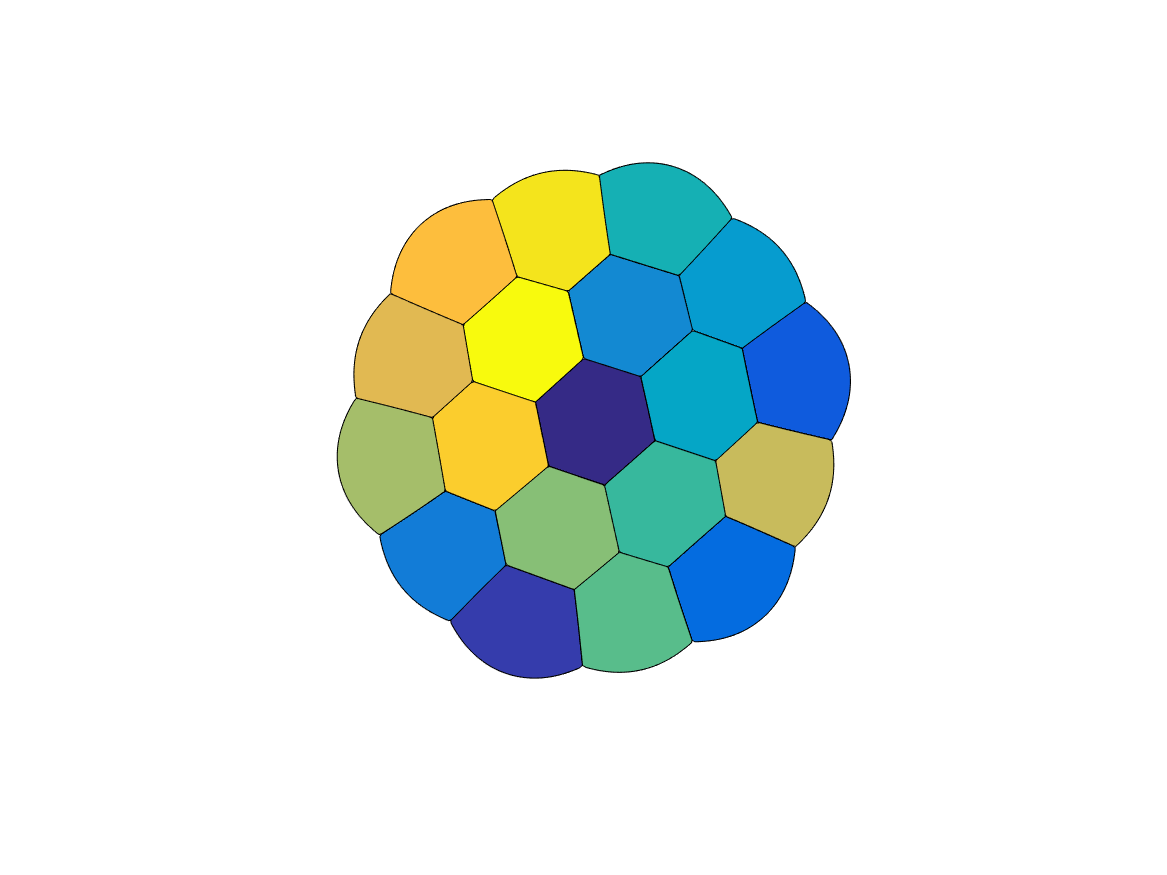}
\includegraphics[width=.22\textwidth,clip,trim= 5cm 2cm 5cm 2cm]{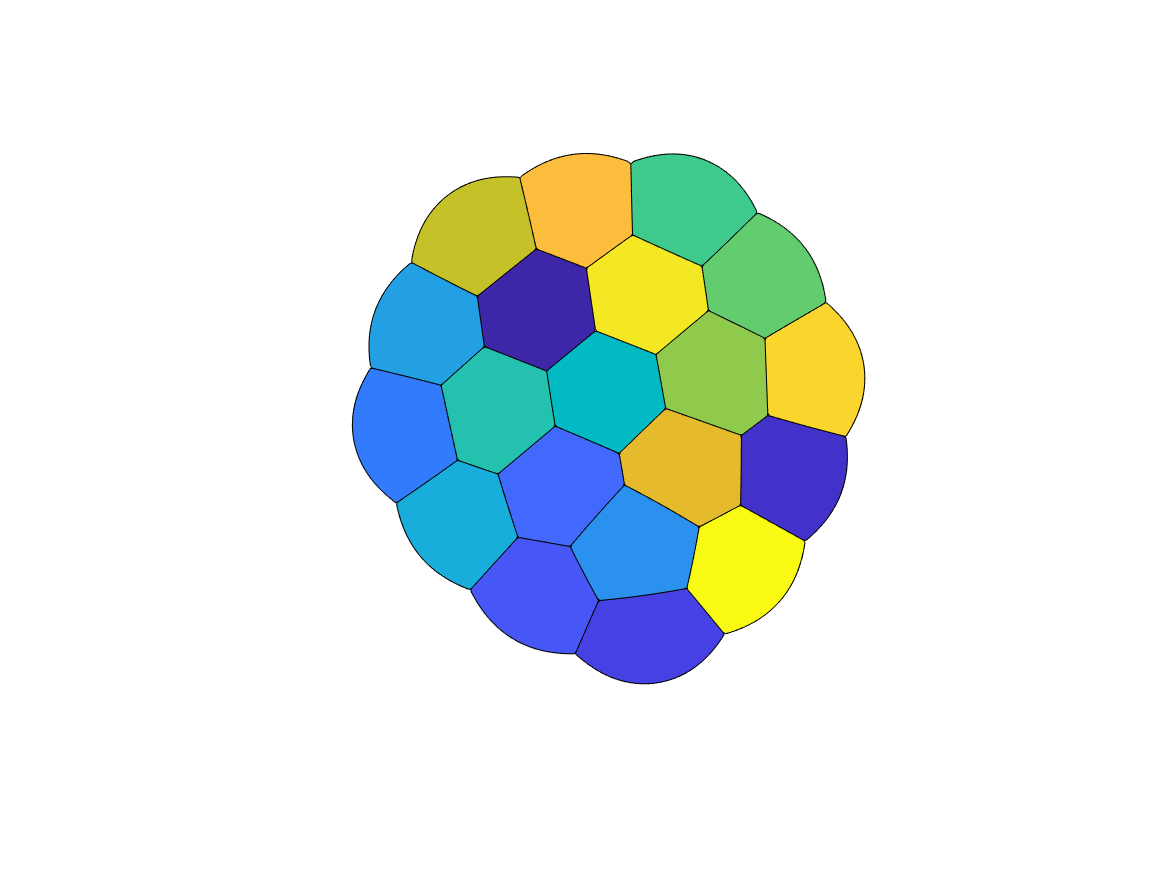}
\includegraphics[width=.22\textwidth,clip,trim= 5cm 2cm 5cm 2cm]{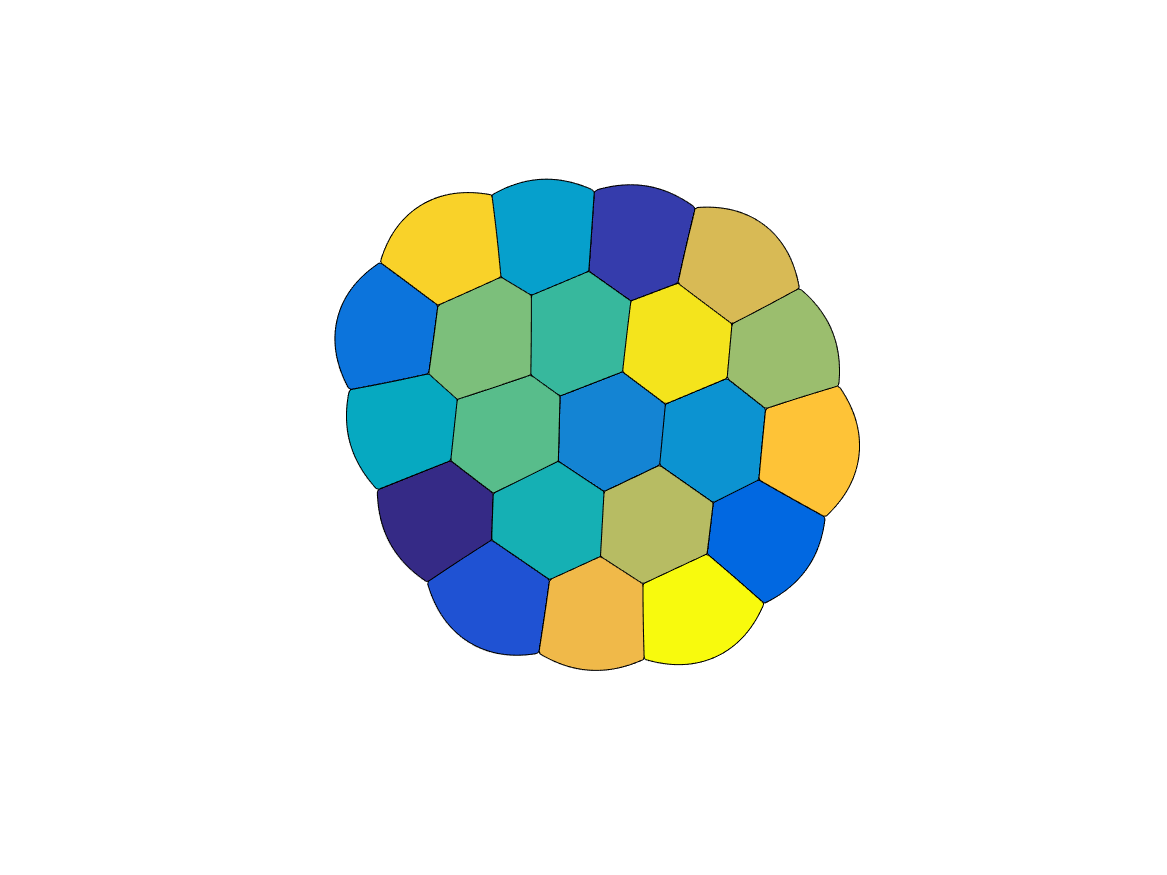}
\caption{Stationary equal-area $n$-foams for $n=2,\ldots  21$ with smallest computed total perimeter. See Section~\ref{s:2dStationary}.}
\label{f:eq-area}
\end{center}
\end{figure}

\begin{figure}[t]
\includegraphics[width=.2\textwidth,clip,trim= 5cm 3cm 4cm 2cm]{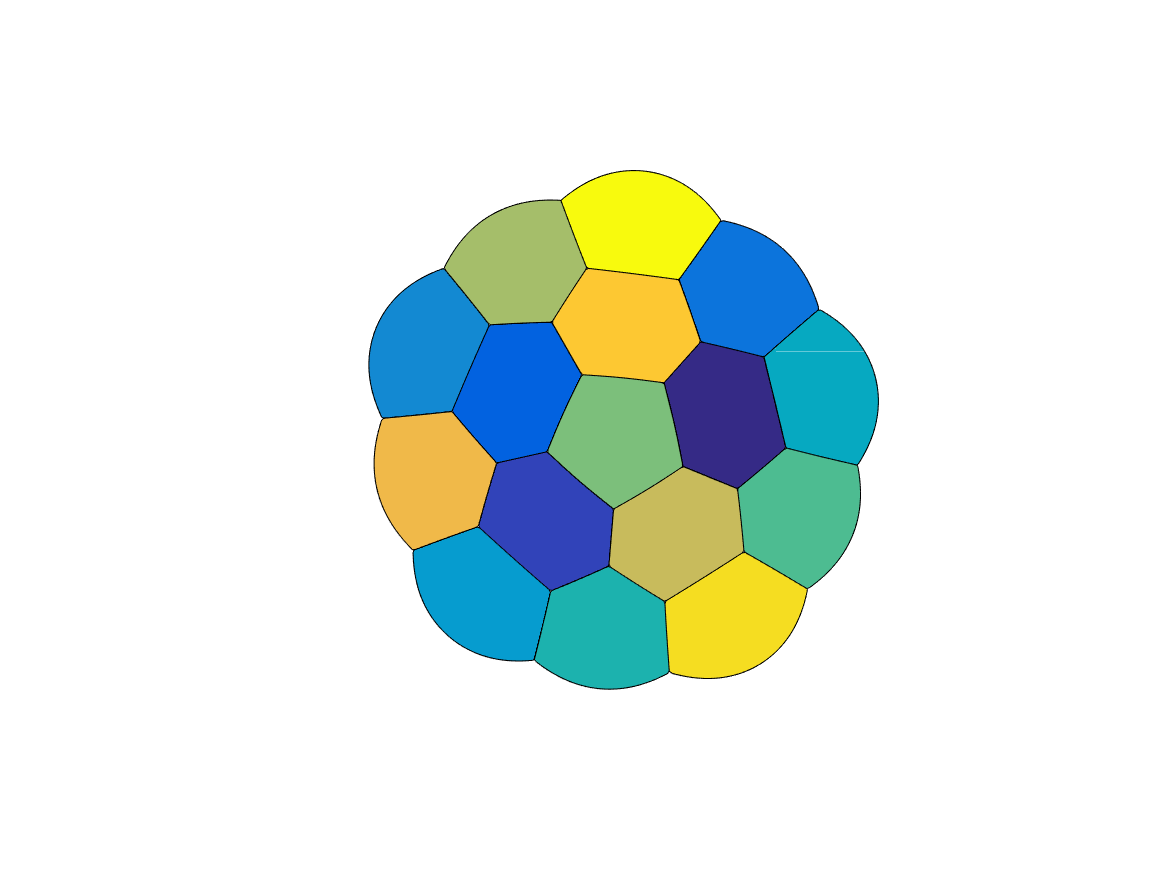}
\caption{Another stationary equal-area $16$-foam with larger total perimeter than the 16-foam displayed in Figure~\eqref{f:eq-area}.  See Section~\ref{s:2dStationary}.}\label{f:162d}
\end{figure}

We  consider two-dimensional equal-area foams and evolve many random initial configurations until we obtain stationary. The random initial configurations are chosen as described in Section~\ref{s:2dDynamics}. 
In Figure~\ref{f:eq-area}, we plot the $n$-foams with the smallest total perimeter obtained for  $n=2,\ldots 21$. These results reproduce the results in  \cite{Cox_2003}.  We make the following observations: 
\begin{enumerate}
\item In all cases, Plateau's necessary conditions for optimality, discussed in Section~\ref{s:b2d}, are satisfied. 

\item For $n=2$ and $n=3$, we obtain the expected double and triple-bubble configurations. 

\item For $n$-foams with $n\leq 5$, there are no interior bubbles and for $n$-foams with $n\geq 6$, there appears to be at least one interior bubble. 

\item For $n=6,7,8$, we obtain $n$-foams with one interior bubble and $n-1$ boundary bubbles. For $n=6$ and $n=8$, due to the $120^\circ$ angle condition, the interior bubble is not a polygon, but has curved boundary.  

\item The configurations for some values of $n$ exhibit more symmetry than others. For example, $n=10,16$, and $20$ display additional symmetries. 

\item In Figure~\ref{f:162d}, another stationary equal-area  $16$-foam is given with slightly larger (numerically computed) total perimeter than the $16$-foam given in Figure~\ref{f:eq-area}. Interestingly, the 16-foam in Figure~\ref{f:162d} has more rotational symmetries than the 16-foam in Figure~\ref{f:eq-area}. It is also more similar to the 17-foam in Figure~\ref{f:eq-area}. 
\end{enumerate}

\subsection{Quasi-stationary flows corresponding to changing bubble size.} \label{s:2dQ}
\begin{figure}[t!]
\includegraphics[width=1 \textwidth, clip, trim=3cm 0cm 4cm 0cm]{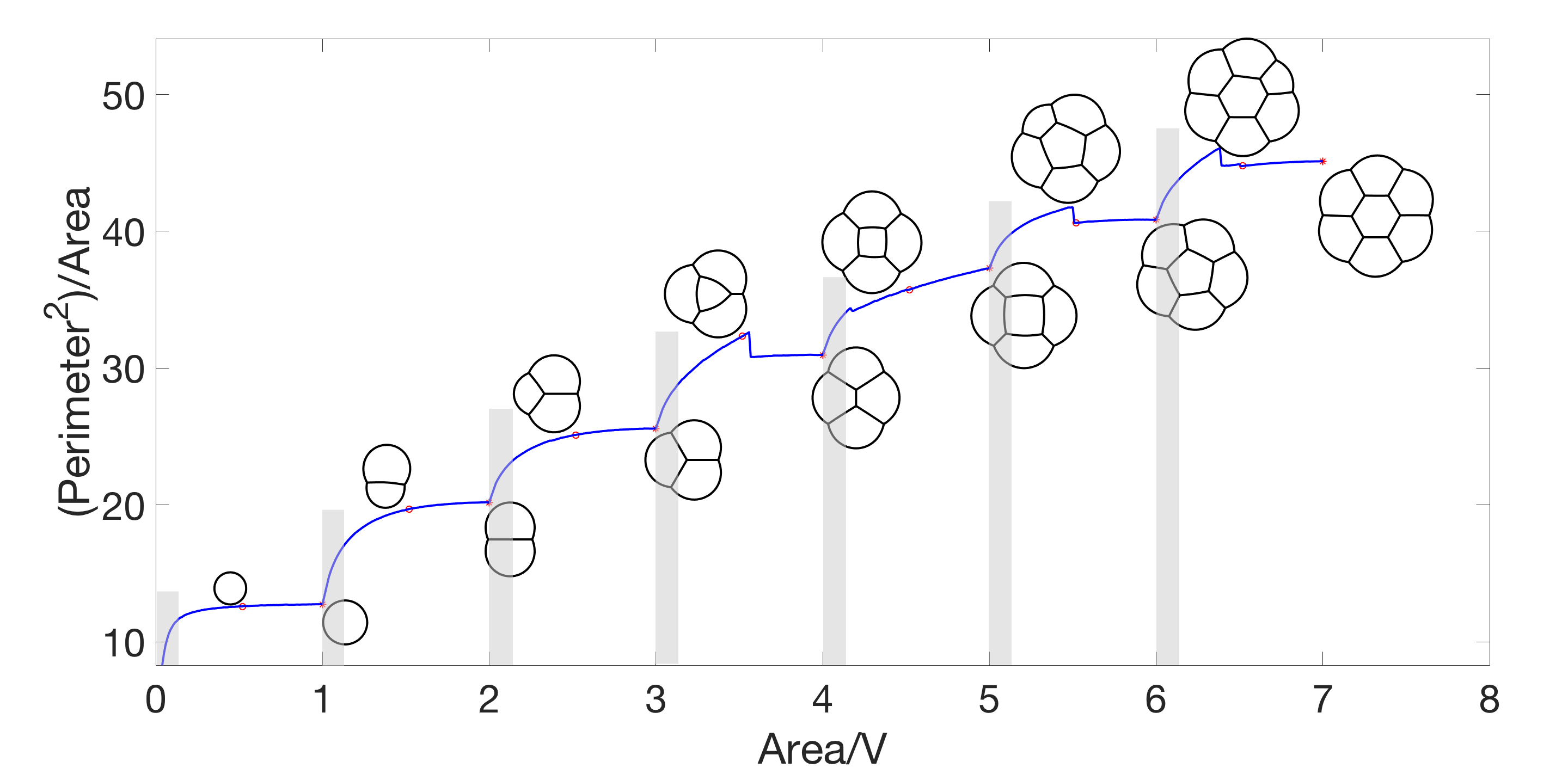} \\
\vspace{.5cm}
\includegraphics[width=1 \textwidth, clip, trim=3cm 0cm 4cm 0cm]{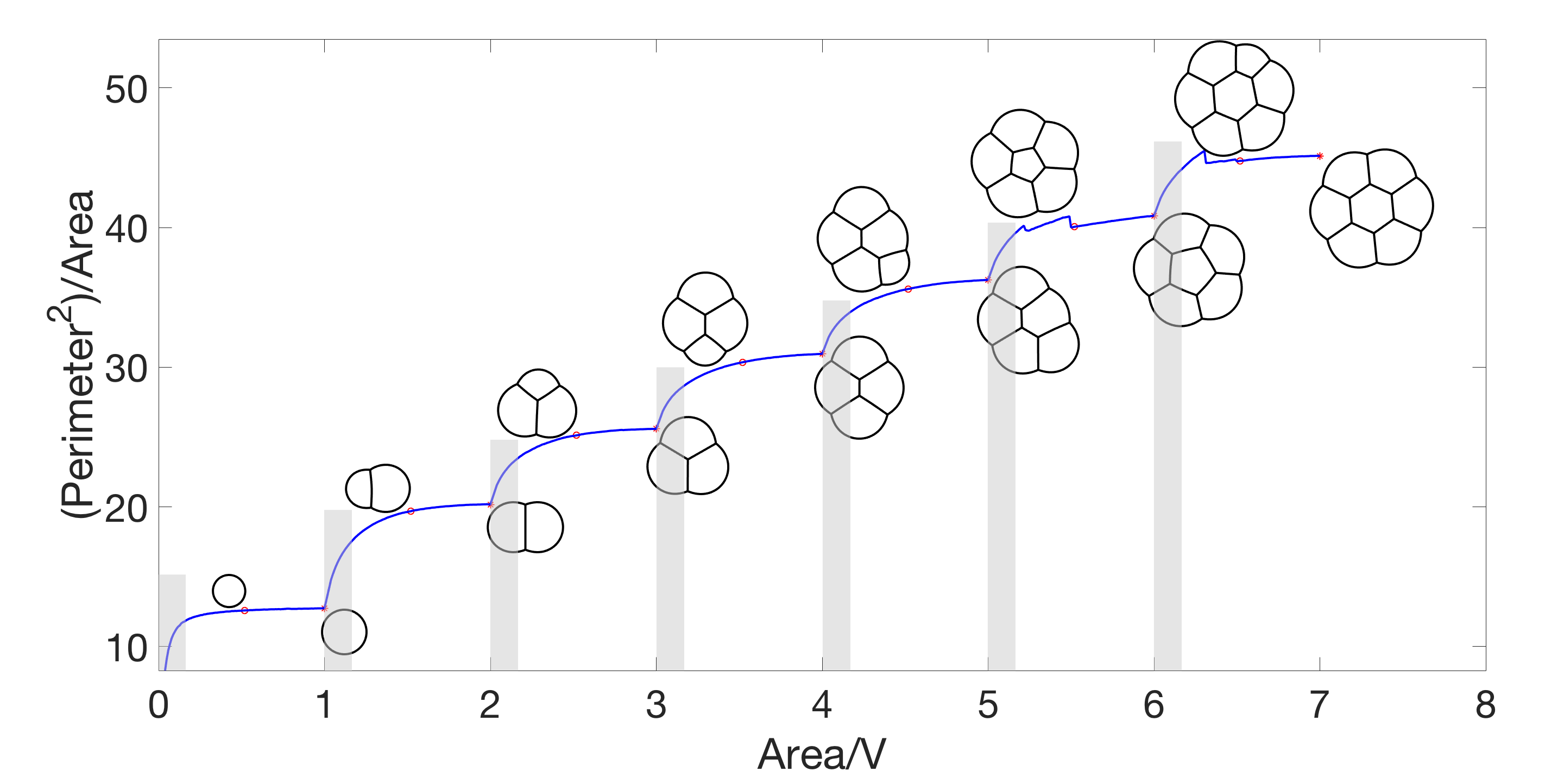} 

\caption{Total perimeter for the  quasi-stationary flow, where the area of one of the bubbles is slowly varied and for each fixed area, the stationary solution is computed. When the area reaches $V$, a new bubble with area $v$ is introduced. 
The top and bottom panels correspond to different positions where the new bubble is introduced. 
The foam configuration at various values of total area is plotted. 
Links to videos for this quasi-stationary flow are given in Table~\ref{t:2dLinks}. 
See Section~\ref{s:2dQ}.} 
\label{fig:increasing}
\end{figure}

We consider the configuration transition by increasing volume by $dV$ from only one bubble with small volume ($v$) to a fixed $V$ gradually. Then, we add another small bubble on the boundary of the cluster and increase the volume of this small bubble to $V$ gradually. With adding same bubbles at different positions, we obtain different paths of configuration transitions. Two example quasi-stationary flows are displayed in Figure \ref{fig:increasing}. In this example, $V = 0.4$, $dV = 0.004$, and $v= 0.016$. Links to corresponding videos are given in Table~\ref{t:2dLinks}.

\begin{rem} \label{rem:SmallVol}
The approximation $\mathcal H^{d-1}( \partial \Omega_i \cap \partial \Omega_j) \approx L(u_i, u_j)$, where $L(u_i, u_j)$ is defined in \eqref{lengthapprox}, has $O(\tau)$ accuracy. When the volume of one bubble is $o(\tau)$, this approximation is not very accurate. 
Of course, to resolve a smaller volume, the accuracy could be improved by using a smaller value of $\tau$. 
However, for a smaller $\tau$, the mesh must also be refined to avoid freezing at some non-stationary configuration, which makes the overall algorithm more computationally expensive. 
In Figure~\ref{fig:increasing}, we use gray rectangular boxes to indicate the regime where the results of the algorithm are not very convincing for the value $\tau=0.0625$ used. 
For example, when there is only one bubble, the isoperimetric quantity, $\frac{\textrm{Perimeter}^2}{\textrm{Area}}$, should be constant ($=4\pi$) and our numerical result agrees well with this value outside of the gray region. 
\end{rem}

\begin{table}[t!]
{\small
\begin{center}
\begin{tabular}{l l }
\multicolumn{2}{l}{Quasi-stationary flows corresponding to decreasing the area of  one bubble.}  \\
\hline 
Evolution from a $3$-foam to a $2$-foam: & \href{https://www.youtube.com/watch?v=LcX9iVE3cEk}{youtu.be/LcX9iVE3cEk}  \\
Evolution from a  $4$-foam to a $3$-foam: & \href{https://www.youtube.com/watch?v=t44JBQ4Cv9E}{youtu.be/t44JBQ4Cv9E}  \\
Evolution from a  $5$-foam to a $4$-foam: & \href{https://www.youtube.com/watch?v=uyRvH9CpQCM}{youtu.be/uyRvH9CpQCM}   \\
Evolution from a  $6$-foam to a $5$-foam: & \href{https://www.youtube.com/watch?v=Fs8XF6aNjEg}{youtu.be/Fs8XF6aNjEg}  \\
Evolution from a  $7$-foam to a $6$-foam: & \href{https://www.youtube.com/watch?v=w7p6E2Vcspg}{youtu.be/w7p6E2Vcspg}  \\
Evolution from a $8$-foam to a $7$-foam: & \href{https://www.youtube.com/watch?v=s0XNdaJP364}{youtu.be/s0XNdaJP364}  \\
Evolution from a $9$-foam to a $8$-foam: & \href{https://www.youtube.com/watch?v=XBiQRvjgDVQ}{youtu.be/XBiQRvjgDVQ} \\
\hline \\
\multicolumn{2}{l}{Quasi-stationary flows corresponding to increasing the area of  one bubble.}  \\
\hline
Evolution from a $2$-foam to a $3$-foam: & \href{https://www.youtube.com/watch?v=dfPmFPD4Atw}{youtu.be/dfPmFPD4Atw}  \\
Evolution from a $3$-foam to a $4$-foam: & \href{https://www.youtube.com/watch?v=cFHXdMwFo7M}{youtu.be/cFHXdMwFo7M}  \\
Evolution from a $4$-foam to a $5$-foam: & \href{https://www.youtube.com/watch?v=j7-5L9ff_xg}{youtu.be/j7-5L9ff\_xg}   \\
Evolution from a $5$-foam to a $6$-foam: & 
\href{https://www.youtube.com/watch?v=m85uyeiQ2BM}{youtu.be/m85uyeiQ2BM}  \\
&\href{https://www.youtube.com/watch?v=0KpHnPKl0tA}{youtu.be/0KpHnPKl0tA}  \\
&\href{https://www.youtube.com/watch?v=jatMSRAxYfQ}{youtu.be/jatMSRAxYfQ}   \\
Evolution from a $6$-foam to a $7$-foam: & 
\href{https://www.youtube.com/watch?v=BP0z93JULCE}{youtu.be/BP0z93JULCE}   \\
\hline
\end{tabular}
\end{center}
\caption{Links to videos showing the quasi-stationary flow as the area of one bubble is either increased or decreased. Example foam configurations from this flow are shown in Figure~\ref{fig:increasing}. See Section~\ref{s:2dQ}.}
\label{t:2dLinks}
}
\end{table}%

\subsection{Configuration transitions}
\label{s:2dConfTran}
The problem of finding minimal total perimeter foams \eqref{e:min} possesses several local solutions corresponding to distinct foam configurations which are well-separated and have almost the same total perimeter. 
When the problem is perturbed (\eg, the volume of one of the bubbles increases or decreases), these local minima vary. 
As we perturb the problem, we observe \emph{configuration transitions} where a local minima rapidly transitions and converges to another local minima. This is demonstrated in Figure~\ref{fig:increasing}, where there are small jumps in the energy curve. In this section we further study this phenomena. 

\begin{figure}[t]
\includegraphics[width=1 \textwidth, clip, trim=5cm 0cm 5cm 2cm]{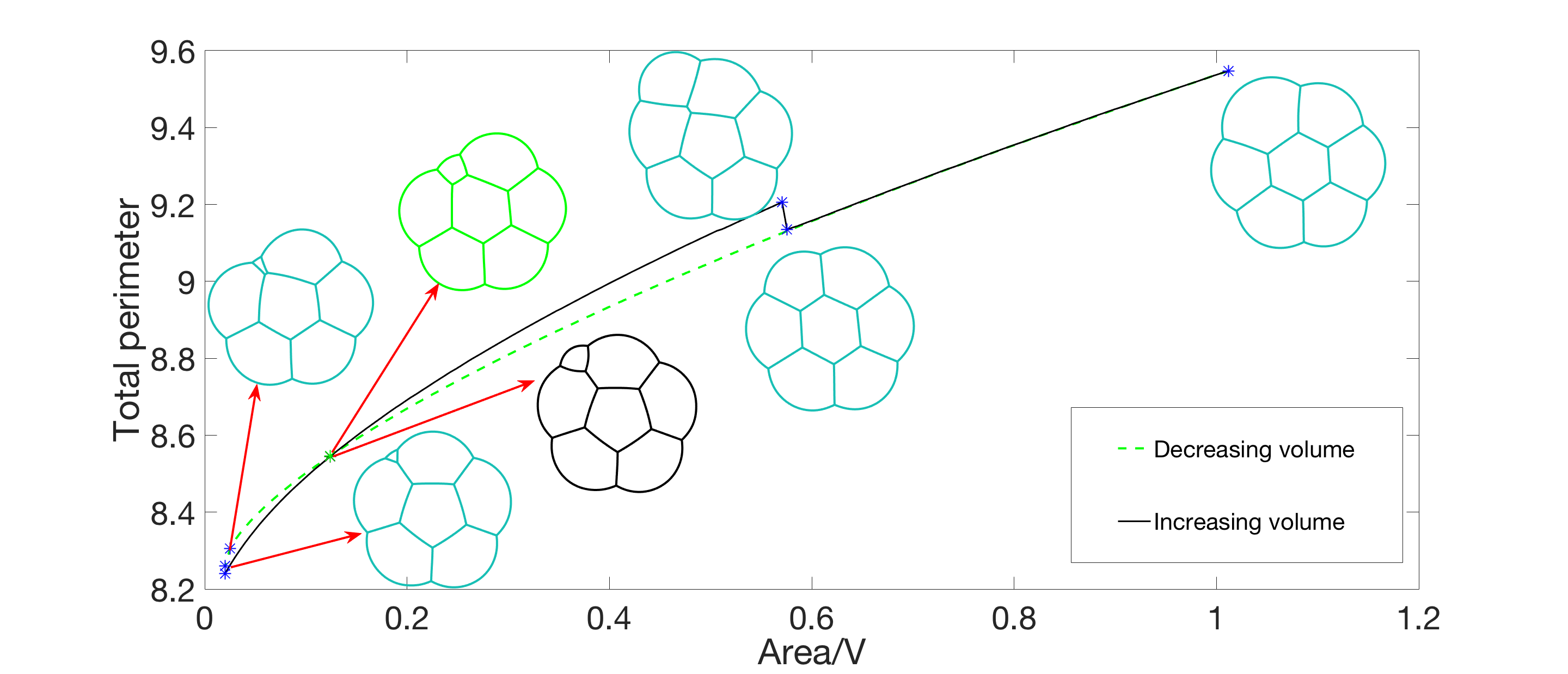}
\caption{Energy plot of increasing and decreasing area with snapshots at different value of area. See Section~\ref{s:2dConfTran}.} \label{fig:hysteresis}
\end{figure}

By considering the system with $6$ bubbles with equal area $V$ and one small bubble with area $v$, we gradually increase the volume of the small bubble to $1.5V$ and then decrease the volume of this bubble to the original area $v$. The energy plot is displayed in Figure~\ref{fig:hysteresis}. The black line is the energy plot for increasing area and the green dashed line is the energy plot for decreasing area. The jumps on the black and dashed green lines are positions of configuration transitions.  We also note that the intersection between the black line and dashed green line correspond two different configurations. These two configurations have the same energy and same areas of bubbles. Interestingly, from this experiment, we see that the process of increasing and decreasing volume are irreversible; one can view this as a type of hysteresis in the sense that the flow depends on the initialization. In this example, $V = 0.677$, $dV = 0.00496$, and $v= 0.0201$.

\begin{figure}[t]
\includegraphics[width=1 \textwidth, clip, trim=5cm 0cm 5cm 2cm]{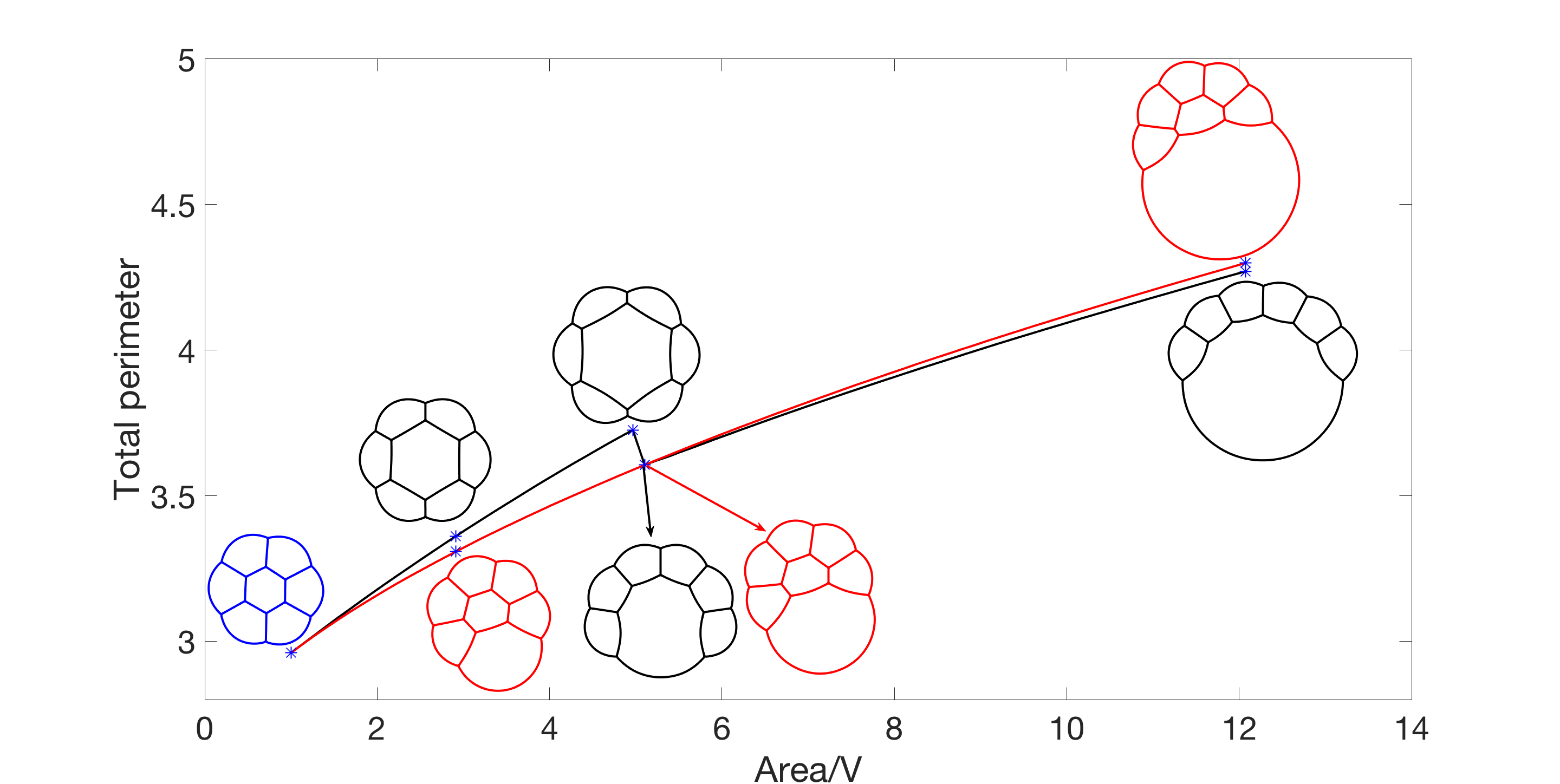} \\

\vspace{.2cm}

\begin{tabular}{l l }
\multicolumn{2}{l}{Quasi-stationary flows corresponding to increasing the area of  one bubble.}  \\
\hline 
Increasing the area of the middle bubble: & \href{https://youtu.be/--HWXssRERk}{youtu.be/--HWXssRERk}  \\
Increasing the area of a border bubble: & \href{https://youtu.be/cJsbU1mtT3E}{youtu.be/cJsbU1mtT3E}  \\
\hline
\end{tabular}
\caption{Energy plot of increasing area at different positions with snapshots at different value of area and links to videos showing the quasi-stationary flow as the area of one bubble is increased either in the middle or the border. See Section~\ref{s:2dConfTran}.} \label{fig:MiddleVsBorder}
\end{figure}

Also, from Figure~\ref{fig:increasing}, we see different computed stationary configurations when we add area to one bubble at different positions. To further study this, starting from the computed stationary configuration for an equal-area $7$-foam (see Figure~\ref{f:eq-area}), $V$,  we gradually add area to one bubble until the area is $12V$. We compare the difference between adding area to the middle bubble and adding the area to the border bubble. In Figure~\ref{fig:MiddleVsBorder}, the black line displays the change in total perimeter when we increase the area of the middle bubble while the red line displays the change in total perimeter when we increase the area of a border bubble starting from the same initial configuration which is plotted in blue lines. The snapshots of increasing the area of the middle bubble are plotted in black and the snapshots of increasing the area of a border bubble are plotted in red. In this example, $V = 0.1474$ and $dV = 0.02$. The links for the corresponding videos are given in Figure~\ref{fig:MiddleVsBorder}. 

\section{Three-dimensional numerical examples} \label{s:3d}

\subsection{Time-evolution of foams}\label{s:3dDynamics} 
\begin{figure}[t]
\centering
\includegraphics[width=0.95 \textwidth, clip, trim=5cm 0cm 5cm 2cm]{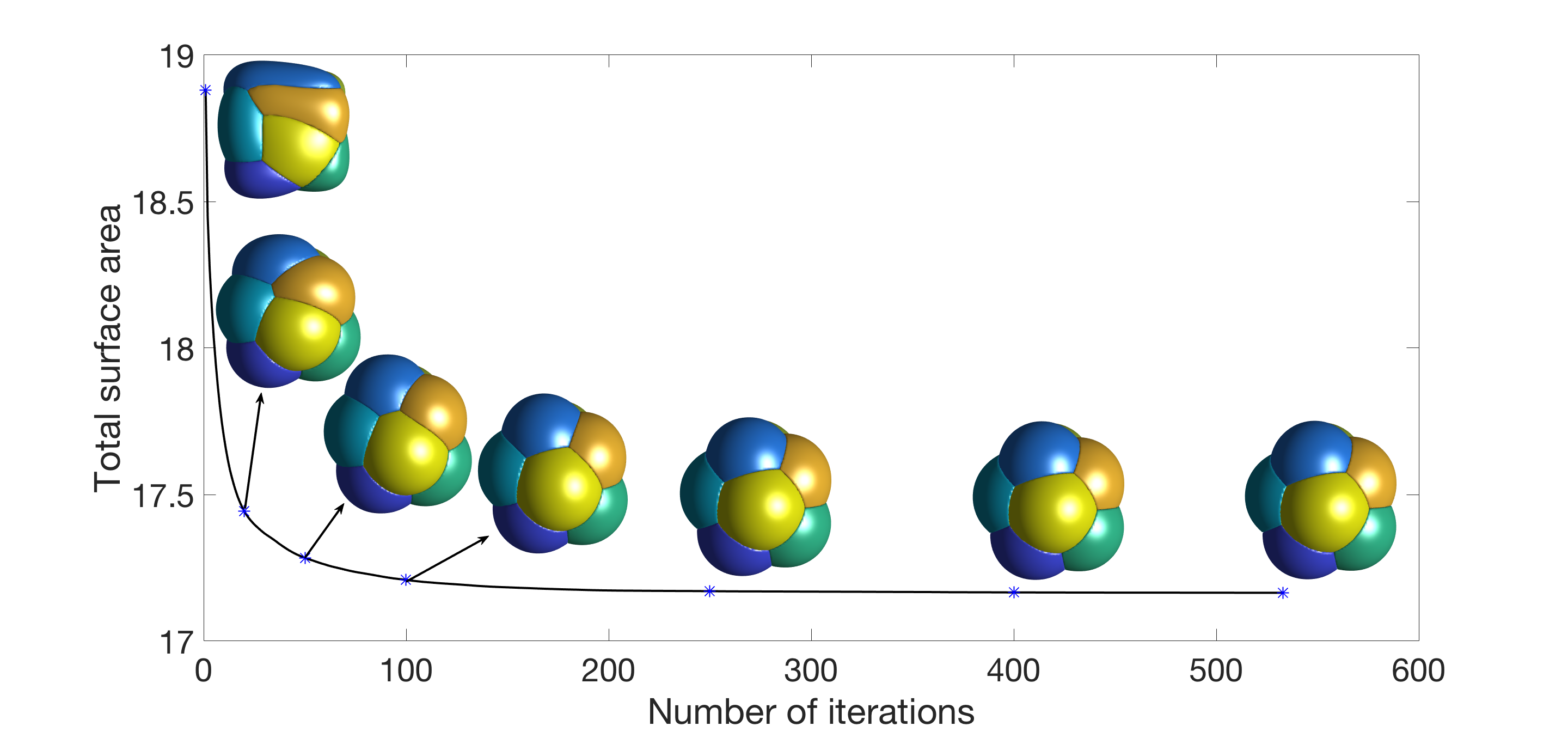}
\caption{A plot of the energy as a $n=8$-foam evolves from a random initialization together with the foam configuration at various iterations. See Section~\ref{s:3dDynamics}.}\label{fig:dynamics3d}
\end{figure}

In Figure~\ref{fig:dynamics3d}, for an equal-volume, $n=8$-foam, we show the time evolution corresponding to the gradient flow of the  total surface area with a random initialization; the initial configuration was chosen as in the two-dimensional flow described in Section~\ref{s:2dDynamics}. The energy at each iteration is plotted together with the foam configuration at various iterations. Note that the energy decays very fast; even in three-dimensional space, after $533$ iterations, the configuration is stationary in the sense that no grid points are changing bubble membership. After $\approx 150$ iterations, the foam configuration changes very little. 

\subsection{Stationary solutions}  \label{s:3dStationary}
\begin{figure}[ht]
\begin{center}
\includegraphics[width=.24\textwidth,clip,trim= 5cm 3cm 5cm 2.5cm]{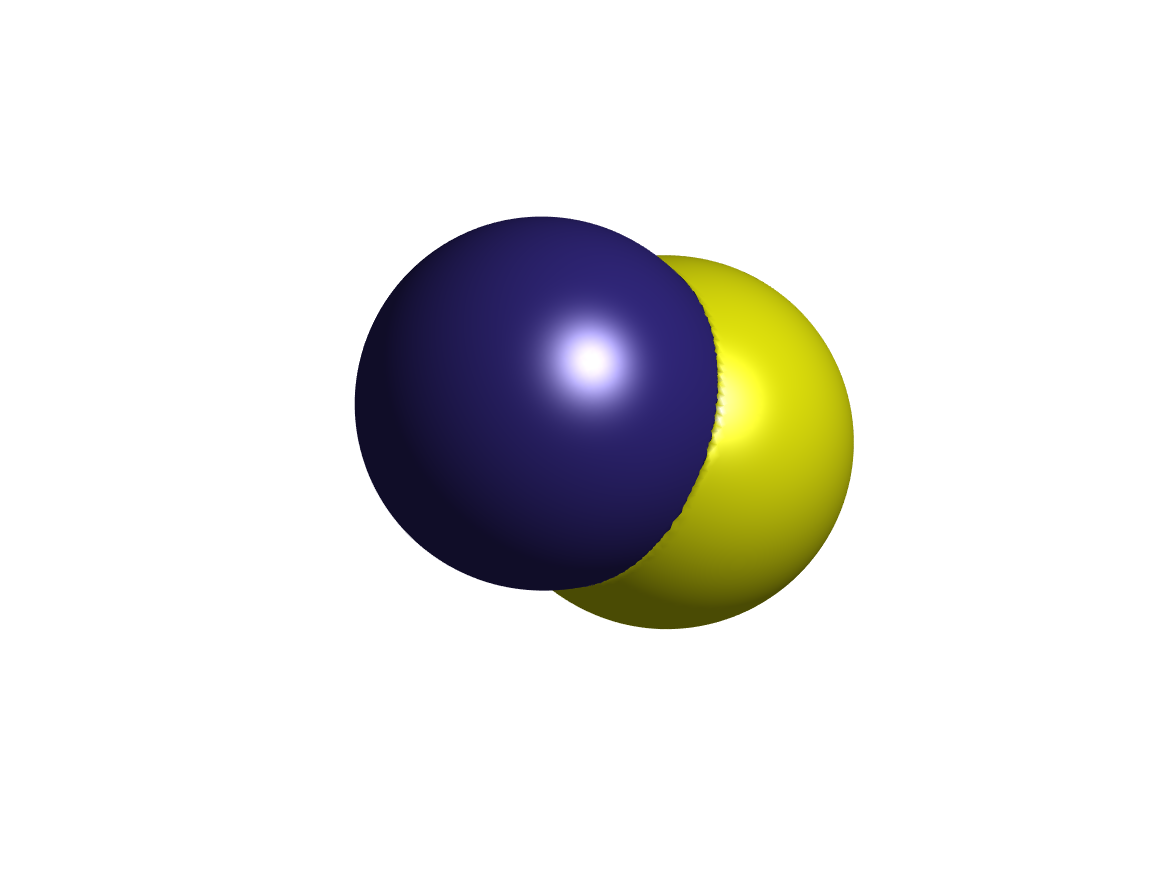}
\includegraphics[width=.24\textwidth,clip,trim= 5cm 3cm 5cm 2.5cm]{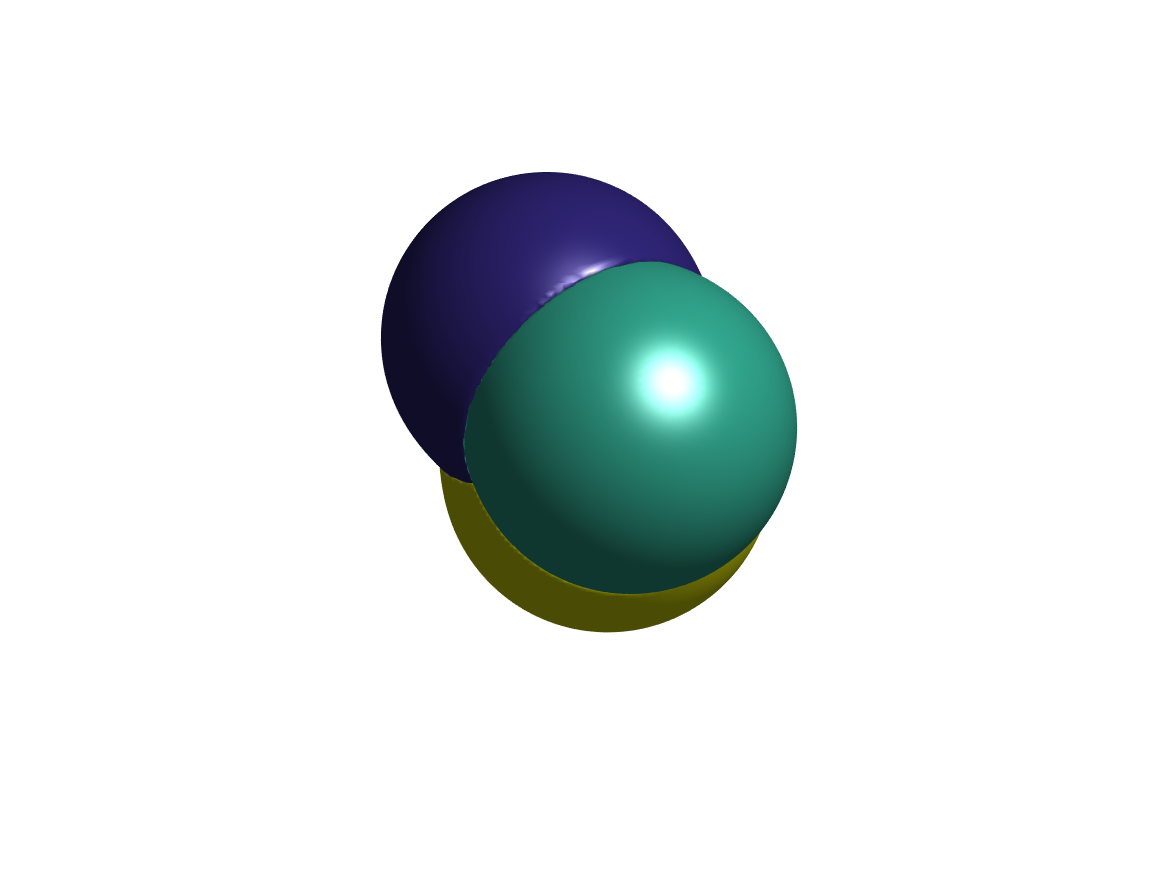}
\includegraphics[width=.24\textwidth,clip,trim= 5cm 3cm 5cm 2.5cm]{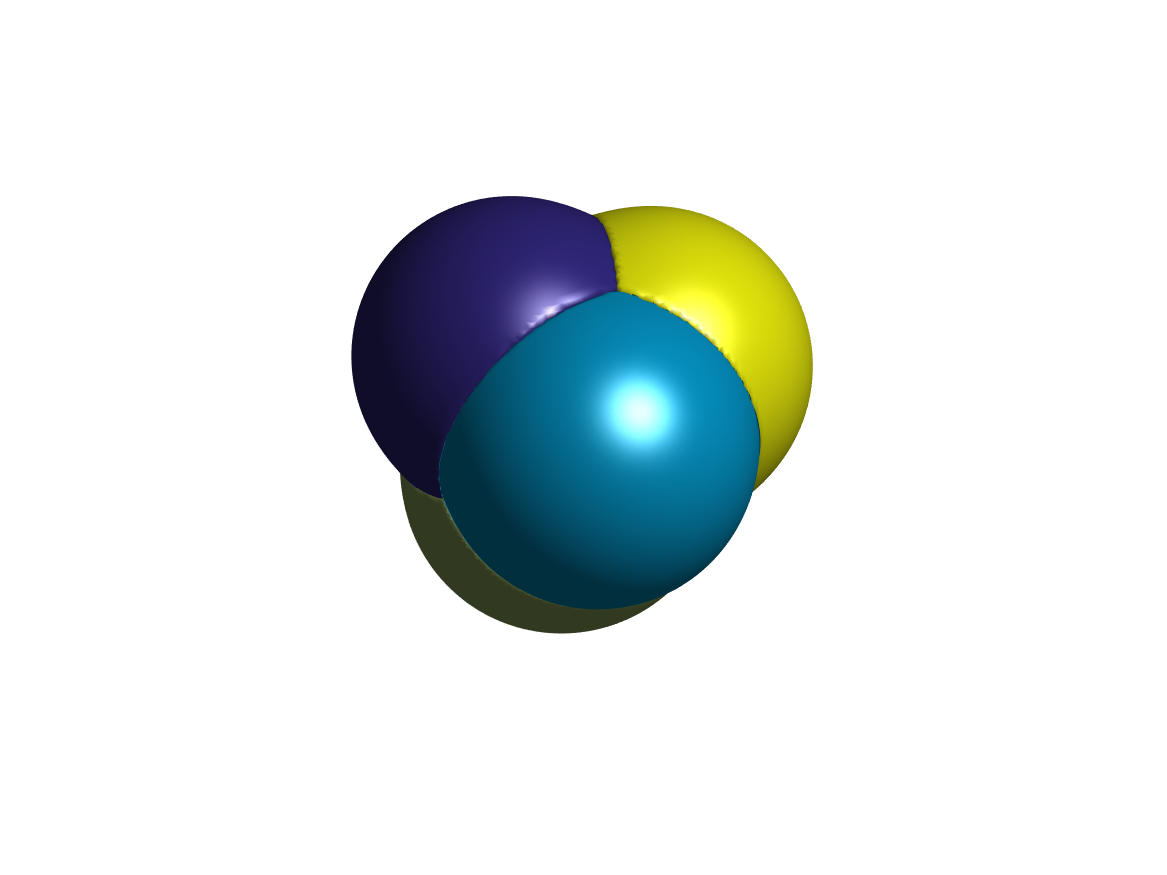}
\includegraphics[width=.24\textwidth,clip,trim= 5cm 3cm 5cm 2.5cm]{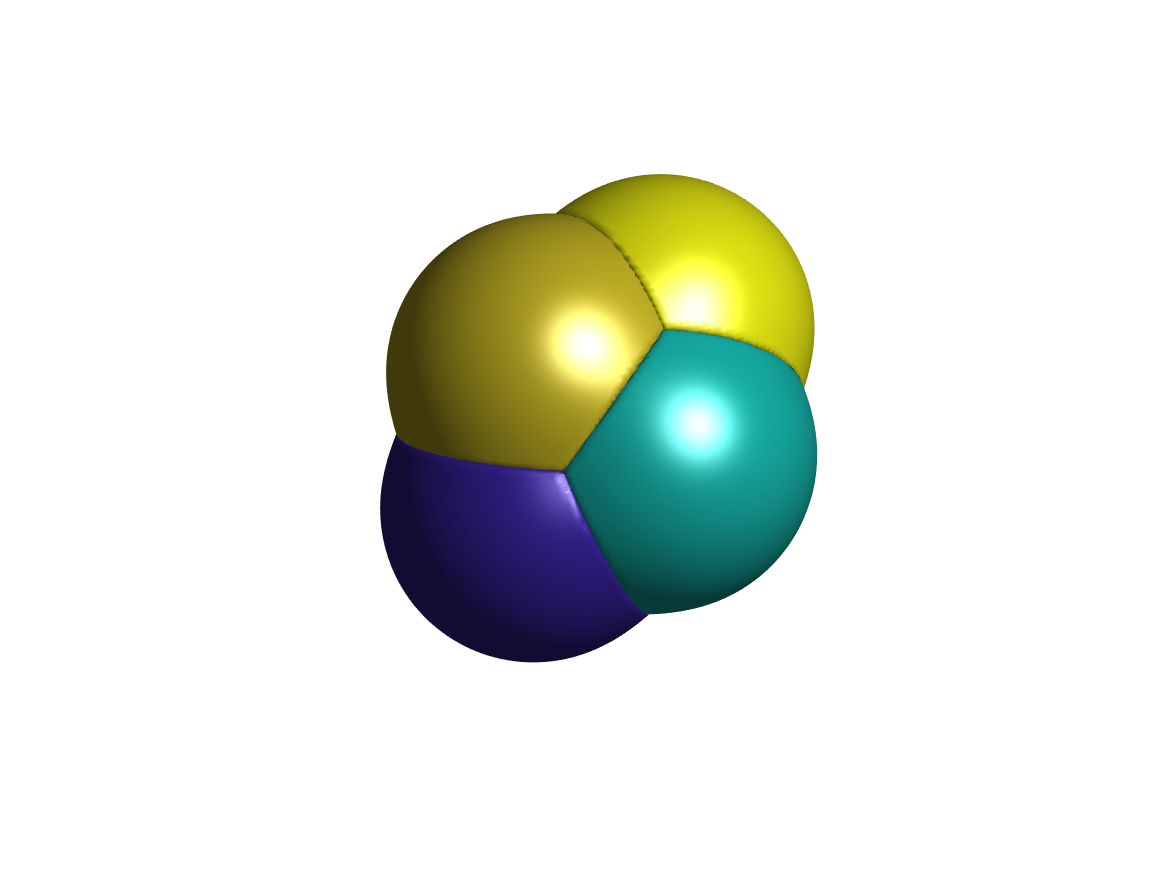}
\includegraphics[width=.24\textwidth,clip,trim= 5cm 3cm 5cm 2.5cm]{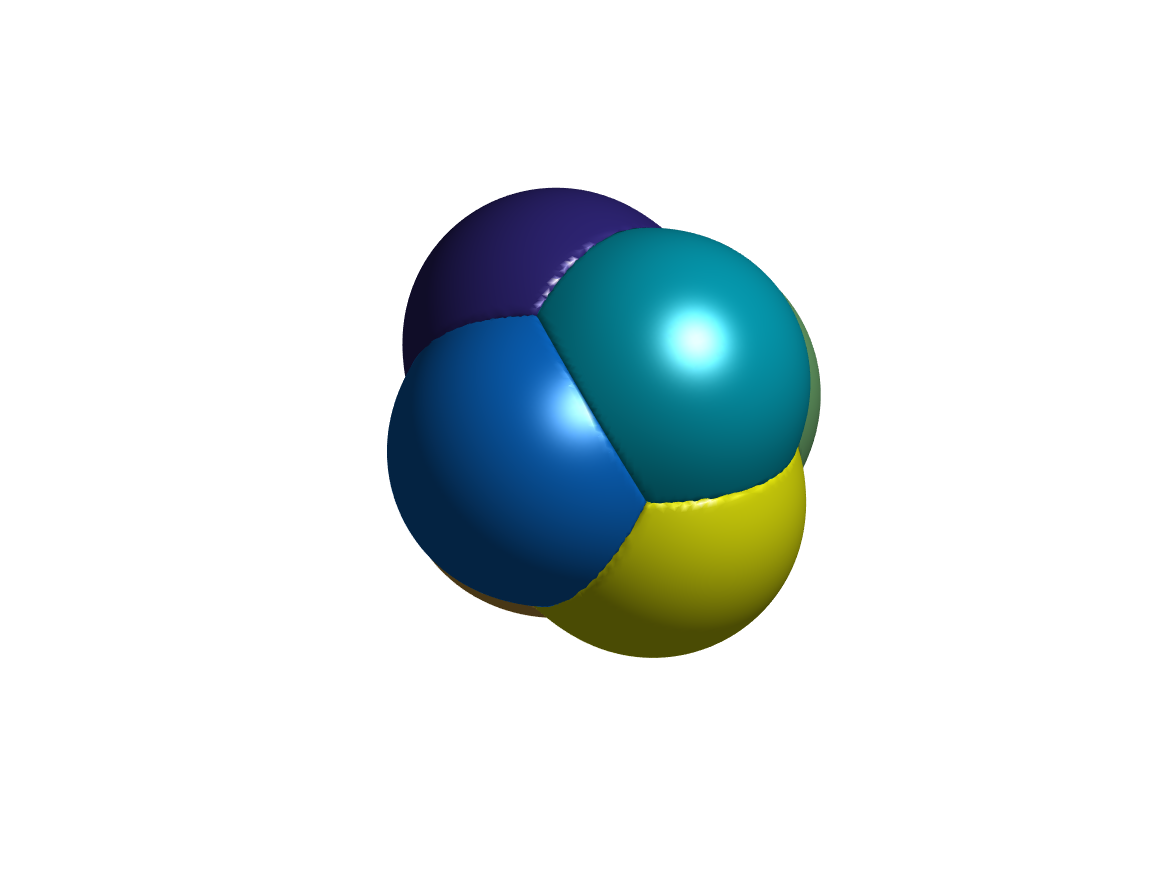}
\includegraphics[width=.24\textwidth,clip,trim= 5cm 3cm 5cm 2.5cm]{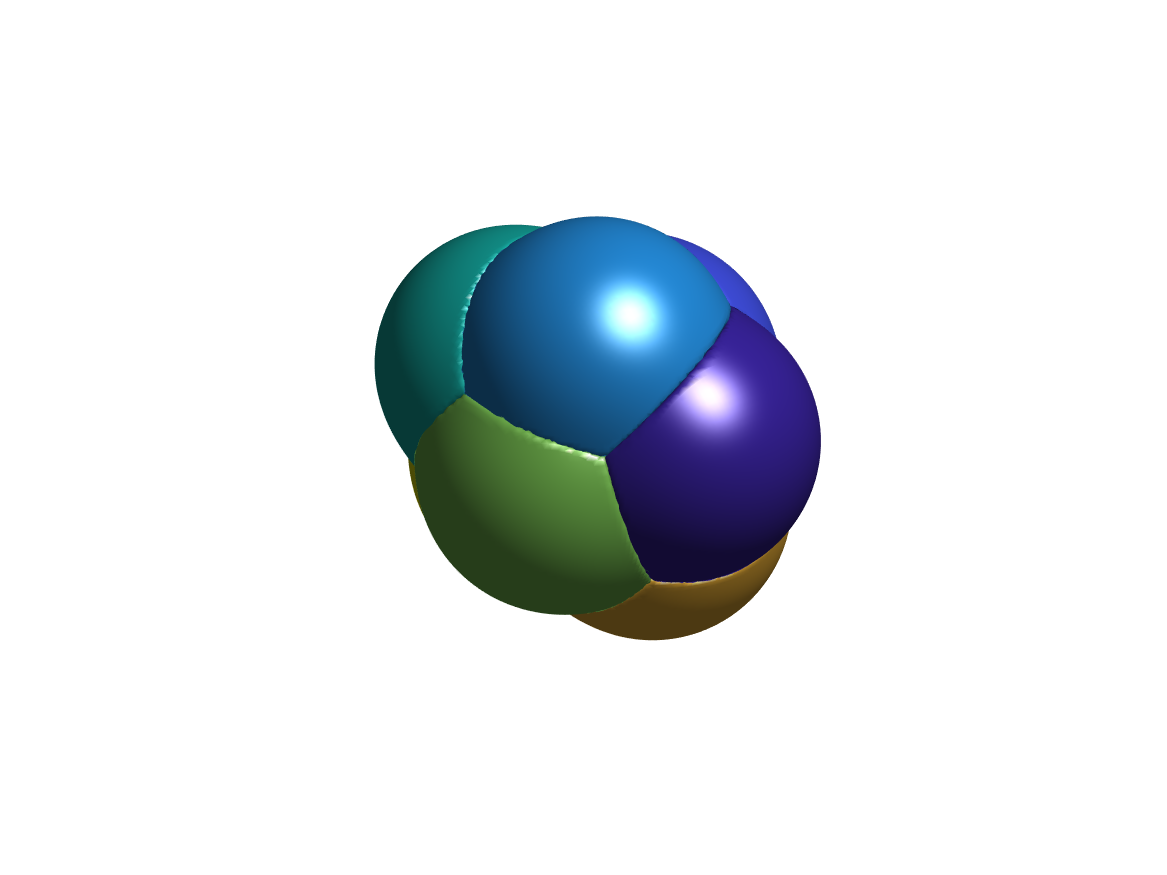}
\includegraphics[width=.24\textwidth,clip,trim= 5cm 3cm 5cm 2.5cm]{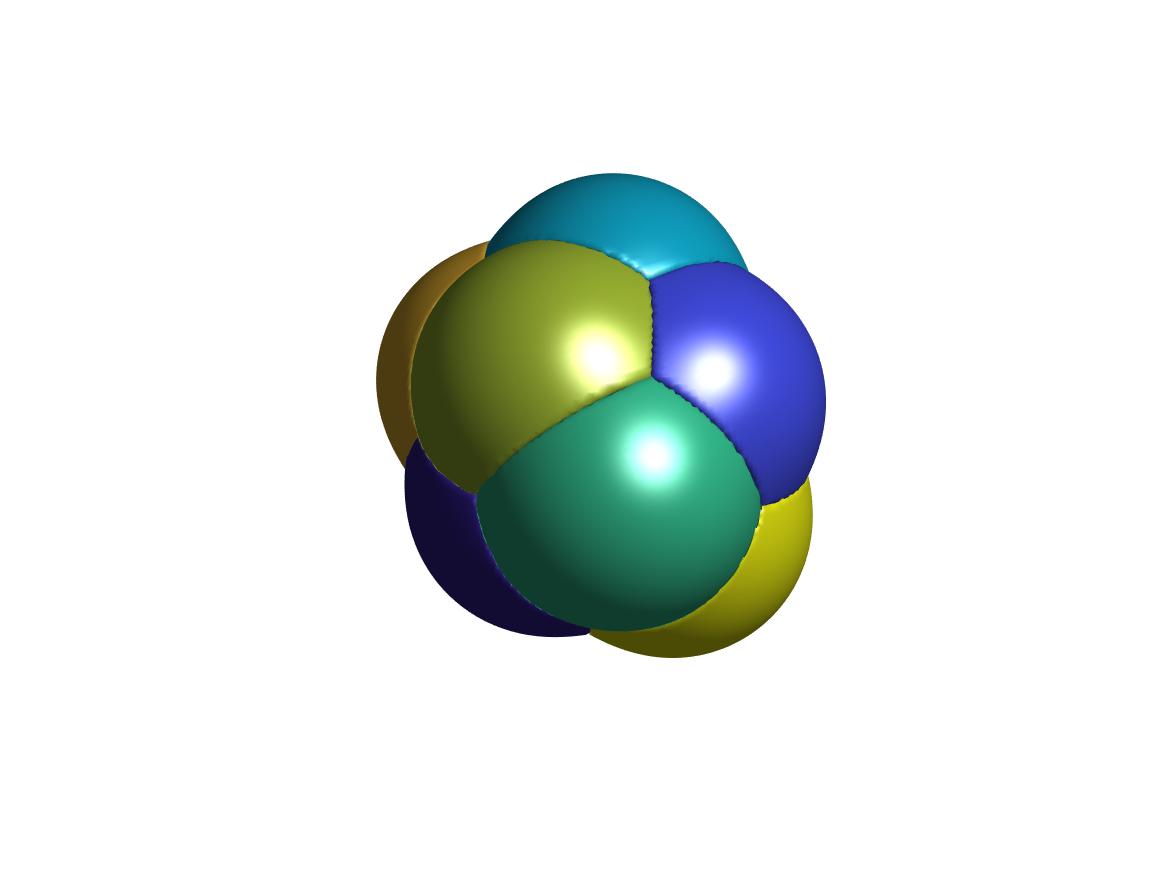}
\includegraphics[width=.24\textwidth,clip,trim= 5cm 3cm 5cm 2.5cm]{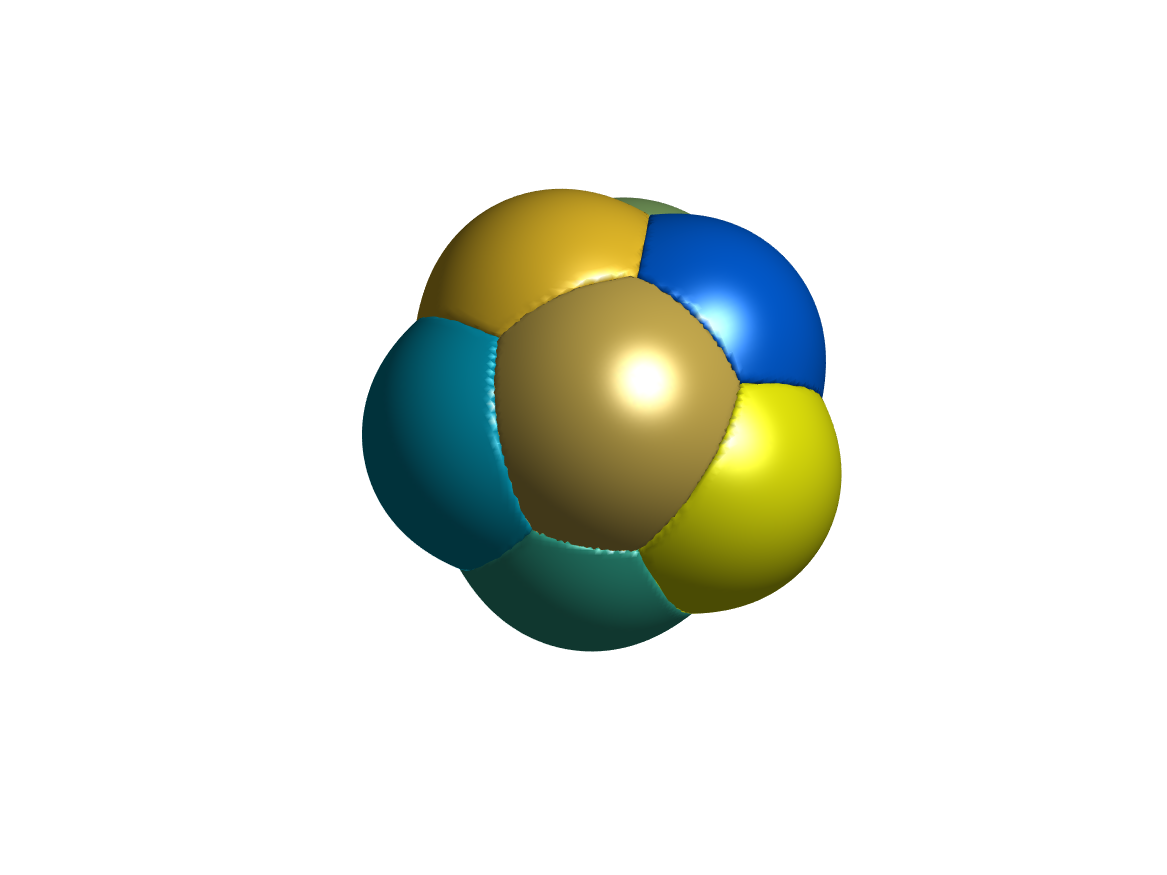}
\includegraphics[width=.24\textwidth,clip,trim= 5cm 3cm 5cm 2.5cm]{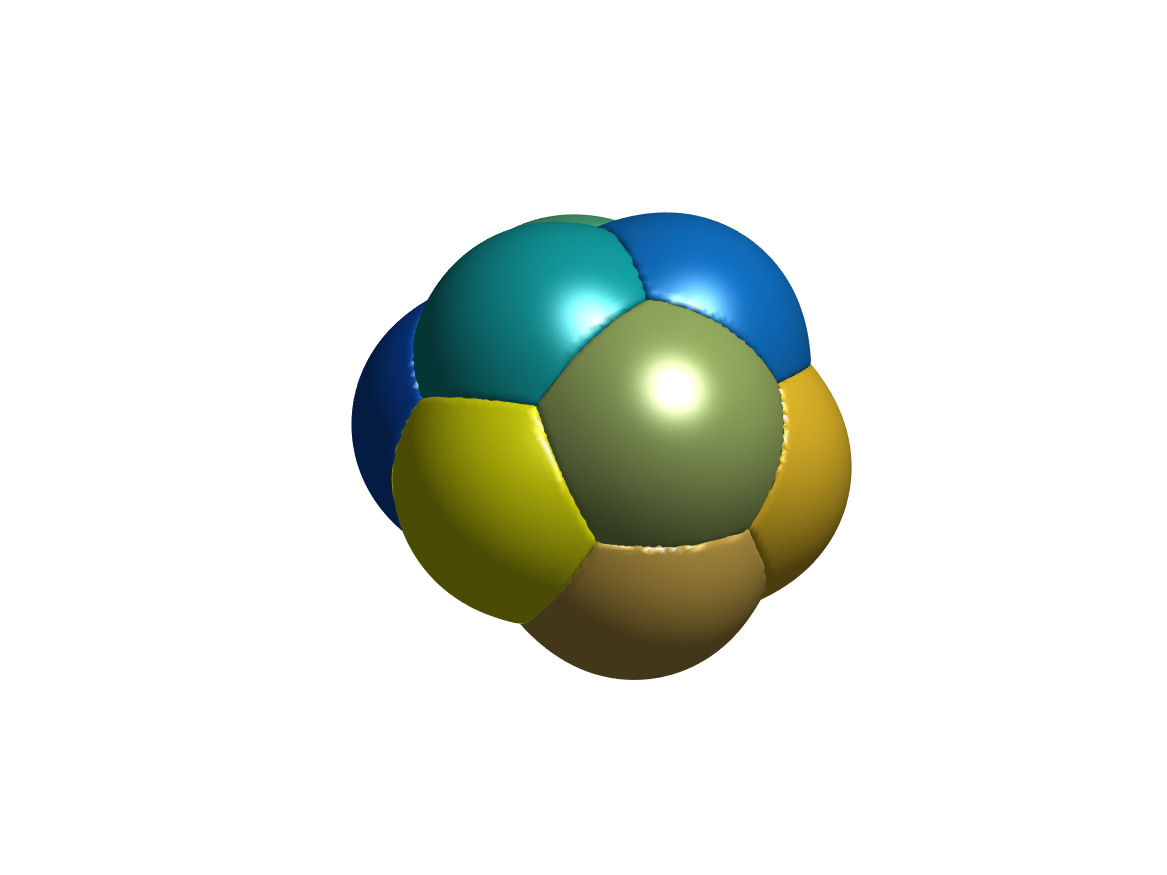}
\includegraphics[width=.24\textwidth,clip,trim= 5cm 3cm 5cm 2.5cm]{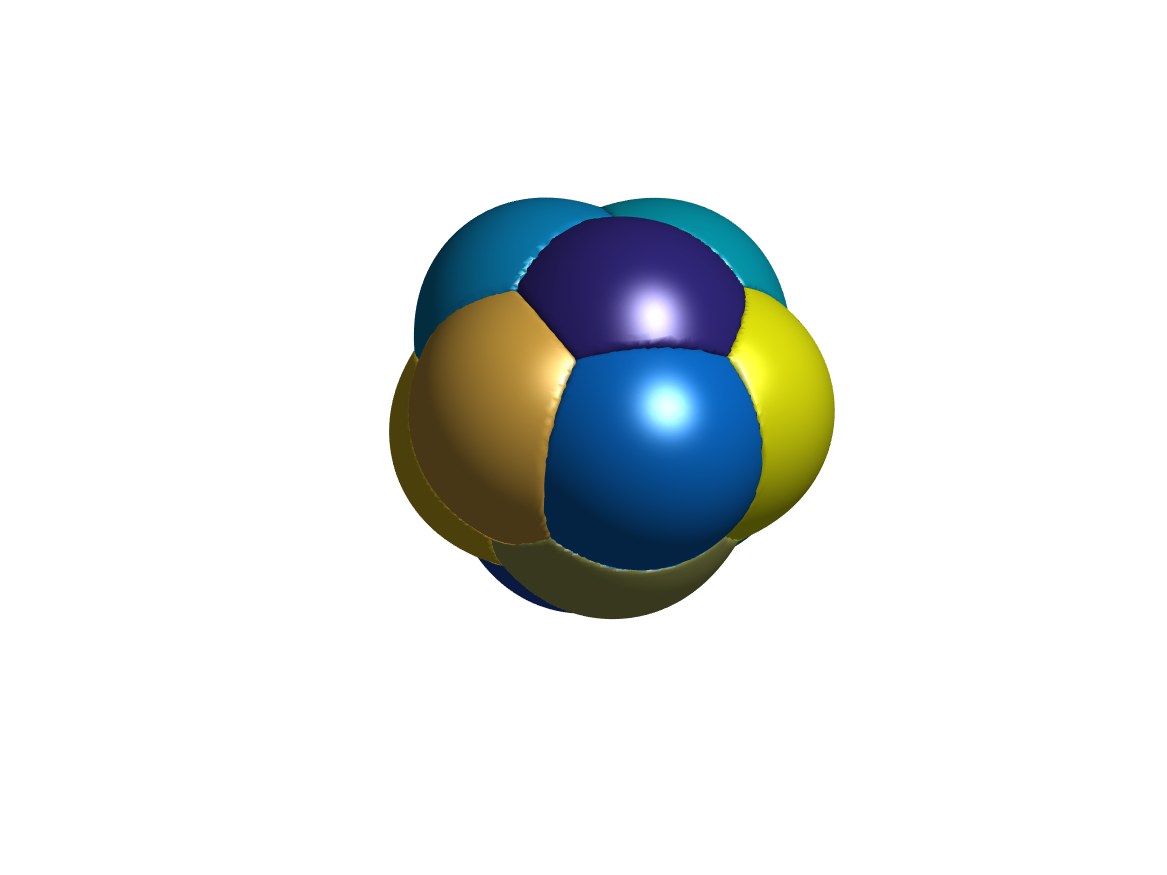}
\includegraphics[width=.24\textwidth,clip,trim= 5cm 3cm 5cm 2.5cm]{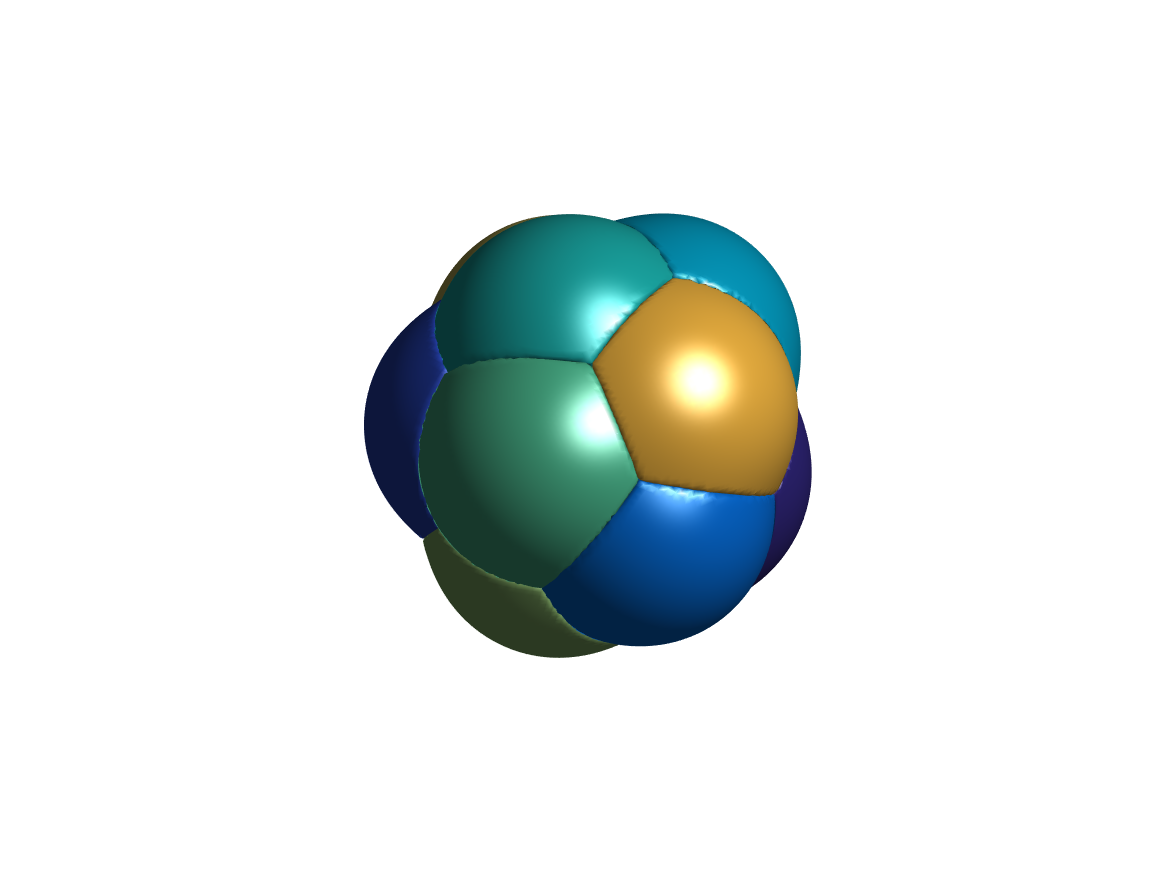}
\includegraphics[width=.24\textwidth,clip,trim= 5cm 3cm 5cm 2.5cm]{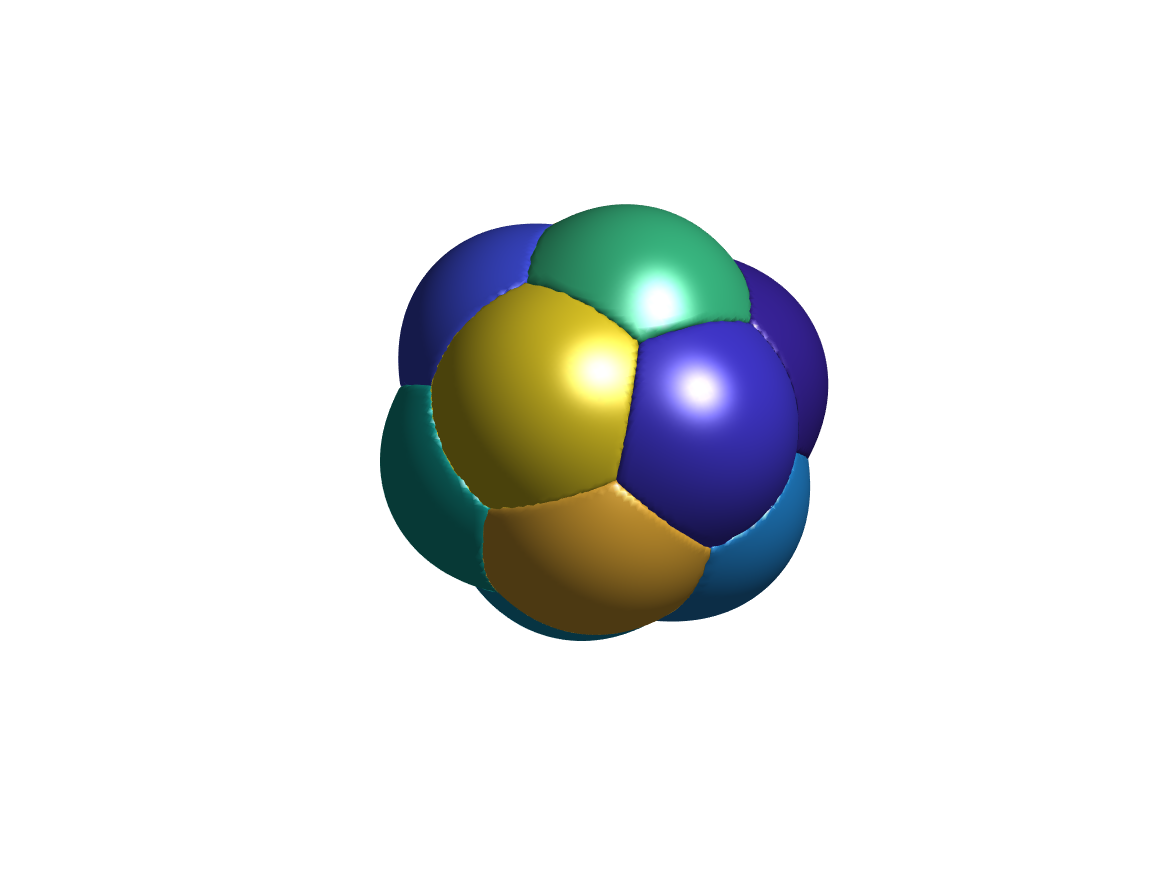}
\includegraphics[width=.24\textwidth,clip,trim= 5cm 3cm 5cm 2.5cm]{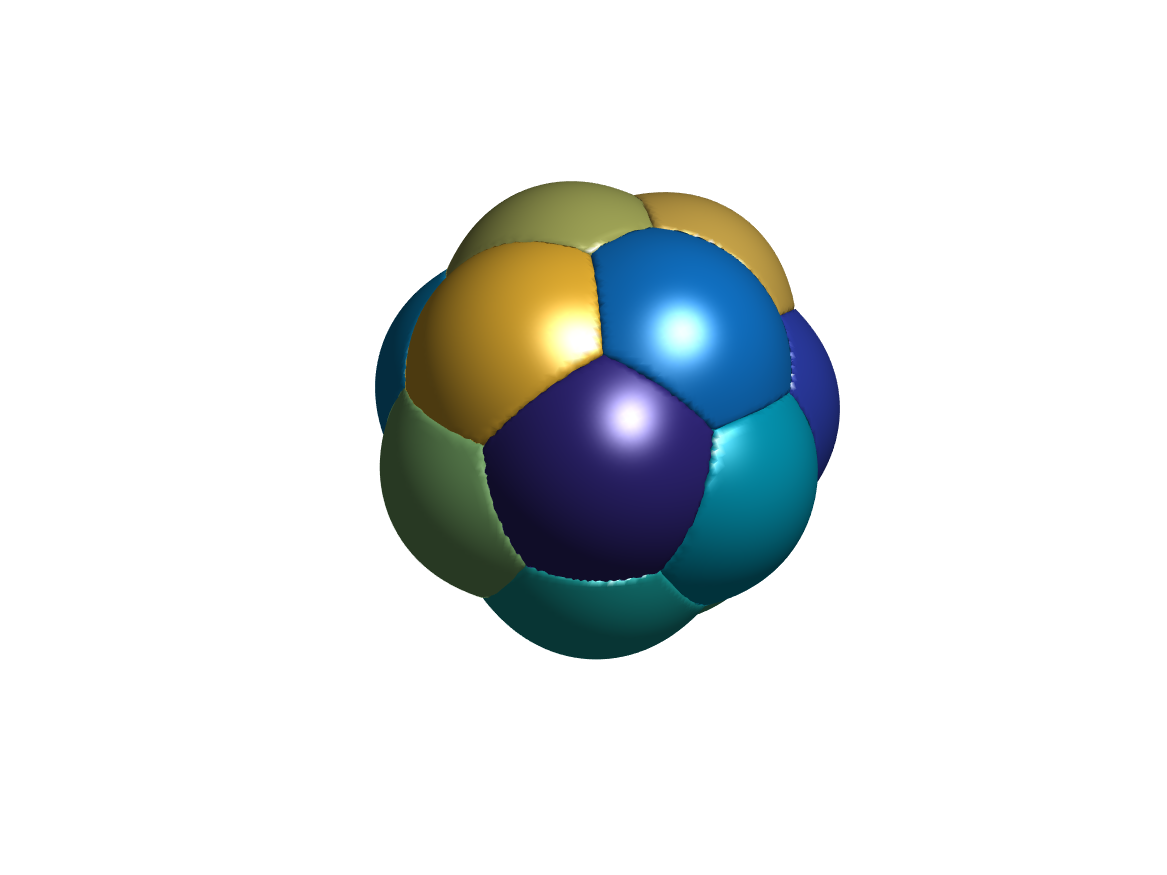}
\includegraphics[width=.24\textwidth,clip,trim= 5cm 3cm 5cm 2.5cm]{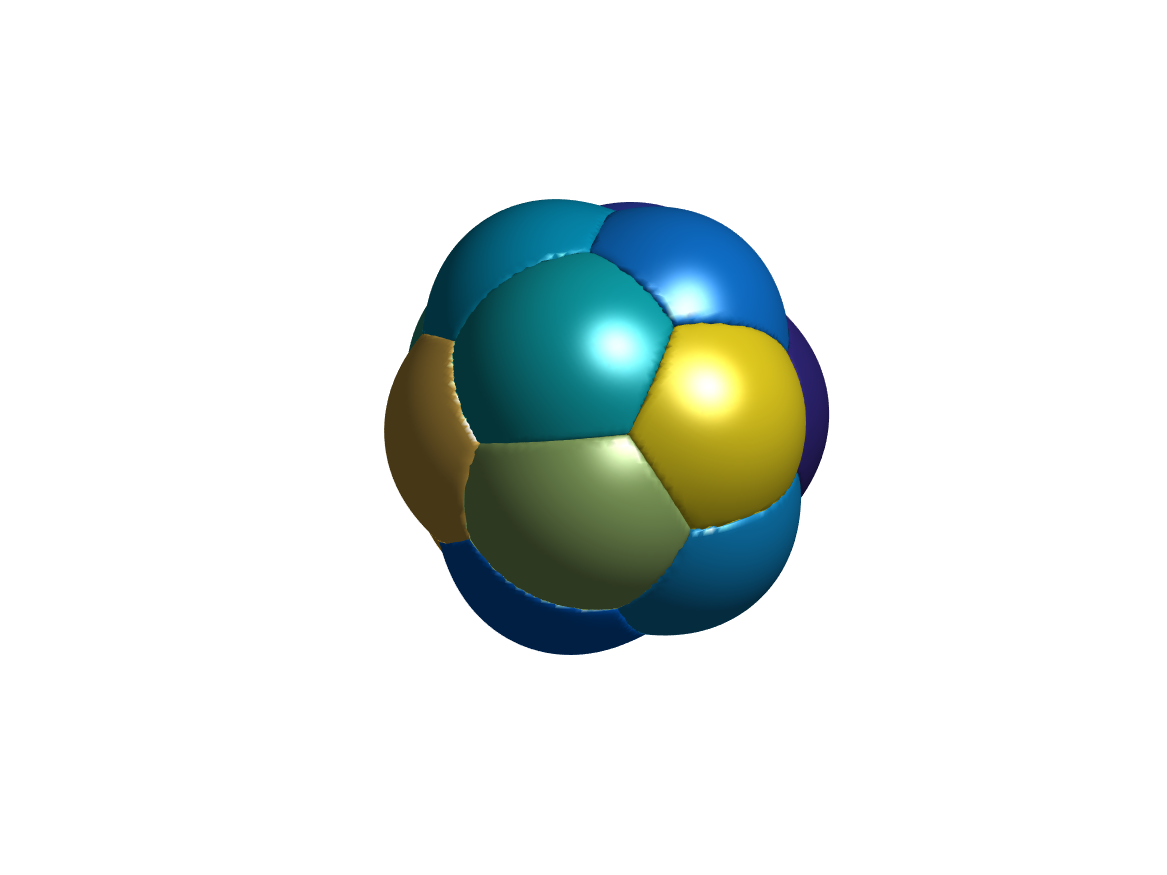}
\includegraphics[width=.24\textwidth,clip,trim= 5cm 3cm 5cm 2.5cm]{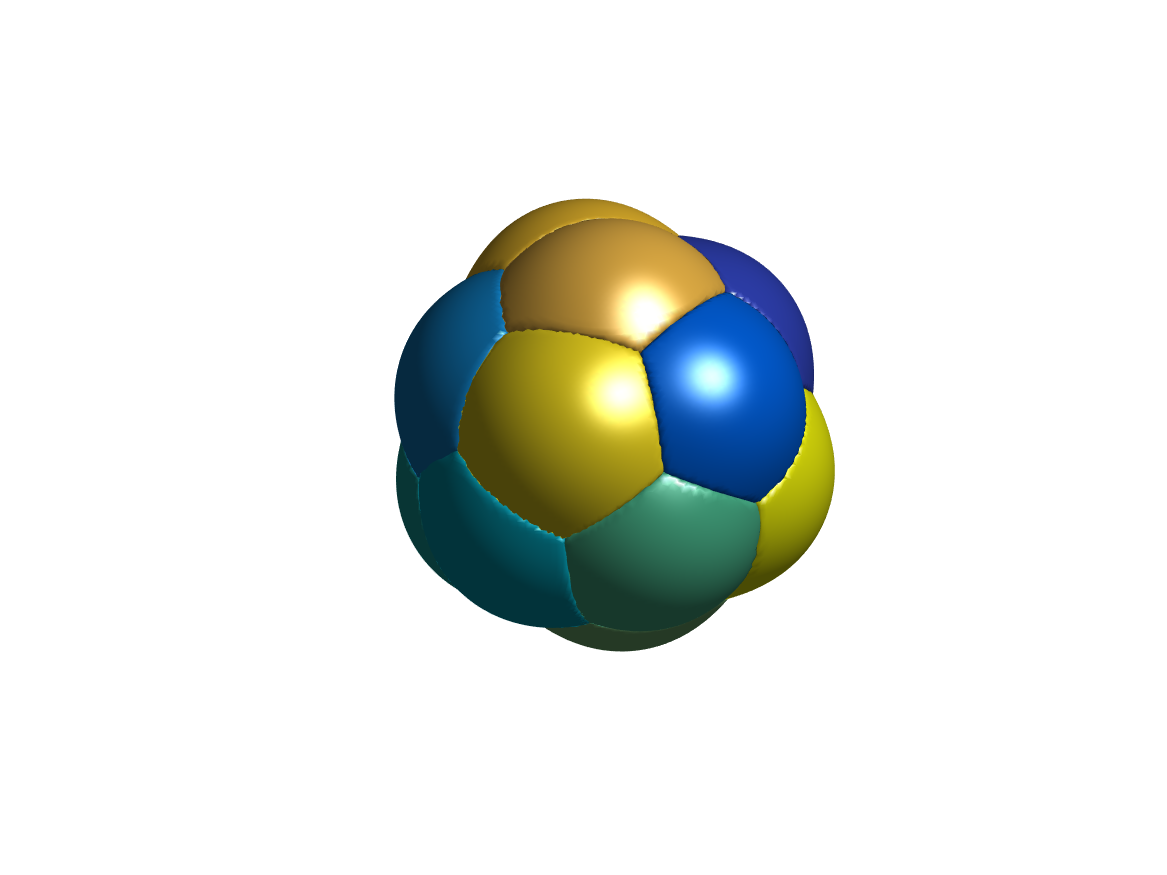}
\includegraphics[width=.24\textwidth,clip,trim= 5cm 3cm 5cm 2.5cm]{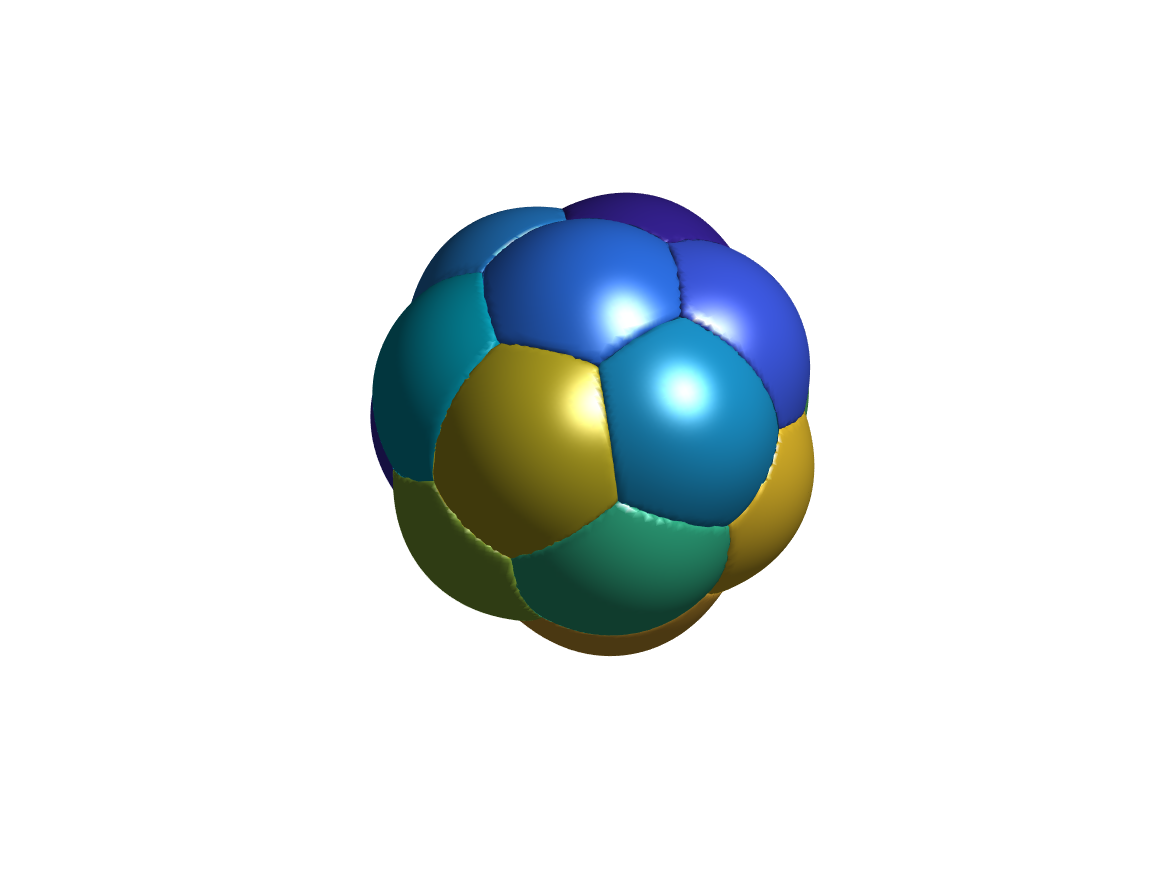} \\
\vspace{.5cm}
\begin{tabular}{ll|ll}
\hline
$2$-foam: &  \href{https://youtu.be/UrqE4vmkxnc}{youtu.be/UrqE4vmkxnc} 
& $3$-foam: &  \href{https://youtu.be/oQxR6Z_fpTE}{youtu.be/oQxR6Z\_fpTE} \\
$4$-foam: &  \href{https://youtu.be/LQEphY2Ctq4}{youtu.be/LQEphY2Ctq4} 
& $5$-foam: &  \href{https://youtu.be/EDEwdMR21Xo}{youtu.be/EDEwdMR21Xo} \\
$6$-foam: &  \href{https://youtu.be/-PlS75F6Ueo}{youtu.be/-PlS75F6Ueo} 
&$7$-foam: &  \href{https://youtu.be/Avke4wfKADY}{youtu.be/Avke4wfKADY} \\
$8$-foam: &  \href{https://youtu.be/tI1zb685MAs}{youtu.be/tI1zb685MAs} 
&$9$-foam: &  \href{https://youtu.be/eu2OYEC7KUE}{youtu.be/eu2OYEC7KUE} \\
$10$-foam: &  \href{https://youtu.be/J9mlXTLNuxc}{youtu.be/J9mlXTLNuxc} 
&$11$-foam: &  \href{https://youtu.be/pMCCIQzVEqk}{youtu.be/pMCCIQzVEqk} \\
$12$-foam: &  \href{https://youtu.be/-cHEncU2a7o}{youtu.be/-cHEncU2a7o} 
&$13$-foam: &  \href{https://youtu.be/bMBFzjJ-3wY}{youtu.be/bMBFzjJ-3wY} \\
$14$-foam: &  \href{https://youtu.be/Rj9VfPd9Trc}{youtu.be/Rj9VfPd9Trc} 
&$15$-foam: &  \href{https://youtu.be/atUOXP0FtcA}{youtu.be/atUOXP0FtcA} \\
$16$-foam: &  \href{https://youtu.be/QZRtyG-fOb0}{youtu.be/QZRtyG-fOb0} 
&$17$-foam: &  \href{https://youtu.be/AHjbckdh5EY}{youtu.be/AHjbckdh5EY} \\
\hline
\end{tabular}
\caption{Stationary equal-volume $n$-foams for $n=2,\ldots, 17  $ with smallest computed total surface area and links to videos  illustrating the foam structure. See Section~\ref{s:3dStationary}.}
\label{f:eq-volume}
\end{center}
\end{figure}

\begin{figure}[ht]
\begin{center}
\includegraphics[width=.18\textwidth,clip,trim= 6cm 4cm 5cm 4cm]{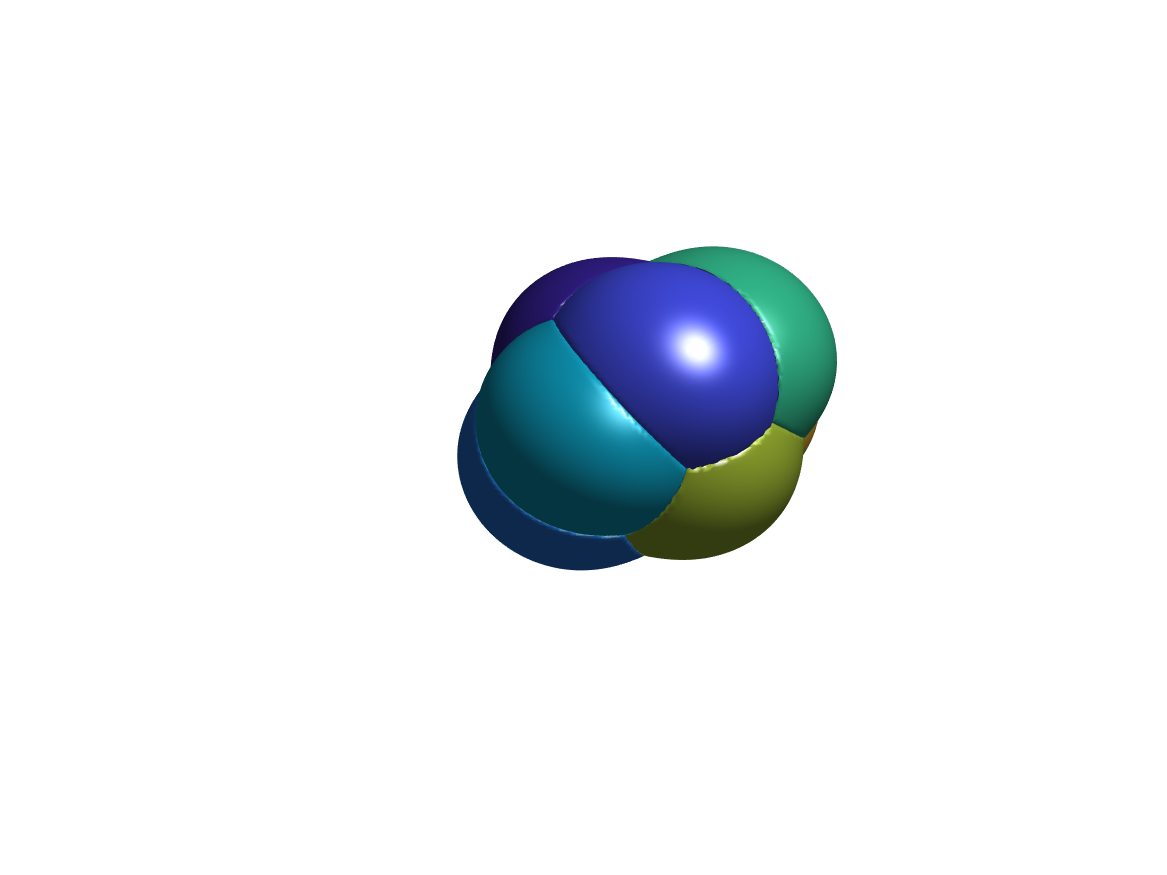}
\includegraphics[width=.18\textwidth,clip,trim= 5cm 3cm 5cm 2cm]{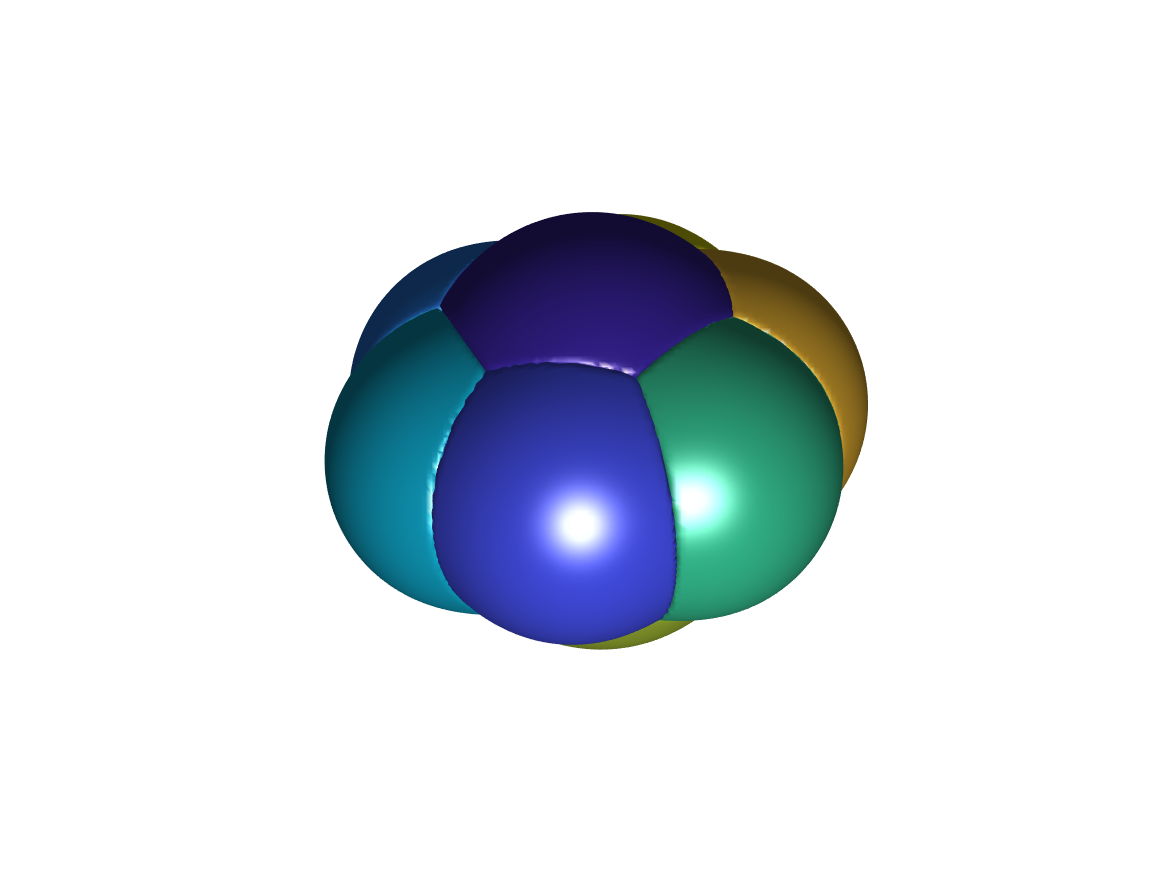}
\includegraphics[width=.18\textwidth,clip,trim= 5cm 3cm 5cm 2cm]{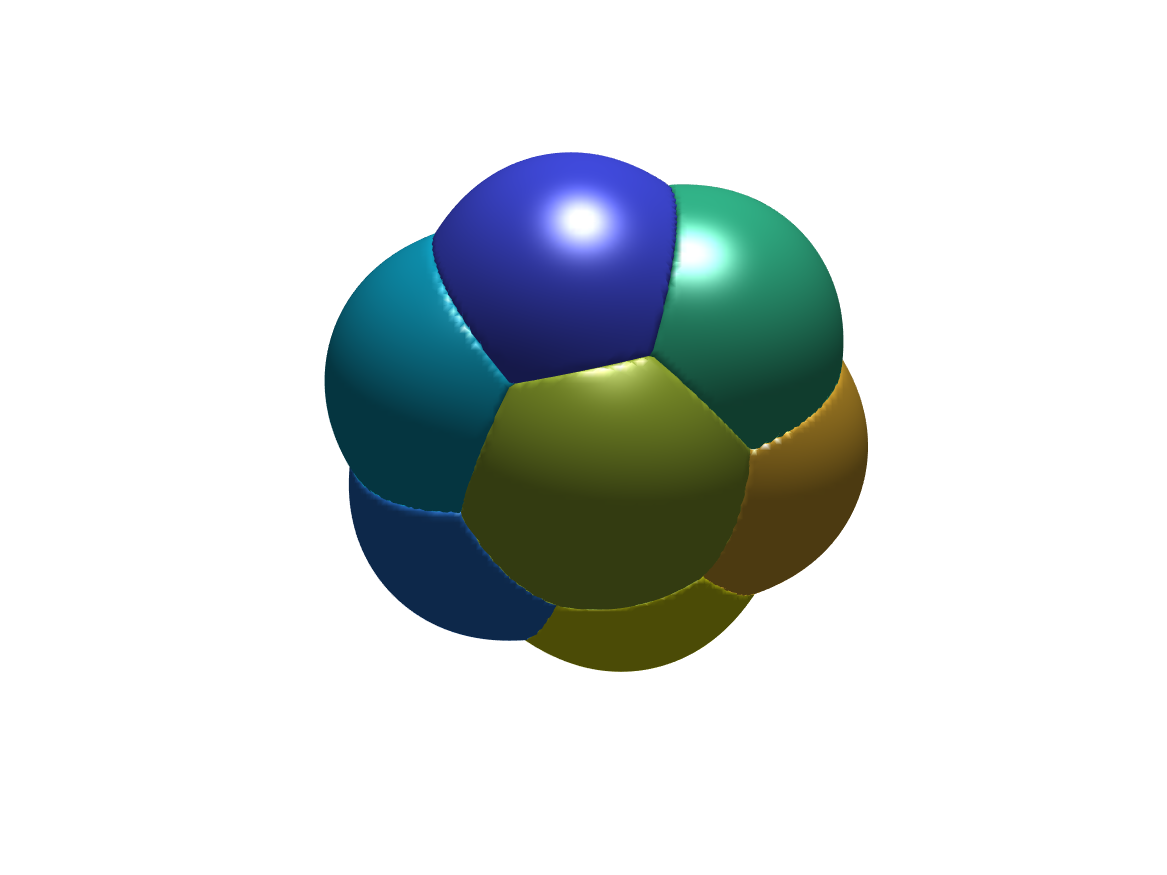}
\includegraphics[width=.18\textwidth,clip,trim= 5cm 3cm 5cm 2cm]{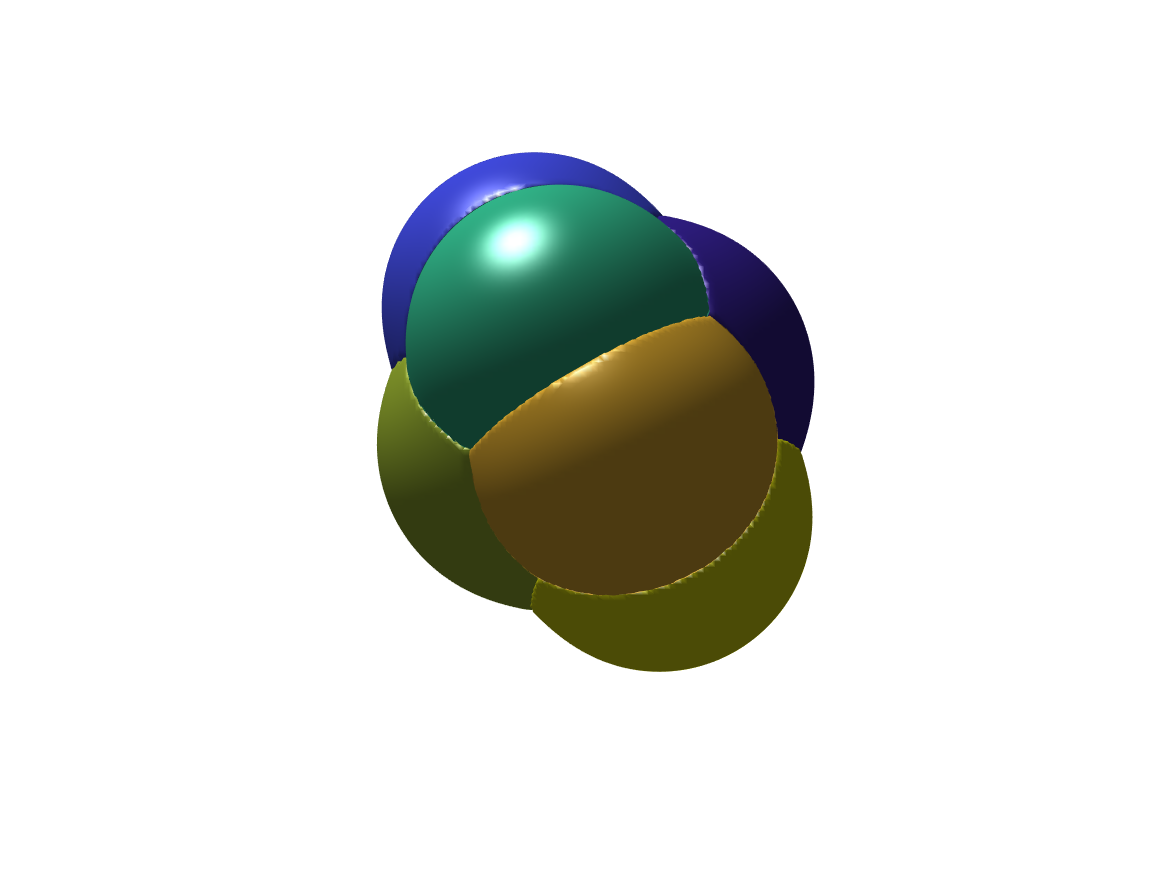}
\includegraphics[width=.18\textwidth,clip,trim= 6cm 4cm 5cm 3cm]{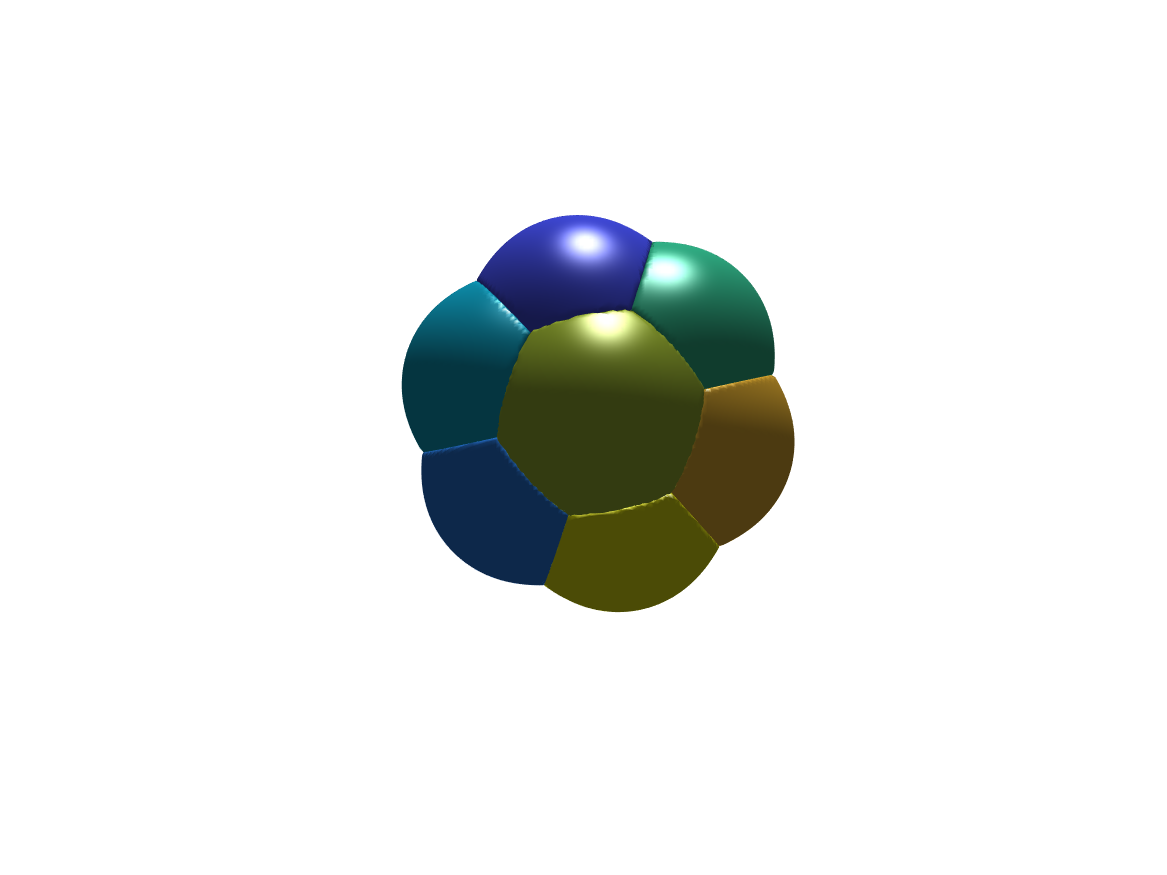}
\caption{The left panel shows another stationary equal-area $8$-foam with larger total surface area than the $8$-foam in Figure~\ref{f:eq-volume}. The middle three panels show $xy$-, $xz$-, and $yz$-views of the $8$-foam. The right panel shows another view showing the hexagonal shaped bubble on the top. A  corresponding video can be found here:  \href{https://youtu.be/4_uAeq19qJY}{youtu.be/4\_uAeq19qJY}.   See Section~\ref{s:3dStationary}. } \label{f:83d}
\end{center}
\end{figure}

\begin{figure}[ht]
\begin{center}
\includegraphics[width=.22\textwidth,clip,trim= 2.5cm 1cm 2.5cm 1cm]{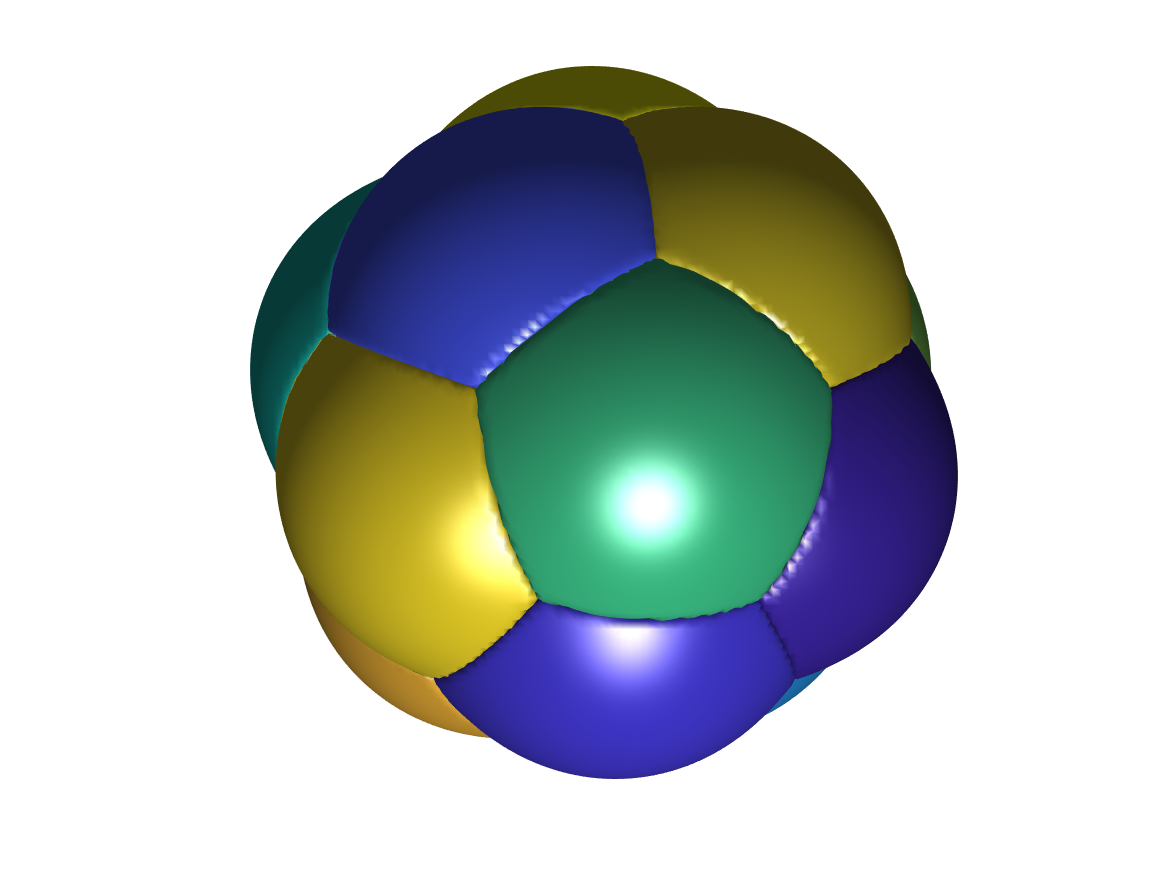}
\includegraphics[width=.22\textwidth,clip,trim= 2.5cm 1cm 2.5cm 1cm]{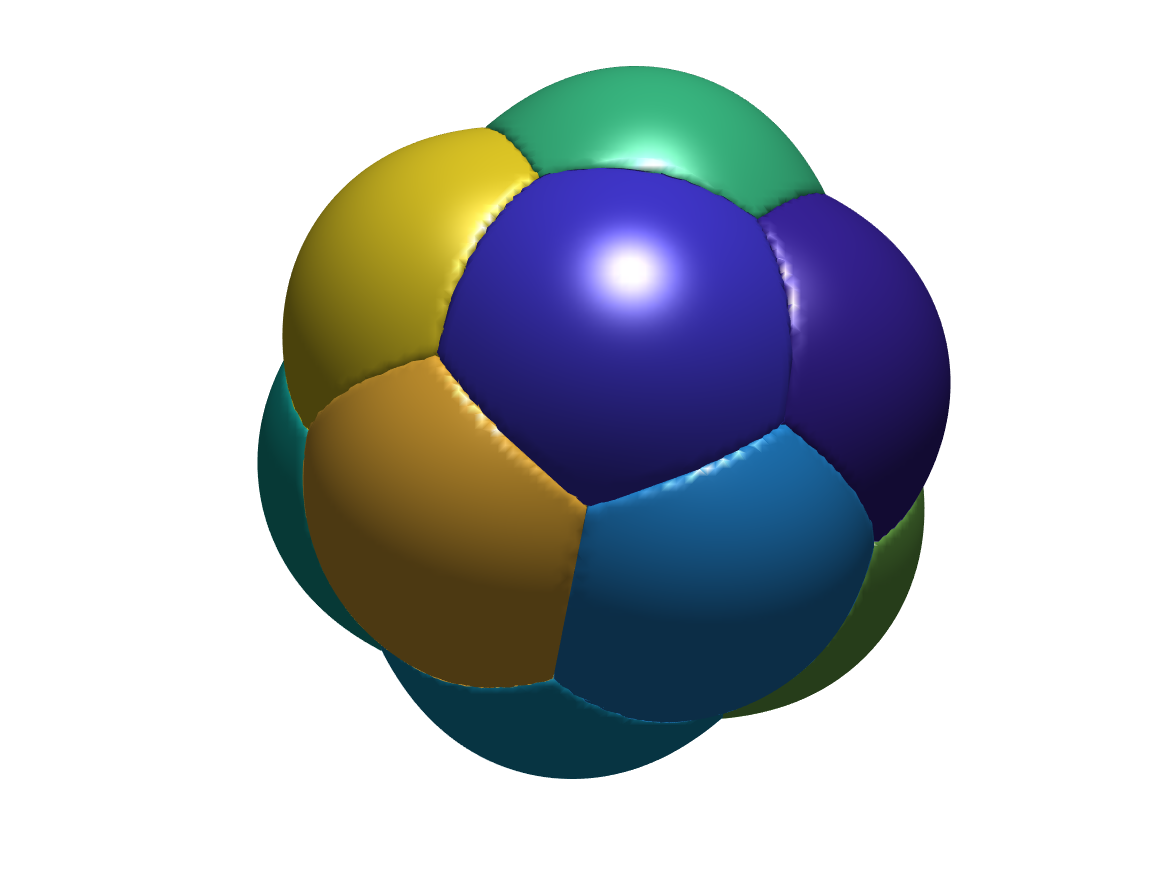}
\includegraphics[width=.22\textwidth,clip,trim= 2.5cm 1cm 2.5cm 1cm]{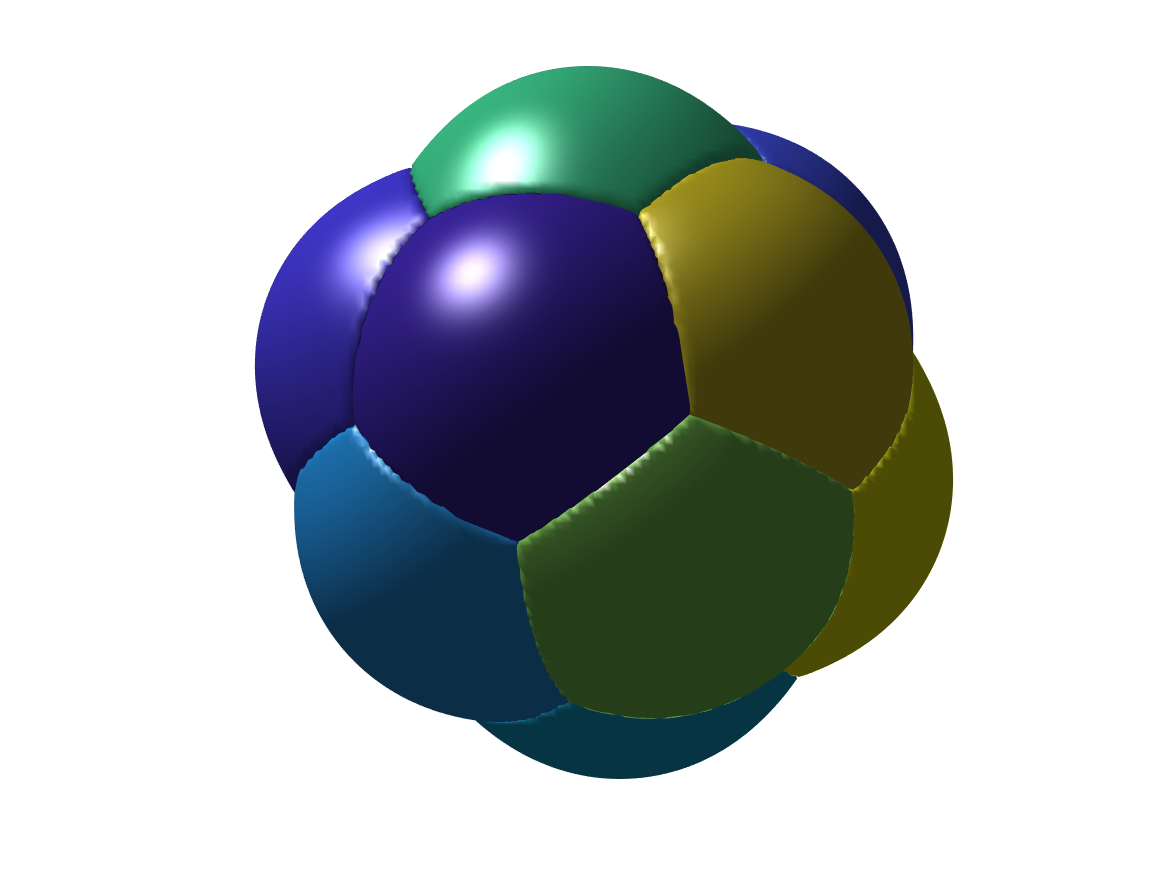}
\includegraphics[width=.22\textwidth,clip,trim= 2.5cm 1cm 2.5cm 1cm]{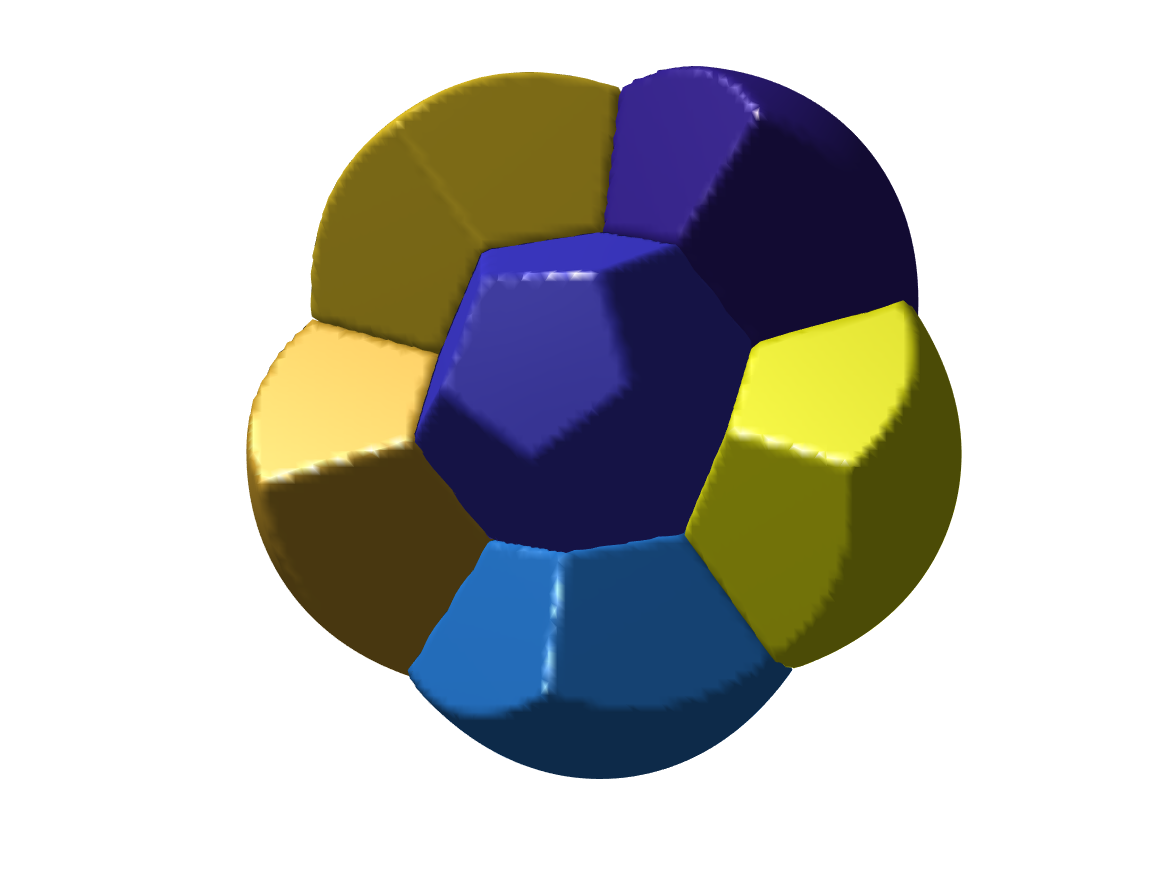}
\caption{The first three panels show $xy$-, $xz$-, and $yz$-views of the $13$-foam in Figure~\ref{f:eq-volume}. The right panel shows a dissection of this foam, exposing the interior bubble, which is a regular dodecahedron. See Section~\ref{s:3dStationary}.} \label{f:133d}
\end{center}
\end{figure}

\begin{figure}[ht]
\begin{center}
\includegraphics[width=.22\textwidth,clip,trim= 5cm 3cm 5cm 2cm]{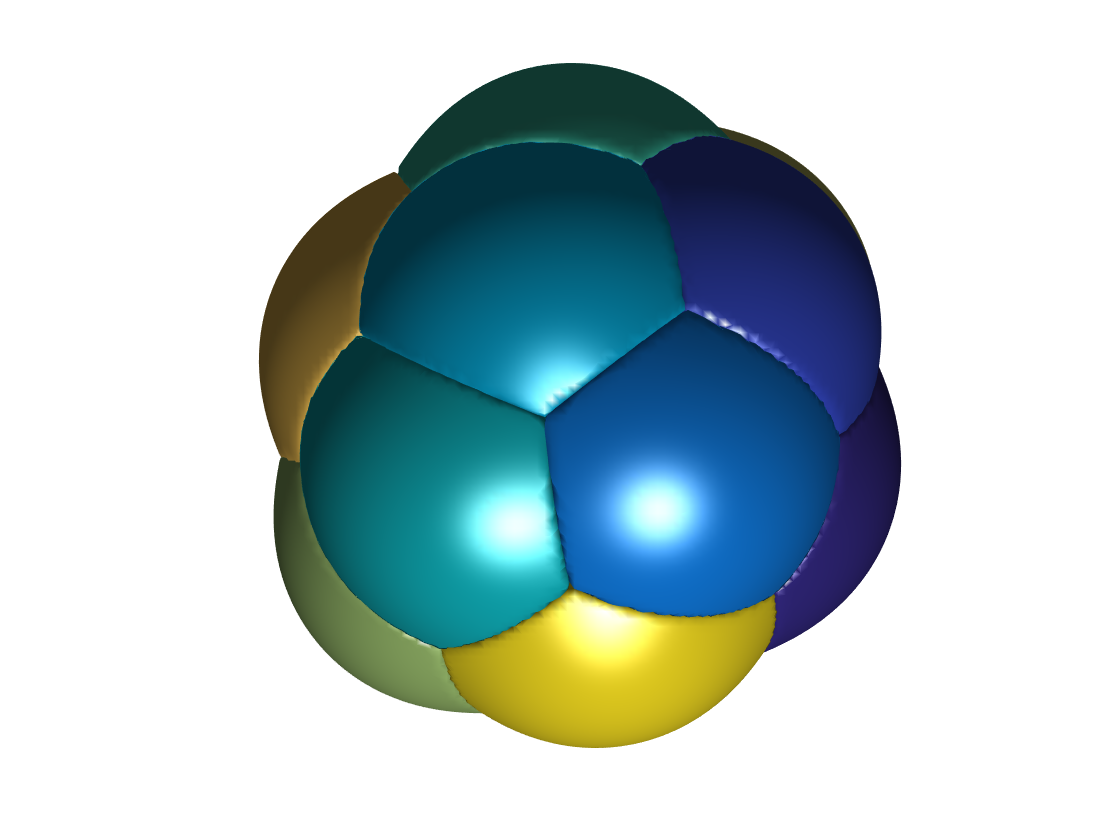}
\includegraphics[width=.22\textwidth,clip,trim= 5cm 3cm 5cm 2cm]{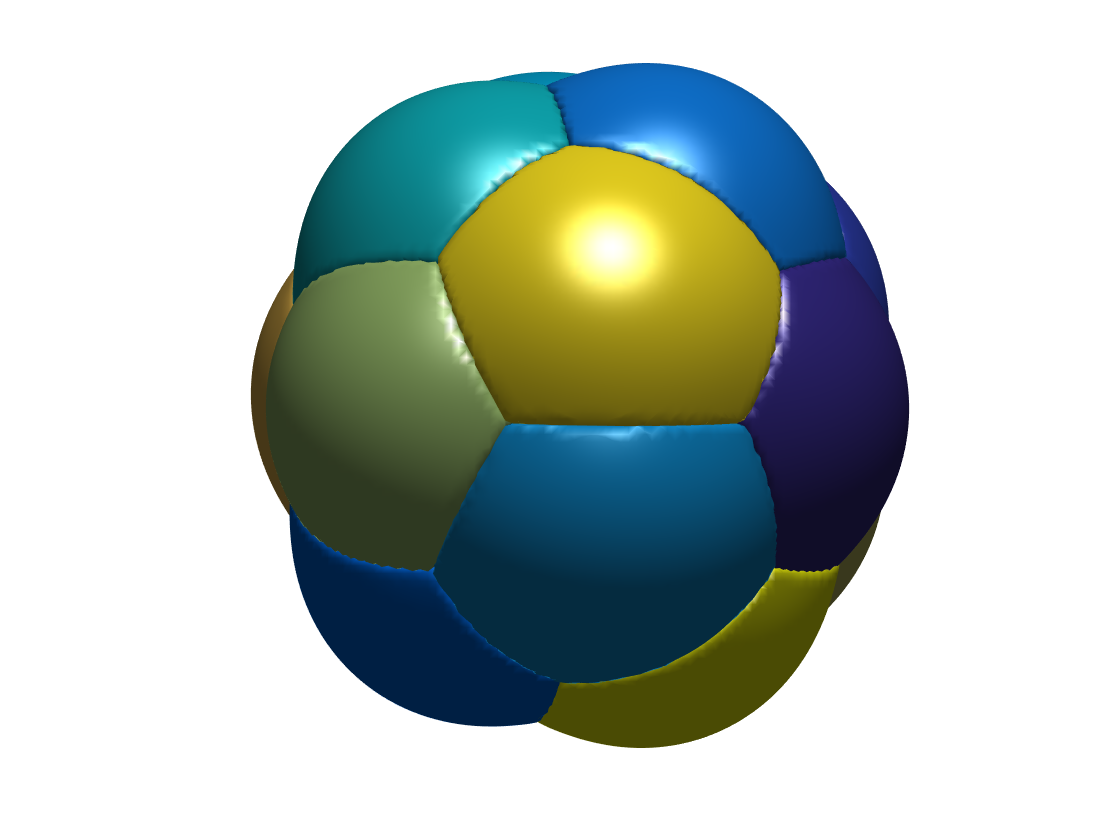}
\includegraphics[width=.22\textwidth,clip,trim= 5cm 3cm 5cm 2cm]{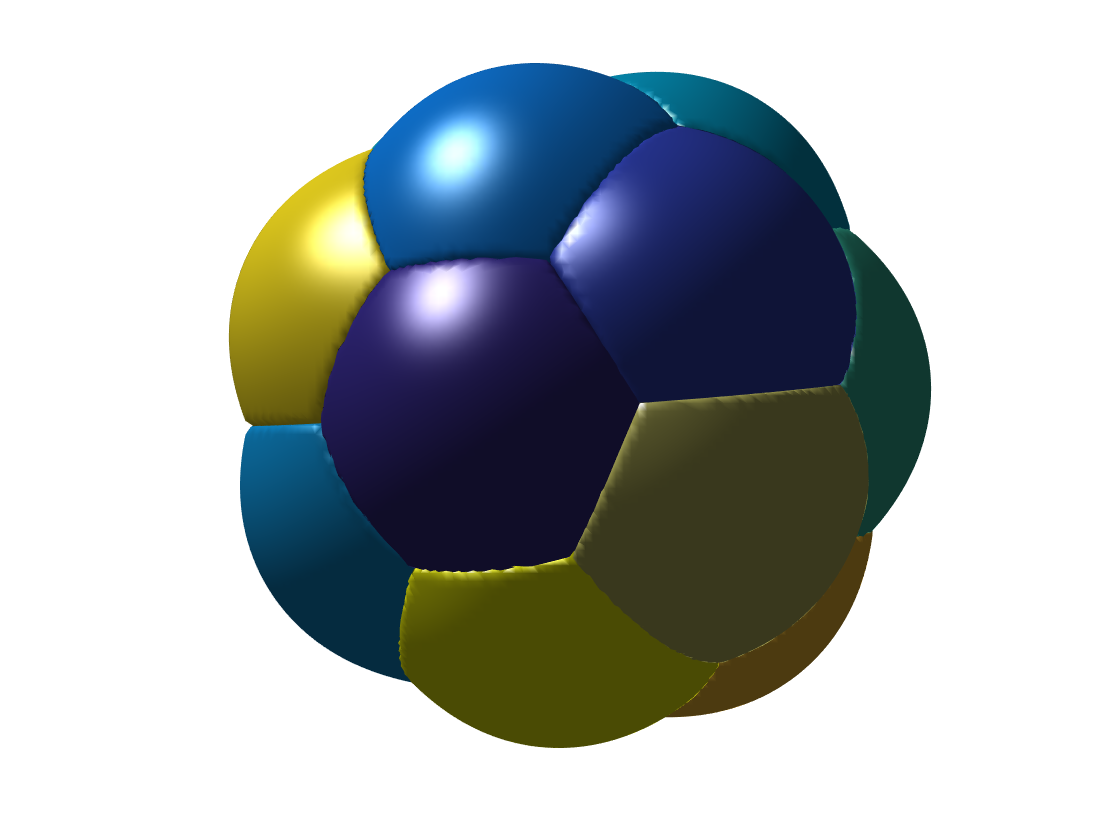}
\includegraphics[width=.22\textwidth,clip,trim= 5cm 3cm 5cm 2cm]{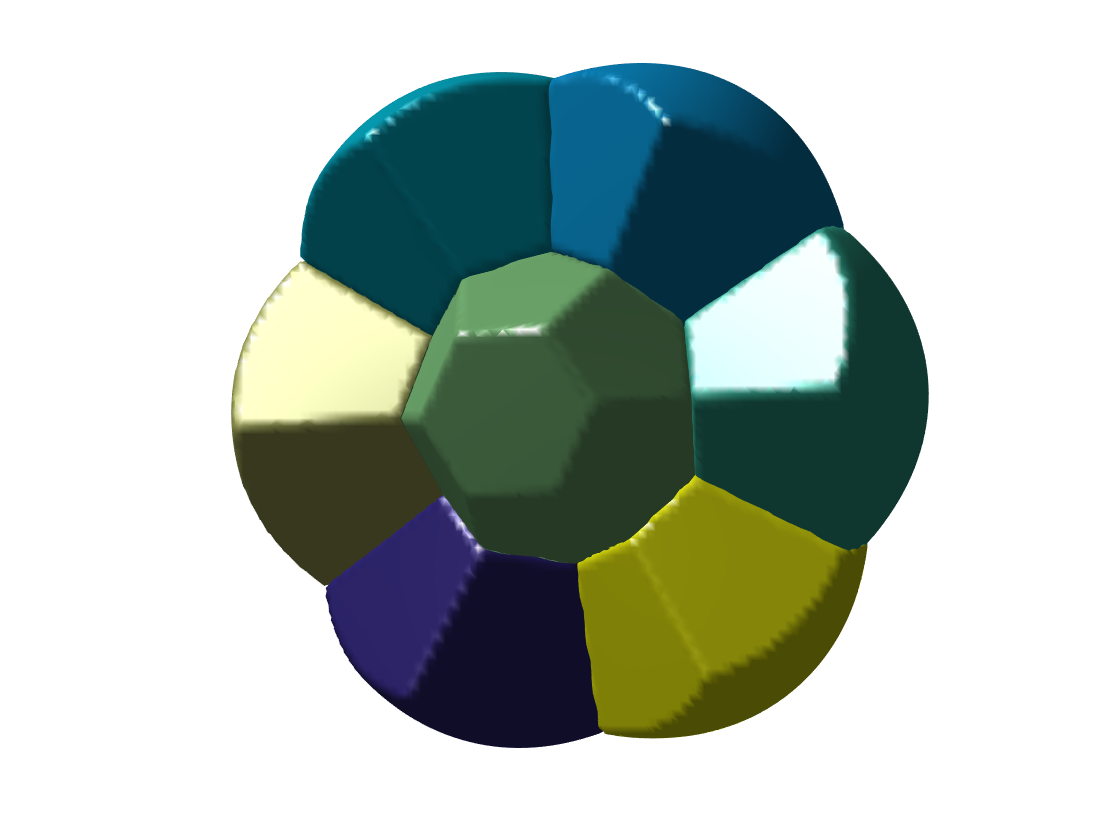}
\caption{The first three panels show $xy$-, $xz$-, and $yz$-views of the $15$-foam in Figure~\ref{f:eq-volume}. The right panel shows a dissection of this foam, exposing the interior bubble, which is similar to the Weaire--Phelan structure. See Section~\ref{s:3dStationary}.} \label{f:153d}
\end{center}
\end{figure}

In Figure~\ref{f:eq-volume}, we plot the three-dimensional  $n$-foams with smallest total surface area found for $n=2,\ldots,17$.  We make the following observations. 
\begin{enumerate}
\item In all cases, Plateau's necessary conditions for optimality, discussed in Section~\ref{s:b3d}, are satisfied. 

\item For $n=2$ and $n=3$, we obtain the expected double and triple-bubble configurations. 

\item For $n=4$, the centers of the bubbles form a tetrahedron. 

\item For $n=5,6,7$, the $n$-foams consist of two vertically-stacked bubbles with $n-2$ bubbles arranged with centers in a regular polygon. 

\item For $n=8$, we repeated the experiment with random initial conditions $100$ times. In $99$ of the experiments, we obtained the $8$-foam as shown in Figure~\ref{f:eq-volume}. In one of the $100$ experiments, we obtained another candidate foam which consists of two vertically-stacked bubbles with $6$ bubbles arranged with centers in a regular hexagon as shown in Figure~\ref{f:83d}. 
The computed total surface area of the configuration in Figure~\ref{f:83d} is $\approx 3.8\%$ higher than the stationary $8$-foam in Figure~\ref{f:eq-volume}. 
It is interesting that the algorithm converges to this local minimizer so infrequently, so the basin of attraction for this local minimum is small. 

\item For $n$-foams with $n\leq 11$, there are no interior bubbles and for $n$-foams with $n\geq 12$, there appears to be at least one interior bubble. 

\item The stationary $13$-foam is very regular and composed of one interior bubble and $12$ bubbles that are on the boundary. In Figure~\ref{f:133d}, we plot $xy$-, $xz$-, and $yz$-views of the $13$-foam and a partial plot of the foam showing the interior bubble. Interestingly, the interior bubble is very similar to a regular dodecahedron. We note that, in a regular dodecahedron, the angle between each two faces is $\approx117^\circ$; we expect the surface of the interior bubble to be slightly curved (non-flat). 

\item The $15$-foam candidate is also very regular and is composed of one interior bubble and $14$ bubbles on the boundary. In Figure~\ref{f:153d}, we plot $xy$-, $xz$-, and $yz$-views of the $15$-foam and a partial plot of the foam showing the interior bubble. The interior bubble is very similar to the truncated hexagonal trapezohedron that appears in the Weaire--Phelan structure. The bubbles on the boundary consist of twelve rounded irregular dodecahedron and two rounded truncated hexagonal trapezohedron. 
\end{enumerate}

\section{Discussion} \label{s:disc}
In this paper, we considered the variational foam model \eqref{e:min}, where the goal is to minimize the total surface area of a collection of bubbles subject to the constraint that the volume of each bubble is prescribed. 
Sharp interface methods together with an approximation of the interfacial surface area using heat diffusion leads to \eqref{pro:lin}, which can be efficiently solved using the auction dynamics method developed in \cite{jacobsauction}. 
This computational method was then used to simulate time dynamics of foams in two- and three-dimensions; 
compute stationary states of foams in two- and three-dimensions; 
and study configurational transitions in the quasi-stationary flow where the volume of one of the bubbles is varied and, for each volume, the stationary state is computed. 
The results from these numerical experiments are described and accompanied by many figures and videos. 

The methods considered in this paper could be used to simulate foams where the bubbles have different surface tensions or different surface mobilities using the modifications developed in \cite{WWX2018}.

In Remark~\ref{rem:SmallVol}, we observed that for small bubbles, a small time step $\tau$ must be used and consequently a fine mesh. Also, the computational cost for this algorithm increases  with the number of bubbles. Finding ways to extend this method to small bubbles and large number of bubbles is challenging and beyond the  scope of this paper. 

One question that we find intriguing is: for fixed $k\in \mathbb N$, how many bubbles in an equal-area stationary foam are needed before there are $k$ in the interior? 
In two-dimensions, we observe that 6 bubbles are needed for one interior bubble, 9 are needed for two, 11 are needed for three, etc\ldots. 
In three-dimensions, 12 bubbles are needed for one interior bubble. 
Numerical evidence suggests that more than 20 bubbles are needed before two interior bubbles appear. 

We hope that the numerical experiments conducted in this paper and further experiments using the methods developed can provide insights for further rigorous geometric results for this foam model. 


\printbibliography

\end{document}